\documentclass[12pt]{article}

\usepackage{amsmath,amssymb,amsthm,amscd,a4wide,upref}
\usepackage[toc]{appendix}

\usepackage{hyperref}
\hypersetup{colorlinks,citecolor=blue,filecolor=black,linkcolor=blue,urlcolor=blue}
\usepackage{enumerate}
\usepackage{mathtools}

\begin{document}

\newcommand{\ad}{{\rm ad}}
\newcommand{\cri}{{\rm cri}}
\newcommand{\End}{{\rm{End}\ts}}
\newcommand{\Rep}{{\rm{Rep}\ts}}
\newcommand{\Hom}{{\rm{Hom}}}
\newcommand{\Mat}{{\rm{Mat}}}
\newcommand{\ch}{{\rm{ch}\ts}}
\newcommand{\chara}{{\rm{char}\ts}}
\newcommand{\diag}{{\rm diag}}
\newcommand{\non}{\nonumber}
\newcommand{\wt}{\widetilde}
\newcommand{\wh}{\widehat}
\newcommand{\ot}{\otimes}
\newcommand{\la}{\lambda}
\newcommand{\La}{\Lambda}
\newcommand{\De}{\Delta}
\newcommand{\al}{\alpha}
\newcommand{\be}{\beta}
\newcommand{\ga}{\gamma}
\newcommand{\Ga}{\Gamma}
\newcommand{\ep}{\epsilon}
\newcommand{\ka}{\kappa}
\newcommand{\vk}{\varkappa}
\newcommand{\si}{\sigma}
\newcommand{\vs}{\varsigma}
\newcommand{\vp}{\varphi}
\newcommand{\de}{\delta}
\newcommand{\ze}{\zeta}
\newcommand{\om}{\omega}
\newcommand{\Om}{\Omega}
\newcommand{\ee}{\epsilon^{}}
\newcommand{\su}{s^{}}
\newcommand{\hra}{\hookrightarrow}
\newcommand{\ve}{\varepsilon}
\newcommand{\ts}{\,}
\newcommand{\pr}{^{\tss\prime}}
\newcommand{\vac}{\mathbf{1}}
\newcommand{\di}{\partial}
\newcommand{\qin}{q^{-1}}
\newcommand{\tss}{\hspace{1pt}}
\newcommand{\Sr}{ {\rm S}}
\newcommand{\U}{ {\rm U}}
\newcommand{\BL}{ {\overline L}}
\newcommand{\BE}{ {\overline E}}
\newcommand{\BP}{ {\overline P}}
\newcommand{\AAb}{\mathbb{A}\tss}
\newcommand{\CC}{\mathbb{C}\tss}
\newcommand{\KK}{\mathbb{K}\tss}
\newcommand{\QQ}{\mathbb{Q}\tss}
\newcommand{\SSb}{\mathbb{S}\tss}
\newcommand{\TT}{\mathbb{T}\tss}
\newcommand{\ZZ}{\mathbb{Z}\tss}
\newcommand{\DY}{ {\rm DY}}
\newcommand{\X}{ {\rm X}}
\newcommand{\Y}{ {\rm Y}}
\newcommand{\Z}{{\rm Z}}
\newcommand{\Ac}{\mathcal{A}}
\newcommand{\Lc}{\mathcal{L}}
\newcommand{\Mc}{\mathcal{M}}
\newcommand{\Pc}{\mathcal{P}}
\newcommand{\Qc}{\mathcal{Q}}
\newcommand{\Rc}{\mathcal{R}}
\newcommand{\Sc}{\mathcal{S}}
\newcommand{\Tc}{\mathcal{T}}
\newcommand{\Bc}{\mathcal{B}}
\newcommand{\Ec}{\mathcal{E}}
\newcommand{\Fc}{\mathcal{F}}
\newcommand{\Gc}{\mathcal{G}}
\newcommand{\Hc}{\mathcal{H}}
\newcommand{\Uc}{\mathcal{U}}
\newcommand{\Vc}{\mathcal{V}}
\newcommand{\Wc}{\mathcal{W}}
\newcommand{\Yc}{\mathcal{Y}}
\newcommand{\Ar}{{\rm A}}
\newcommand{\Br}{{\rm B}}
\newcommand{\Ir}{{\rm I}}
\newcommand{\Fr}{{\rm F}}
\newcommand{\Jr}{{\rm J}}
\newcommand{\Or}{{\rm O}}
\newcommand{\GL}{{\rm GL}}
\newcommand{\Spr}{{\rm Sp}}
\newcommand{\Rr}{{\rm R}}
\newcommand{\Zr}{{\rm Z}}
\newcommand{\ZX}{{\rm ZX}}
\newcommand{\gl}{\mathfrak{gl}}
\newcommand{\middd}{{\rm mid}}
\newcommand{\ev}{{\rm ev}}
\newcommand{\Pf}{{\rm Pf}}
\newcommand{\Norm}{{\rm Norm\tss}}
\newcommand{\oa}{\mathfrak{o}}
\newcommand{\spa}{\mathfrak{sp}}
\newcommand{\osp}{\mathfrak{osp}}
\newcommand{\f}{\mathfrak{f}}
\newcommand{\g}{\mathfrak{g}}
\newcommand{\h}{\mathfrak h}
\newcommand{\n}{\mathfrak n}
\newcommand{\z}{\mathfrak{z}}
\newcommand{\Zgot}{\mathfrak{Z}}
\newcommand{\p}{\mathfrak{p}}
\newcommand{\sll}{\mathfrak{sl}}
\newcommand{\agot}{\mathfrak{a}}
\newcommand{\qdet}{ {\rm qdet}\ts}
\newcommand{\Ber}{ {\rm Ber}\ts}
\newcommand{\HC}{ {\mathcal HC}}
\newcommand{\cdet}{{\rm cdet}}
\newcommand{\rdet}{{\rm rdet}}
\newcommand{\tr}{ {\rm tr}}
\newcommand{\gr}{ {\rm gr}\ts}
\newcommand{\str}{ {\rm str}}
\newcommand{\loc}{{\rm loc}}
\newcommand{\Gr}{{\rm G}}
\newcommand{\sgn}{ {\rm sgn}\ts}
\newcommand{\sign}{{\rm sgn}}
\newcommand{\ba}{\bar{a}}
\newcommand{\bb}{\bar{b}}
\newcommand{\eb}{\bar{e}}
\newcommand{\bi}{\bar{\imath}}
\newcommand{\bj}{\bar{\jmath}}
\newcommand{\bk}{\bar{k}}
\newcommand{\bl}{\bar{l}}
\newcommand{\hb}{\mathbf{h}}
\newcommand{\Sym}{\mathfrak S}
\newcommand{\fand}{\quad\text{and}\quad}
\newcommand{\Fand}{\qquad\text{and}\qquad}
\newcommand{\For}{\qquad\text{or}\qquad}
\newcommand{\OR}{\qquad\text{or}\qquad}
\newcommand{\grpr}{{\rm gr}^{\tss\prime}\ts}
\newcommand{\degpr}{{\rm deg}^{\tss\prime}\tss}

\numberwithin{equation}{section}

\newtheorem{thm}{Theorem}[section]
\newtheorem{lem}[thm]{Lemma}
\newtheorem{prop}[thm]{Proposition}
\newtheorem{cor}[thm]{Corollary}
\newtheorem{conj}[thm]{Conjecture}
\newtheorem*{mthm}{Main Theorem}
\newtheorem*{mthma}{Theorem A}
\newtheorem*{mthmb}{Theorem B}
\newtheorem*{mthmc}{Theorem C}
\newtheorem*{mthmd}{Theorem D}

\theoremstyle{definition}
\newtheorem{defin}[thm]{Definition}

\theoremstyle{remark}
\newtheorem{remark}[thm]{Remark}
\newtheorem{example}[thm]{Example}

\newcommand{\bth}{\begin{thm}}
\renewcommand{\eth}{\end{thm}}
\newcommand{\bpr}{\begin{prop}}
\newcommand{\epr}{\end{prop}}
\newcommand{\ble}{\begin{lem}}
\newcommand{\ele}{\end{lem}}
\newcommand{\bco}{\begin{cor}}
\newcommand{\eco}{\end{cor}}
\newcommand{\bde}{\begin{defin}}
\newcommand{\ede}{\end{defin}}
\newcommand{\bex}{\begin{example}}
\newcommand{\eex}{\end{example}}
\newcommand{\bre}{\begin{remark}}
\newcommand{\ere}{\end{remark}}
\newcommand{\bcj}{\begin{conj}}
\newcommand{\ecj}{\end{conj}}

\newcommand{\bal}{\begin{aligned}}
\newcommand{\eal}{\end{aligned}}
\newcommand{\beq}{\begin{equation}}
\newcommand{\eeq}{\end{equation}}
\newcommand{\ben}{\begin{equation*}}
\newcommand{\een}{\end{equation*}}

\newcommand{\bpf}{\begin{proof}}
\newcommand{\epf}{\end{proof}}

\def\beql#1{\begin{equation}\label{#1}}

\title{\Large\bf Isomorphism between the $R$-matrix and Drinfeld presentations
of Yangian in types $B,C$ and $D$}

\author{{Naihuan Jing,\quad Ming Liu\footnote{Corresponding author} \quad and\quad Alexander Molev}}

\date{} % Start January 2017
\maketitle

\vspace{5 mm}

\begin{abstract}
It is well-known that the Gauss decomposition of the generator matrix
in the $R$-matrix presentation of the Yangian in type $A$ yields generators of
its Drinfeld presentation. Defining relations between these generators
are known in an explicit form thus providing an isomorphism between
the presentations. It has been an open problem since the pioneering work
of Drinfeld to extend this result to the remaining types.
We give a solution for the classical types $B$, $C$ and $D$
by constructing an explicit isomorphism between
the $R$-matrix and Drinfeld presentations of the Yangian.
It is based on an embedding theorem which allows us to consider
the Yangian of rank $n-1$ as a subalgebra of the Yangian of rank $n$ of the same type.

%\medskip

%Mathematics Subject Classification 2010: 17B37, 17B69

\end{abstract}

\vspace{5 mm}
%%%
%{\it Key words:}
%%%
%\footnote{Corresponding author: mamliu@scut.edu.cn}

%\tableofcontents
%

\section{Introduction}
\label{sec:int}

According to the original definition of Drinfeld~\cite{d:ha}, the {\em Yangian}
associated to a simple Lie algebra $\g$ is a Hopf algebra with
a finite set of generators. Another presentation of the Yangian
was given by him in \cite{d:nr} and it is known as the {\em new realization} or
{\em Drinfeld presentation}; see also book by Chari and Pressly~\cite[Chapter~12]{cp:gq}
for an exposition.
The Hopf algebra which coincides with the Yangian in type $A$
was considered previously in the work of
Faddeev and the St.~Petersburg (Leningrad) school; see
the expository paper \cite{ks:qs} by Kulish and Sklyanin.
The defining relations of this algebra are written
in the form of a single $RTT$ relation involving the Yang $R$-matrix $R$.
An explicit isomorphism between the $R$-matrix and Drinfeld presentations
of the Yangian in type $A$ is constructed with the use of the Gauss decomposition
of the generator matrix $T(u)$. Complete proofs were given
by Brundan and Kleshchev~\cite{bk:pp}; see also \cite[Section~3.1]{m:yc} for
an exposition.

At least for the classical types,
the $R$-matrix presentation is convenient for describing the coproduct of the Yangian
and allows one to develop tensor techniques to investigate its algebraic structure
and representations;
cf. Arnaudon {\em et al.}~\cite{aacfr:rp}, \cite{amr:rp},
Guay {\em et al.}~\cite{gr:ty}, \cite{grw:rt} (types $B,C$ and $D$)
and \cite[Chapter~1]{m:yc} (type~$A$).
On the other hand, finite-dimensional irreducible
representations of the Yangian associated with any simple Lie algebra $\g$
are classified uniformly in terms of its
Drinfeld presentation. An explicit isomorphism between the presentations
is therefore important for bringing together the two approaches and enhancing
algebraic tools for understanding the Yangian and its representations.
Our main result is a construction of such an isomorphism
in the remaining classical types $B$, $C$ and $D$ which thus solves the open problem
going back to Drinfeld's work \cite{d:nr}.

To explain our construction, suppose that the simple Lie
algebra $\g$ is
associated with the Cartan matrix $A=[a_{ij}]_{i,j=1}^n$. Let $\al_1,\dots,\al_n$
be the corresponding simple roots (normalized
as in \eqref{ali} and \eqref{aln} below
for types $B_n$, $C_n$ and $D_n$).
In accordance to \cite{d:nr},
the {\em Drinfeld Yangian} $\Y^{D}(\g)$ is generated by
elements $\kappa^{}_{i\tss r}$, $\xi_{i\tss r}^{+}$ and
$\xi_{i\tss r}^{-}$ with $i=1,\dots,n$ and $r=0,1,\dots$
subject to the defining relations
\begin{align}
\non%\label{kapkap}
[\kappa^{}_{i\tss r},\kappa^{}_{j\tss s}]&=0,\\[0.3em]
\non%\label{xipxim}
[\xi_{i\tss r}^{+},\xi_{j\tss s}^{-}]&=\de_{ij}\ts\kappa^{}_{i\ts r+s},\\[0.3em]
\non%\label{kapoxi}
[\kappa^{}_{i\tss 0},\xi_{j\tss s}^{\pm}]&=\pm\ts (\al_i,\al_j)\ts\xi_{j\tss s}^{\pm},
\\
\non%\label{kapxi}
[\kappa^{}_{i\ts r+1},\xi_{j\tss s}^{\pm}]-[\kappa^{}_{i\tss r},\xi_{j\ts s+1}^{\pm}]&=
\pm\ts\frac{(\al_i,\al_j)}{2}\big(\kappa^{}_{i\tss r}\ts \xi_{j\tss s}^{\pm}
+\xi_{j\tss s}^{\pm}\ts \kappa^{}_{i\tss r}\big),\\
\non%\label{xixi}
[\xi_{i\ts r+1}^{\pm},\xi_{j\tss s}^{\pm}]-[\xi_{i\tss r}^{\pm},\xi_{j\ts s+1}^{\pm}]&=
\pm\ts\frac{(\al_i,\al_j)}{2}\big(\xi_{i\tss r}^{\pm}\ts \xi_{j\tss s}^{\pm}
+\xi_{j\tss s}^{\pm}\ts \xi_{i\tss r}^{\pm}\big),\\[0.3em]
\label{xisym}
\sum_{p\in\Sym_m}
\,
[\xi_{i\tss r_{p(1)}}^{\pm},[\xi_{i\tss r_{p(2)}}^{\pm},\dots
&[\xi_{i\tss r_{p(m)}}^{\pm},\xi_{j\tss s}^{\pm}]\dots ]]=0,
\end{align}
where the last relation holds for
all $i\ne j$, and we denoted\/ $m=1-a_{ij}$.

If $\g=\g_N$
is the orthogonal Lie algebra $\oa_N$ (with $N=2n$ or $N=2n+1$)
or symplectic Lie algebra $\spa_N$
(with even $N=2n$) then the algebra $\Y^{R}(\g_N)$
(the {\em Yangian in the $R$-matrix or $RTT$ presentation})
can be defined with the use of the rational $R$-matrix
first discovered in \cite{zz:rf}. The defining relations take the form
of the $RTT$ {\em relation}
\beql{RTTint}
R(u-v)\ts T_1(u)\ts T_2(v)=T_2(v)\ts T_1(u)\ts R(u-v)
\eeq
together with the {\em unitarity condition}
\beql{unitaryint}
T^{\tss\prime}(u+\ka)\ts T(u)=1,
\eeq
with the notation explained below in Section~\ref{sec:nd}. Here
$T(u)$ is a square matrix of size $N$ whose $(i,j)$ entry is the formal series
\ben
t_{ij}(u)=\de_{ij}+\sum_{r=1}^{\infty}t_{ij}^{(r)}\ts u^{-r}
\een
so that the algebra $\Y^{R}(\g_N)$ is generated by all coefficients
$t_{ij}^{(r)}$ subject to the defining relations
\eqref{RTTint} and \eqref{unitaryint}.
Apply the Gauss decomposition
to the matrix $T(u)$,
\beql{gd}
T(u)=F(u)\ts H(u)\ts E(u),
\eeq
where $F(u)$, $H(u)$ and $E(u)$ are uniquely determined matrices of the form
\ben
F(u)=\begin{bmatrix}
1&0&\dots&0\ts\\
f_{21}(u)&1&\dots&0\\
\vdots&\vdots&\ddots&\vdots\\
f_{N1}(u)&f_{N2}(u)&\dots&1
\end{bmatrix},
\qquad
E(u)=\begin{bmatrix}
\ts1&e_{12}(u)&\dots&e_{1N}(u)\ts\\
\ts0&1&\dots&e_{2N}(u)\\
\vdots&\vdots&\ddots&\vdots\\
0&0&\dots&1
\end{bmatrix},
\een
and $H(u)=\diag\ts\big[h_1(u),\dots,h_N(u)\big]$. Define the series
with coefficients in $\Y^{R}(\g_N)$ by
\ben
\kappa^{}_i(u)=h_i\big(u-(i-1)/2\big)^{-1}\ts h_{i+1}\big(u-(i-1)/2\big)
\een
for $i=1,\dots,n-1$, and
\ben
\kappa^{}_n(u)=\begin{cases}
h_n\big(u-(n-1)/2\big)^{-1}\ts h_{n+1}\big(u-(n-1)/2\big)
\qquad&\text{for}\quad \oa_{2n+1}\\[0.6em]
\tss h_n\big(u-n/2\big)^{-1}\ts h_{n+1}\big(u-n/2\big)
\qquad&\text{for}\quad \spa_{2n}\\[0.6em]
h_{n-1}\big(u-(n-2)/2\big)^{-1}\ts h_{n+1}\big(u-(n-2)/2\big)
\qquad&\text{for}\quad \oa_{2n}.
\end{cases}
\een
Furthermore, set
\ben
\xi_i^+(u)=f_{i+1\ts i}\big(u-(i-1)/2\big),\qquad
\xi_i^-(u)=e_{i\ts i+1}\big(u-(i-1)/2\big)
\een
for $i=1,\dots,n-1$,
\ben
\xi_n^+(u)=\begin{cases}
f_{n+1\ts n}\big(u-(n-1)/2\big)
\qquad&\text{for}\quad \oa_{2n+1}\\[0.3em]
f_{n+1\ts n}\big(u-n/2\big)
\qquad&\text{for}\quad \spa_{2n}\\[0.3em]
f_{n+1\ts n-1}\big(u-(n-2)/2\big)
\qquad&\text{for}\quad \oa_{2n}
\end{cases}
\een
and
\ben
\xi_n^-(u)=\begin{cases}
e_{n\ts n+1}\big(u-(n-1)/2\big)
\qquad&\text{for}\quad \oa_{2n+1}\\[0.3em]
\frac{1}{2}e_{n\ts n+1}\big(u-n/2\big)
\qquad&\text{for}\quad \spa_{2n}\\[0.3em]
e_{n-1\ts n+1}\big(u-(n-2)/2\big)
\qquad&\text{for}\quad \oa_{2n}.
\end{cases}
\een
Introduce elements of $\Y^{R}(\g_N)$ by the respective expansions
into power series in $u^{-1}$,
\beql{kaxicoeff}
\kappa^{}_i(u)=1+\sum_{r=0}^{\infty}\kappa^{}_{i\tss r}\ts u^{-r-1}\Fand
\xi_i^{\pm}(u)=\sum_{r=0}^{\infty}\xi_{i\tss r}^{\pm}\ts u^{-r-1}
\eeq
for $i=1,\dots,n$. Our main result is the following.

\begin{mthm}
The mapping which sends the generators
$\kappa^{}_{i\tss r}$ and $\xi_{i\tss r}^{\pm}$ of $\Y^{D}(\g_N)$ to the elements
of $\Y^{R}(\g_N)$ with the same names defines an isomorphism
$\Y^{D}(\g_N)\cong \Y^{R}(\g_N)$.
\end{mthm}

As was pointed out in \cite{aacfr:rp}, the existence of such an isomorphism
follows from the fact that the classical limits of the Hopf algebras
$\Y^{D}(\g_N)$ and $\Y^{R}(\g_N)$ define the same bialgebra structure
on the Lie algebra of polynomial currents $\g_N[x]$. Therefore,
the Hopf algebras must be isomorphic due to Drinfeld's uniqueness theorem on quantization.
The main obstacle for constructing an isomorphism explicitly in types $B$, $C$ and $D$
is that, unlike type $A$, there is no {\em natural} embedding
of the Yangian of rank $n-1$ into the Yangian of rank $n$ in their $R$-matrix presentations.
To overcome this difficulty, we first prove
an embedding theorem
allowing us to consider
the Yangian $\Y^{R}(\g_{N-2})$ as a subalgebra of $\Y^{R}(\g_N)$.
When restricted to the universal enveloping algebras, this coincides with the natural
embedding $\U(\g_{N-2})\hra \U(\g_N)$.
This theorem effectively reduces the isomorphism problem to the case of rank $2$.

The second ingredient of the proof of the Main Theorem
is another presentation (called {\em minimalistic}) of the Yangian
$\Y^{D}(\g_N)$
which goes back to Levendorski\v{\i}~\cite{l:gd} and was
recently given in a modified form by Guay~{\em et al.}~\cite{gu:co}.
Its use eliminates the need to verify complicated Serre-type relations
in $\Y^{D}(\g)$ for the proof that our map is a homomorphism.

We will mainly work with the {\em extended Yangian} $\X(\g_N)$
which is defined by the relation \eqref{RTTint}
omitting \eqref{unitaryint}. We give a Drinfeld presentation
for $\X(\g_N)$ which we believe is of independent interest.
This presentation given in Theorem~\ref{thm:dp} is analogous to
the one in type $A$; see \cite{bk:pp}.
Furthermore, we describe the center $\ZX(\g_N)$ of the extended Yangian
in its Drinfeld presentation by providing explicit formulas for generators
of $\ZX(\g_N)$ (Theorem~\ref{thm:Center}) which are then used in the proof of the Main Theorem.

We expect that the results of this paper
can be extended to the quantum affine algebras of types $B$, $C$ and $D$
to get corresponding analogues of the Ding--Frenkel isomorphism~\cite{df:it}.

After posting the first version of this paper to the arXiv, we received the
preprint \cite{grw:eb} from the authors where a closely related results were obtained.
In particular, an isomorphism between the original presentation \cite{d:ha} of the Yangian in type
$B,C$ and $D$ and its $R$-matrix presentation is produced. However, the approach and
constructions of \cite{grw:eb} are
different from ours, as they do not rely on the Gauss decomposition.

\medskip

%We acknowledge the support of the School of Mathematical Sciences
%at the South China University of Technology
%and the Australian Research Council.

\section{Notation and definitions}
\label{sec:nd}

Define the orthogonal Lie algebras $\oa_N$ with $N=2n+1$ and $N=2n$
(corresponding to types $B$ and $D$, respectively)
and symplectic Lie algebra $\spa_N$ with $N=2n$ (of type $C$)
as subalgebras of $\gl_N$ spanned by all elements $F_{i\tss j}$,
\beql{fij}
F_{i\tss j}=E_{i\tss j}-E_{j\pr i\pr}\Fand F_{i\tss j}
=E_{i\tss j}-\ve_i\ts\ve_j\ts E_{j\pr i\pr},
\eeq
respectively, for $\oa_N$ and $\spa_N$, where the $E_{i\tss j}$ denote
the standard basis elements of $\gl_N$.
Here we use the notation $i\pr=N-i+1$, and
in the symplectic case we set
$\ve_i=1$ for $i=1,\dots,n$ and
$\ve_i=-1$ for $i=n+1,\dots,2n$.
To consider the three cases $B$, $C$ and $D$
simultaneously we will use
the notation $\g_N$ for any of the Lie algebras $\oa_N$
or $\spa_N$.

To introduce the $R$-matrix presentation of the Yangian, we will need a standard
tensor notation. By taking the canonical basis $e_1,\dots,e_N$ of $\CC^N$, we
will identify the endomorphism algebra
$\End\CC^N$
with the algebra of $N\times N$ matrices.
The matrix units
$e_{ij}$ with $i,j\in\{1,\dots,N\}$ form a basis of $\End\CC^N$.
We will work with tensor product algebras of the form
\beql{tenprkea}
\End(\CC^{N})^{\ot m}\ot\Ac=\underbrace{\End\CC^{N}\ot\dots\ot\End\CC^{N}}_m{}\ot\Ac,
\eeq
where $\Ac$ is a unital associative algebra.
For any element
\beql{mata}
X=\sum_{i,j=1}^N e_{ij}\ot X_{ij}\in \End\CC^N\ot \Ac
\eeq
and any $a\in\{1,\dots,m\}$ we will denote by $X_a$ the element \eqref{mata}
associated with the $a$-th copy of $\End\CC^{N}$ so that
\beql{xa}
X_a=\sum_{i,j=1}^N 1^{\ot(a-1)}\ot e_{ij}\ot 1^{\ot(m-a)}\ot X_{ij}
\in \End(\CC^{N})^{\ot m}\ot\Ac,
\eeq
where $1$
is the identity endomorphism. Moreover, given any element
\ben
C=\sum_{i,j,k,l=1}^N c^{}_{ijkl}\ts e_{ij}\ot e_{kl}\in
\End \CC^N\ot\End \CC^N,
\een
for any two indices $a,b\in\{1,\dots,m\}$ such that $a<b$,
we set
\beql{ars}
C_{ab}=\sum_{i,j,k,l=1}^N c^{}_{ijkl}\ts
1^{\ot(a-1)}\ot e_{ij}\ot 1^{\ot(b-a-1)}\ot e_{kl}\ot 1^{\ot(m-b)}\in\End(\CC^{N})^{\ot m}.
\eeq
We will keep the same
notation $C_{ab}$ for the element $C_{ab}\ot1$ of the algebra
\eqref{tenprkea}.

Set
\ben
\kappa=\begin{cases} N/2-1&\qquad\text{in the orthogonal case,}\\[0.3em]
N/2+1&\qquad\text{in the symplectic case}.
\end{cases}
\een
As defined in \cite{zz:rf},
the $R$-{\em matrix} $R(u)$ is a rational function in a complex parameter $u$
with values in the tensor product algebra
$\End\CC^N\ot\End\CC^N$ given by
\beql{zamolr}
R(u)=1-\frac{P}{u}+\frac{Q}{u-\kappa},
\eeq
where
\beql{p}
P=\sum_{i,j=1}^N e_{ij}\ot e_{ji},
\eeq
while $Q$ is defined by the formulas
\beql{qq}
Q=\sum_{i,j=1}^N e_{ij}
\ot e_{i'j'}
\Fand
Q=\sum_{i,j=1}^N \ve_i\tss\ve_j\ts e_{ij}
\ot e_{i'j'},
\eeq
in the
orthogonal and symplectic case, respectively. Note the relations
$
P^2=1,\ Q^2=N\ts Q
$
and
\ben
PQ=QP=\begin{cases} \phantom{-}Q&\qquad\text{in the orthogonal case,}\\
-Q&\qquad\text{in the symplectic case}.
\end{cases}
\een
The rational function \eqref{zamolr} satisfies the Yang--Baxter equation
\beql{yberep}
R_{12}(u-v)\ts R_{13}(u)\ts R_{23}(v)
=R_{23}(v)\ts R_{13}(u)\ts R_{12}(u-v).
\eeq

The {\em extended Yangian}
$\X(\g_N)$
is a unital associative algebra with generators
$t_{ij}^{(r)}$, where $1\leqslant i,j\leqslant N$ and $r=1,2,\dots$,
satisfying certain quadratic relations.
Introduce the formal series
\beql{tiju}
t_{ij}(u)=\de_{ij}+\sum_{r=1}^{\infty}t_{ij}^{(r)}\ts u^{-r}
\in\X(\g_N)[[u^{-1}]]
\eeq
and set
\ben
T(u)=\sum_{i,j=1}^N e_{ij}\ot t_{ij}(u)
\in \End\CC^N\ot \X(\g_N)[[u^{-1}]].
\een
The defining relations for the algebra $\X(\g_N)$ are then written in the form
\beql{RTTbcd}
R_{12}(u-v)\ts T_1(u)\ts T_2(v)=T_2(v)\ts T_1(u)\ts R_{12}(u-v).
\eeq
The {\em Yangian}\footnote{Since we will work with this $R$-matrix presentation
of the Yangian most of the time, we will
suppress the superscript $R$ in the notation $\Y^{R}(\g_N)$ used in the Introduction.}
$\Y(\g_N)$ is defined as the subalgebra of
$\X(\g_N)$ which
consists of the elements stable under
the automorphisms
\beql{muf}
\mu_f:T(u)\mapsto f(u)\ts T(u),
\eeq
for all series
$f(u)=1+f_1u^{-1}+f_2 u^{-2}+\cdots$
with $f_i\in\CC$.

The following tensor product decomposition holds
\beql{tensordecom}
\X(\g_N)=\ZX(\g_N)\ot \Y(\g_N),
\eeq
where $\ZX(\g_N)$ is the center of the extended Yangian $\X(\g_N)$.
The center is generated by the coefficients of the series
\beql{zn}
z_N(u)=1+\sum_{r=1}^{\infty} z^{(r)}_N\ts u^{-r}
\eeq
found by
\beql{zcenter}
T^{\tss\prime}(u+\ka)\ts T(u)=T(u)\ts T^{\tss\prime}(u+\ka)=z_N(u)\tss 1,
\eeq
where
the prime denotes the matrix transposition defined
by
\beql{transp}
(X^{\tss\prime})_{ij}=\begin{cases} X_{j'i'}\qquad&\text{in the orthogonal case,}\\
\ve_i\tss\ve_j\ts X_{j'i'}\qquad&\text{in the symplectic case.}
\end{cases}
\eeq
Equivalently,
the Yangian $\Y(\g_N)$ is the quotient of the algebra $\X(\g_N)$
by the relation $z_N(u)=1$, that is,
\beql{unita}
T^{\tss\prime}(u+\ka)\ts T(u)=1;
\eeq
see \cite{aacfr:rp} and \cite{amr:rp} for more details on the structure of the Yangian.

In terms of the series \eqref{tiju} the defining relations \eqref{RTTbcd}
can be written as
\begin{align}\label{defrel}
[\tss t_{ij}(u),t_{kl}(v)]&=\frac{1}{u-v}
\Big(t_{kj}(u)\ts t_{il}(v)-t_{kj}(v)\ts t_{il}(u)\Big)\\
{}&-\frac{1}{u-v-\kappa}
\Big(\de_{k i\pr}\sum_{p=1}^N\theta_{ip}\ts t_{pj}(u)\ts t_{p'l}(v)-
\de_{l j\pr}\sum_{p=1}^N\theta_{jp}\ts t_{k\tss p'}(v)\ts t_{ip}(u)\Big),
\non
\end{align}
where we set $\theta_{ij}\equiv 1$ in the orthogonal case, and
$\theta_{ij}=\ve_i\tss\ve_j$ in the symplectic case. Similarly,
relation \eqref{unita} reads as
\beql{unitagen}
\sum_{i=1}^N\theta_{ki}\ts t_{i\pr k'}(u+\kappa)\ts t_{il}(u)=\de_{kl}.
\eeq

\section{Embedding theorems}
\label{sec:et}

Let $A=[a_{ij}]$ be an $N\times N$ matrix over a ring with $1$.
Denote by $A^{ij}$ the matrix obtained from $A$
by deleting the $i$-th row
and $j$-th column. Suppose that the matrix
$A^{ij}$ is invertible.
The $ij$-{\em th quasideterminant of} $A$
is defined by the formula
\ben
|A|_{ij}=a_{ij}-r^{\tss j}_i(A^{ij})^{-1}\ts c^{\tss i}_j,
\een
where $r^{\tss j}_i$ is the row matrix obtained from the $i$-th
row of $A$ by deleting the element $a_{ij}$, and $c^{\tss i}_j$
is the column matrix obtained from the $j$-th
column of $A$ by deleting the element $a_{ij}$; see
\cite{gr:dm}. In particular,
the four quasideterminants of a $2\times 2$ matrix $A$ are
\ben
\bal
&|A|_{11}=a^{}_{11}-a^{}_{12}\ts a_{22}^{-1}\ts a^{}_{21},\qquad
|A|_{12}=a^{}_{12}-a^{}_{11}\ts a_{21}^{-1}\ts a^{}_{22},\\
&|A|_{21}=a^{}_{21}-a^{}_{22}\ts a_{12}^{-1}\ts a^{}_{11},\qquad
|A|_{22}=a^{}_{22}-a^{}_{21}\ts a_{11}^{-1}\ts a^{}_{12}.
\eal
\een
The quasideterminant $|A|_{ij}$ is also denoted
by boxing the entry $a_{ij}$,
\ben
|A|_{ij}=\left|\begin{matrix}a_{11}&\dots&a_{1j}&\dots&a_{1N}\\
                                   &\dots&      &\dots&      \\
                             a_{i1}&\dots&\boxed{a_{ij}}&\dots&a_{iN}\\
                                   &\dots&      &\dots&      \\
                             a_{N1}&\dots&a_{Nj}&\dots&a_{NN}
                \end{matrix}\right|.
\een

Now suppose that $n\geqslant 1$ in the case $B$, and $n\geqslant 2$ in the cases $C$ and $D$.
With the given value of $n$, consider the algebra $\X(\g_{N-2})$ and let
the indices
of the generators $t_{ij}^{(r)}$ range over the sets
$2\leqslant i,j\leqslant 2\pr$ and $r=1,2,\dots$ (the prime refers to $\g_N$ so that
$i\pr=N-i+1$).

The following is our first main result.

\bth\label{thm:embed}
The mapping
\beql{embedgen}
t_{ij}(u)\mapsto \left|\begin{matrix}
t_{11}(u)&t_{1j}(u)\\
t_{i1}(u)&\boxed{t_{ij}(u)}
\end{matrix}\right|,\qquad 2\leqslant i,j\leqslant 2\pr,
\eeq
defines an injective algebra homomorphism $\X(\g_{N-2})\to \X(\g_{N})$.
Moreover, its restriction to the subalgebra $\Y(\g_{N-2})$ defines
an injective algebra homomorphism $\Y(\g_{N-2})\to \Y(\g_{N})$.
\eth

\bpf
Denote by $s_{ij}(u)$ the quasideterminant appearing in \eqref{embedgen} so that
\beql{ttil}
s_{ij}(u)=t_{ij}(u)-t_{i1}(u)\tss t_{11}(u)^{-1}\tss t_{1j}(u).
\eeq
We start by verifying that the series $s_{ij}(u)$ satisfy the defining relations
for $\X(\g_{N-2})$. We will do this by
connecting the quasideterminants with
quantum minors of the matrix $T(u)$. Introduce power series
$\tau^{\tss a_1 a_2}_{\tss b_1 b_2}(u)$ in $u^{-1}$
with coefficients in $\X(\g_{N})$
as the matrix elements
of either side of \eqref{RTTbcd} with $v=u-1$:
\beql{quamintau}
R_{12}(1)\ts T_1(u)\ts T_2(u-1)=\sum_{a_i,b_i}e_{a_1b_1}\ot e_{a_2b_2}
\ot \tau^{\tss a_1 a_2}_{\tss b_1 b_2}(u).
\eeq

\ble\label{lem:skewsymm}\quad
(i)\quad
If $a_1\ne a'_2$ then $\tau^{\tss a_1 a_2}_{\tss b_1 b_2}(u)
=-\tau^{\tss a_2 a_1}_{\tss b_1 b_2}(u)$.

\medskip
\noindent
(ii)\quad
If $b_1\ne b'_2$ then $\tau^{\tss a_1 a_2}_{\tss b_1 b_2}(u)
=-\tau^{\tss a_1 a_2}_{\tss b_2 b_1}(u)$.
\ele

\bpf
The operator $R_{12}(1)$ remains unchanged if we multiply it from
the left or from the right by $(1-P_{12})/2+Q_{12}/N$ in the orthogonal case, and by
$(1-P_{12})/2$ in the symplectic case. By applying multiplication from the left
to the left hand side of \eqref{quamintau} we derive
part $(i)$. Applying \eqref{RTTbcd} and using
the respective multiplications of
the left hand side of \eqref{quamintau} from the right we get part $(ii)$.
\epf

\bre\label{rem:unnecc}
As the proof of Lemma~\ref{lem:skewsymm} shows, the assumptions are not necessary
in the symplectic case for the skew-symmetry properties to hold.
\qed
\ere

\ble\label{lem:toneone}
For any $2\leqslant i,j\leqslant 2\pr$ we have
\beql{wttfo}
s_{ij}(u)=t_{11}(u+1)^{-1}\tss \tau^{\tss 1 i}_{\tss 1 j}(u+1).
\eeq
Moreover,
\beql{commtoo}
\big[t_{11}(u), \tau^{\tss 1 i}_{\tss 1 j}(v)\big]=0.
\eeq
\ele

\bpf
By \eqref{ttil} we have
\ben
t_{11}(u+1)\ts s_{ij}(u)=t_{11}(u+1)\tss
t_{ij}(u)-t_{11}(u+1)\tss t_{i1}(u)\tss t_{11}(u)^{-1}\tss t_{1j}(u).
\een
However, $t_{11}(u+1)\tss t_{i1}(u)=t_{i1}(u+1)\tss t_{11}(u)$ by \eqref{defrel},
so that
\ben
t_{11}(u+1)\ts s_{ij}(u)=t_{11}(u+1)\tss
t_{ij}(u)-t_{i1}(u+1)\tss t_{1j}(u).
\een
The definition \eqref{quamintau} implies that this
coincides with $\tau^{\tss 1 i}_{\tss 1 j}(u+1)$ hence \eqref{wttfo} follows.
Relation \eqref{commtoo} follows easily from \eqref{defrel}. It can also be
derived by noting that the commutation relations between the series involved
in this calculation are the same as for
the Yangian $\Y(\gl_N)$. Therefore, \eqref{commtoo} holds due to
the corresponding properties of the quantum minors for $\Y(\gl_N)$;
see, e.g., \cite[Section~1.7]{m:yc}.
\epf

Now we will need some simplified expressions for both sides of \eqref{yberep}
when $v=u-1$.

\ble\label{lem:ybesi}
We have the relations
\begin{multline}\label{fus}
R_{12}(1)\ts R_{13}(u)\ts R_{23}(u-1)\\[0.3em]
=R_{12}(1)\ts
\Big(1-\frac{P_{13}+P_{23}}{u-1}+\frac{Q_{13}+Q_{23}}{u-\ka}
-\frac{P_{23}\tss Q_{12}}{(u-1)(u-\ka)}-P_{13}\tss Q_{23}\ts\vp(u)\Big)
\end{multline}
and
\begin{multline}\label{fusopp}
R_{23}(u-1)\ts R_{13}(u)\ts R_{12}(1)\\[0.3em]
=\Big(1-\frac{P_{13}+P_{23}}{u-1}+\frac{Q_{13}+Q_{23}}{u-\ka}
-\frac{Q_{12}\tss P_{23}}{(u-1)(u-\ka)}-Q_{23}\tss P_{13}\ts\vp(u)\Big)\ts R_{12}(1),
\end{multline}
where
\ben
\vp(u)=\begin{cases} \dfrac{1-4/N}{(u-1)(u-\ka)}&\qquad\text{in the orthogonal case,}\\[0.9em]
\dfrac{1}{u(u-\ka-1)}&\qquad\text{in the symplectic case}.
\end{cases}
\een
\ele

\bpf
The product of $R$-matrices $R_{13}(u)\ts R_{23}(u-1)$ on
the left hand side of \eqref{fus} equals
\begin{multline}
1-\frac{P_{13}}{u}-\frac{P_{23}}{u-1}+\frac{Q_{13}}{u-\ka}+
\frac{Q_{23}}{u-\ka-1}+\frac{P_{13}\tss P_{23}}{u(u-1)}\\[0.3em]
{}-
\frac{P_{13}\tss Q_{23}}{u(u-\ka-1)}-\frac{Q_{13}\tss P_{23}}{(u-1)(u-\ka)}
+\frac{Q_{13}\tss Q_{23}}{(u-\ka)(u-\ka-1)}.
\non
\end{multline}
We have
\ben
(1-P_{12})\tss P_{13}\tss P_{23}=(1-P_{12})\tss P_{12}\tss P_{13}=-(1-P_{12})\tss P_{13}
\een
and
\ben
(1-P_{12})\tss Q_{13}\tss Q_{23}=(1-P_{12})\tss Q_{13}\tss P_{12}=-(1-P_{12})\tss Q_{23}.
\een
Now continuing with the symplectic case, we also have
\ben
Q_{12}\tss P_{13}\tss P_{23}=Q_{12}\tss P_{12}\tss P_{13}=-Q_{12}\tss P_{13}
\een
and
\ben
Q_{12}\tss Q_{13}\tss Q_{23}=Q_{12}\tss \tss Q_{13}\tss P_{12}=-Q_{12}\tss Q_{23}.
\een
Since $Q_{13}\tss P_{23}=P_{23}\tss Q_{12}$, we can obtain equation \eqref{fus} for the symplectic case.

In the orthogonal case, the last two relations take the different form
\ben
Q_{12}\tss P_{13}\tss P_{23}=Q_{12}\tss P_{13}\Fand
Q_{12}\tss Q_{13}\tss Q_{23}=Q_{12}\tss Q_{23}.
\een
Write $P_{13}\tss Q_{23}=Q_{12}\tss P_{13}$ and use the relations
$(1-P_{12})\tss Q_{12}=0$ and $Q_{12}^2=N\tss Q_{12}$
together with $Q_{12}\tss Q_{23}=Q_{12}\tss P_{13}$
to express the left hand side
of \eqref{fus} as
\begin{multline}
(1-P_{12})\Big(1-\frac{P_{13}+P_{23}}{u-1}+\frac{Q_{13}+Q_{23}}{u-\ka}
-\frac{Q_{13}\tss P_{23}}{(u-1)(u-\ka)}\Big)
\non\\[0.4em]
{}+\frac{Q_{12}}{1-\ka}\Big(1-\frac{P_{13}+P_{23}}{u-1}+\frac{Q_{13}+Q_{23}}{u-\ka}
-\frac{Q_{13}\tss P_{23}+2\tss(\ka-1)\tss
P_{13}}{(u-1)(u-\ka)}\Big).
\end{multline}
This coincides with the right hand side of \eqref{fus}.

The proof of \eqref{fusopp} is obtained by
using the same arguments and
writing the products
of the $P$ and $Q$ operators in the reverse order.
\epf

\ble\label{lem:sylv}
The mapping
\ben
t_{ij}(u)\mapsto \tau^{\tss 1 i}_{\tss 1 j}(u),\qquad 2\leqslant i,j\leqslant 2\pr,
\een
defines a homomorphism $\X(\g_{N-2})\to \X(\g_{N})$.
\ele

\bpf
Consider the tensor product algebra \eqref{tenprkea} with $m=4$, which is associated
with the extended Yangian. We have the relation
\begin{multline}
\label{longrel}
R_{23}(a-1)\tss R_{13}(a)\tss R_{24}(a)\tss R_{14}(a+1)\tss
R_{12}(1)\tss T_1(u)\tss T_2(u-1)
R_{34}(1)\tss T_3(v)\tss T_4(v-1)\\[0.4em]
{}=R_{34}(1)\tss T_3(v)\tss T_4(v-1)\tss R_{12}(1)\tss T_1(u)\tss T_2(u-1)\tss
R_{14}(a+1)\tss R_{24}(a)\tss R_{13}(a)\tss R_{23}(a-1),
\end{multline}
where $a=u-v$. It follows easily by a repeated application of
the Yang--Baxter relation \eqref{yberep} and
the $RTT$ relation \eqref{RTTbcd}. We will
transform the operators on both sides of \eqref{longrel} by using Lemma~\ref{lem:ybesi}
and then equate some matrix elements.
We begin with the right hand side and apply first
\eqref{RTTbcd} with $v=u-1$ to write the product $R_{12}(1)\tss T_1(u)\tss T_2(u-1)$
in the reverse order. Next, use \eqref{yberep} to write
\ben
R_{12}(1)\tss R_{14}(a+1)\tss R_{24}(a)\tss R_{13}(a)\tss R_{23}(a-1)
=R_{24}(a)\tss R_{14}(a+1)\tss R_{12}(1)\tss R_{13}(a)\tss R_{23}(a-1)
\een
and apply \eqref{fus} with $u$ replaced by $a$
to the last three factors. Then use \eqref{yberep}
again, and apply \eqref{fus} with $u$ replaced by $a+1$
to the product
$
R_{12}(1)\tss R_{14}(a+1)\tss R_{24}(a).
$
As a result, the right hand side of \eqref{longrel} is transformed
in such a way that the last four factors are replaced with the product
\begin{multline}
\label{longpro}
\Big(1-\frac{P_{14}+P_{24}}{a}+\frac{Q_{14}+Q_{24}}{a-\ka+1}
-\frac{P_{24}\tss Q_{12}}{a(a-\ka+1)}-P_{14}\tss Q_{24}\ts\vp(a+1)\Big)\\[0.5em]
{}\times\Big(1-\frac{P_{13}+P_{23}}{a-1}+\frac{Q_{13}+Q_{23}}{a-\ka}
-\frac{P_{23}\tss Q_{12}}{(a-1)(a-\ka)}-P_{13}\tss Q_{23}\ts\vp(a)\Big).
\end{multline}
Now apply the operator on the right hand side of \eqref{longrel} to a
basis vector of the form $e_1\ot e_j\ot e_1\ot e_l$
for certain
$j,l\in\{2,\dots,2\pr\}$.
Each of the operators $Q_{12}, Q_{13}$ and $Q_{23}$
annihilates the vector, so that the application of the second factor in \eqref{longpro}
gives
\beql{vectocc}
\frac{a-2}{a-1}\ts e_1\ot e_j\ot e_1\ot e_l-\frac{1}{a-1}\ts e_1\ot e_1\ot e_j\ot e_l.
\eeq
Next apply the first factor in \eqref{longpro} to each of the vectors
occurring in \eqref{vectocc}. The operators $Q_{12}$ and $Q_{14}$ annihilate
the vector $e_1\ot e_j\ot e_1\ot e_l$, while
\ben
P_{14}\ts(e_1\ot e_j\ot e_1\ot e_l)=e_l\ot e_j\ot e_1\ot e_1.
\een
The vector $e_l\ot e_j\ot e_1\ot e_1$ will be annihilated by a subsequent
application of the operator $R_{34}(1)\tss T_3(v)\tss T_4(v-1)$ due to
Lemma~\ref{lem:skewsymm}\tss$(ii)$. The same property holds for the vector
$P_{14}\tss Q_{24}\ts(e_1\ot e_j\ot e_1\ot e_l)$.
The application of the first factor in \eqref{longpro} to the second vector
$e_1\ot e_1\ot e_j\ot e_l$ in \eqref{vectocc} gives
\ben
e_1\ot e_1\ot e_j\ot e_l-\frac{1}{a}\ts(e_l\ot e_1\ot e_j\ot e_1+e_1\ot e_l\ot e_j\ot e_1).
\een
By Lemma~\ref{lem:skewsymm}\tss$(ii)$, this expression will be annihilated by a subsequent
application of the operator $R_{12}(1)\tss T_1(u)\tss T_2(u-1)$.

We may conclude that the restriction of the operator on the right hand side of
\eqref{longrel} to the subspace spanned by the
basis vectors of the form $e_1\ot e_j\ot e_1\ot e_l$ with $j,l\in\{2,\dots,2\pr\}$
coincides with the operator
\beql{oprhs}
\frac{a-2}{a-1}\ts R_{34}(1)\tss T_3(v)\tss T_4(v-1)\tss
R_{12}(1)\tss T_1(u)\tss T_2(u-1)\tss \Big(1-\frac{P_{24}}{a}+\frac{Q_{24}}{a-\ka+1}\Big).
\eeq

Now consider the operator on the left hand side of \eqref{longrel}. We will apply
it to basis vectors of $(\CC^N)^{\ot 4}$ and look at the coefficients of the basis
vectors of the form $e_1\ot e_i\ot e_1\ot e_k$ with $i,k\in\{2,\dots,2\pr\}$
in the image. The same argument as for the right hand side, with the use
of \eqref{fusopp} and Lemma~\ref{lem:skewsymm}\tss$(i)$ instead, and with
reversed factors in the operators, implies that the coefficients
of such basis vectors coincide with those of the operator
\beql{oplhs}
\frac{a-2}{a-1}\ts \Big(1-\frac{P_{24}}{a}+\frac{Q_{24}}{a-\ka+1}\Big)\ts
R_{12}(1)\tss T_1(u)\tss T_2(u-1)\tss
R_{34}(1)\tss T_3(v)\tss T_4(v-1).
\eeq
Furthermore, the application to the vectors $e_1\ot e_j\ot e_1\ot e_l$ with $j,l\in\{2,\dots,2\pr\}$
gives
\begin{multline}
R_{12}(1)\tss T_1(u)\tss T_2(u-1)\tss
R_{34}(1)\tss T_3(v)\tss T_4(v-1)\tss(e_1\ot e_j\ot e_1\ot e_l)\\
{}\equiv\sum_{c,d=1}^{1\pr}\tau^{\tss 1 c}_{\tss 1 j}(u)\tss \tau^{\tss 1 d}_{\tss 1 l}(v)
\tss(e_1\ot e_c\ot e_1\ot e_d),
\label{tautwo}
\end{multline}
where we only keep the basis vectors which can give a nonzero contribution
to the coefficient of $e_1\ot e_i\ot e_1\ot e_k$ after the subsequent
application of the operators $1, P_{24}$ or $Q_{24}$. Moreover, by
Lemma~\ref{lem:skewsymm}\tss$(i)$, the values $c=1$ and $d=1$ can also be
excluded from the range of the summation indices. This implies that the values
$c=1\pr$ and $d=1\pr$ can be excluded as well, and so we can write
an operator equality
\ben
1-\frac{P_{24}}{a}+\frac{Q_{24}}{a-\ka+1}=R_{24}(u-v),
\een
which is the $R$-matrix
associated with the algebra $\X(\g_{N-2})$.
The same argument with the use of Lemma~\ref{lem:skewsymm}\tss$(ii)$
shows that this equality can also be used
for the operator \eqref{oprhs}. In other words, by equating
the matrix elements of the operators \eqref{oprhs} and \eqref{oplhs}
we get the $R$-matrix form of the defining relations for the algebra
$\X(\g_{N-2})$ satisfied by the series
$\tau^{\tss 1 i}_{\tss 1 j}(u)$, as required.
\epf

Returning to the proof of the theorem, we can now show that the map \eqref{embedgen}
defines a homomorphism. By taking the composition of the
homomorphism of Lemma~\ref{lem:sylv} with the shift automorphism $T(u)\mapsto T(u+1)$
we get another homomorphism $\X(\g_{N-2})\to \X(\g_{N})$
defined by
$t_{ij}(u)\mapsto \tau^{\tss 1 i}_{\tss 1 j}(u+1)$. It remains to apply Lemma~\ref{lem:toneone}
and note the commutation relations $\big[t_{11}(u),t_{11}(v)\big]=0$.

Next we show that the homomorphism \eqref{embedgen} is injective. For each $N$ introduce
an ascending filtration on the extended Yangian $\X(\g_{N})$ by setting $\deg t_{ij}^{(r)}=r-1$
for all $r\geqslant 1$. Denote by $\bar t_{ij}^{\ts(r)}$ the image of $t_{ij}^{(r)}$
in the $(r-1)$-th component of the associated graded algebra $\gr\X(\g_{N})$.
The map \eqref{embedgen} defines a
homomorphism of the associated graded algebras $\gr\X(\g_{N-2})\to \gr\X(\g_{N})$.
It takes the generator $\bar t_{ij}^{\ts(r)}\in \gr\X(\g_{N-2})$ to the element of $\gr\X(\g_{N})$
denoted by the same symbol. As shown
in the proof of the Poincar\'e--Birkhoff--Witt theorem for the extended
Yangian \cite[Corollary~3.10]{amr:rp}, the mapping
\beql{isomgrpol}
\bar t_{ij}^{\ts(r)}\mapsto F_{ij}\tss x^{r-1}+\frac12\ts\de_{ij}\ts\ze_r
\eeq
defines an isomorphism
\ben
\gr\X(\g_{N})\cong \U(\g_N[x])\ot\CC[\ze_1,\ze_2,\dots],
\een
where $\CC[\ze_1,\ze_2,\dots]$ is the algebra of polynomials in variables $\ze_i$.
These variables correspond to the images of the central elements $z^{(r)}_N$
defined in \eqref{zn},
\beql{znim}
\bar z^{\ts(r)}_N\mapsto \ze_r.
\eeq
Therefore the homomorphism $\gr\X(\g_{N-2})\to \gr\X(\g_{N})$
is injective and so is the homomorphism \eqref{embedgen}.

Finally, observe that the homomorphism \eqref{embedgen} commutes with the
automorphism $\mu_f$ defined in \eqref{muf} associated with an arbitrary
series $f(u)=1+f_1u^{-1}+f_2 u^{-2}+\cdots$ in the sense that the following
diagram commutes:
\ben
\begin{CD}
\X(\g_{N-2})@>>> \X(\g_{N})\\
@V\mu_f VV     @VV\mu_f V\\
\X(\g_{N-2})@>>> \X(\g_{N}),
\end{CD}
\een
where the horizontal arrows denote the homomorphism \eqref{embedgen}. Therefore,
the image of the restriction of this homomorphism to the subalgebra $\Y(\g_{N-2})$
of $\X(\g_{N-2})$
is contained in the subalgebra $\Y(\g_{N})$ of $\X(\g_{N})$. This restriction thus
defines an injective homomorphism $\Y(\g_{N-2})\to \Y(\g_{N})$.
\epf

The following generalization of
Theorem~\ref{thm:embed} will be useful for our arguments below.
Fix a positive integer $m$ such that
$m\leqslant n$ for type $B$ and $m\leqslant n-1$ for types $C$ and $D$.
Suppose that the generators $t_{ij}^{(r)}$ of the algebra $\X(\g_{N-2m})$ are
labelled by the indices
$m+1\leqslant i,j\leqslant (m+1)\pr$ and $r=1,2,\dots$ with $i\pr=N-i+1$ as before.

\bth\label{thm:red}
The mapping
\beql{redu}
\psi_m:t_{ij}(u)\mapsto \left|\begin{matrix}
t_{11}(u)&\dots&t_{1m}(u)&t_{1j}(u)\\
\dots&\dots&\dots&\dots\\
t_{m1}(u)&\dots&t_{mm}(u)&t_{mj}(u)\\
t_{i1}(u)&\dots&t_{im}(u)&\boxed{t_{ij}(u)}
\end{matrix}\right|,\qquad m+1\leqslant i,j\leqslant (m+1)\pr,
\eeq
defines an injective homomorphism $\X(\g_{N-2m})\to \X(\g_{N})$.
Moreover, its restriction to the subalgebra $\Y(\g_{N-2m})$ defines
an injective homomorphism $\Y(\g_{N-2m})\to \Y(\g_{N})$.
\eth

\bpf
We argue by induction on $m$. The case $m=1$ is Theorem~\ref{thm:embed}.
Suppose that $m\geqslant 2$. By the Sylvester theorem for quasideterminants
\cite{gr:dm} (see also \cite{kl:mi} for a proof),
we have the identity
\ben
\begin{vmatrix}
t_{11}(u)&\dots&t_{1m}(u)&t_{1j}(u)\\
\dots&\dots&\dots&\dots\\
t_{m1}(u)&\dots&t_{mm}(u)&t_{mj}(u)\\
t_{i1}(u)&\dots&t_{im}(u)&\boxed{t_{ij}(u)}
\end{vmatrix}=\begin{vmatrix}
s_{22}(u)&\dots&s_{2m}(u)&s_{2j}(u)\\
\dots&\dots&\dots&\dots\\
s_{m2}(u)&\dots&s_{mm}(u)&s_{mj}(u)\\
s_{i2}(u)&\dots&s_{im}(u)&\boxed{s_{ij}(u)}
\end{vmatrix},
\een
where
\ben
s_{ab}(u)=\left|\begin{matrix}
t_{11}(u)&t_{1b}(u)\\
t_{a1}(u)&\boxed{t_{ab}(u)}
\end{matrix}\right|.
\een
By Theorem~\ref{thm:embed}, the mapping $t_{ab}(u)\mapsto s_{ab}(u)$
with $2\leqslant a,b\leqslant 2\pr$ defines a homomorphism $\X(\g_{N-2})\to \X(\g_{N})$.
Furthermore, by the induction hypothesis, the map
\ben
t_{ij}(u)\mapsto \left|\begin{matrix}
s_{22}(u)&\dots&s_{2m}(u)&s_{2j}(u)\\
\dots&\dots&\dots&\dots\\
s_{m2}(u)&\dots&s_{mm}(u)&s_{mj}(u)\\
s_{i2}(u)&\dots&s_{im}(u)&\boxed{s_{ij}(u)}
\end{matrix}\right|,\qquad m+1\leqslant i,j\leqslant (m+1)\pr,
\een
defines a homomorphism $\X(\g_{N-2m})\to \X(\g_{N-2})$ thus proving that \eqref{redu}
is a homomorphism. Its injectivity and the remaining parts of the theorem are
verified in the same way as for Theorem~\ref{thm:embed}.
\epf

We will point out a consistence property of the embeddings \eqref{redu}
whose particular case $l=1$ was already used in the proof of Theorem~\ref{thm:red};
cf.~\cite[eq.~(1.85)]{m:yc} for its counterpart in type $A$.
We will write $\psi_m=\psi^{(N)}_m$ to indicate the dependence of $N$.
For a parameter $l$ we have the corresponding embedding
\ben
\psi^{(N-2l)}_m:\X(\g_{N-2m-2l})\hra \X(\g_{N-2l})
\een
provided by \eqref{redu}.

\bpr\label{prop:consist}
We have the equality of maps
\ben
\psi^{(N)}_l\circ\psi^{(N-2l)}_m=\psi^{(N)}_{l+m}.
\een
\epr

\bpf
For all $l+1\leqslant a,b\leqslant (l+1)\pr$ introduce the series $s_{ab}(u)$
with coefficients in $\X(\g_N)$ by
\ben
s_{ab}(u)=\begin{vmatrix}
t_{11}(u)&\dots&t_{1l}(u)&t_{1b}(u)\\
\dots&\dots&\dots&\dots\\
t_{l1}(u)&\dots&t_{l\tss l}(u)&t_{lb}(u)\\
t_{a1}(u)&\dots&t_{al}(u)&\boxed{t_{ab}(u)}
\end{vmatrix}.
\een
The desired equality amounts to the identity for series
with coefficients in $\X(\g_N)$,
\ben
\begin{vmatrix}
s_{l+1\ts l+1}(u)&\dots&s_{l+1\ts l+m}(u)&s_{l+1\ts j}(u)\\
\dots&\dots&\dots&\dots\\
s_{l+m\ts l+1}(u)&\dots&s_{l+m\ts l+m}(u)&s_{l+m\ts j}(u)\\
s_{i\ts l+1}(u)&\dots&s_{i\ts l+m}(u)&\boxed{s_{ij}(u)}
\end{vmatrix}
=\begin{vmatrix}
t_{11}(u)&\dots&t_{1\ts l+m}(u)&t_{1j}(u)\\
\dots&\dots&\dots&\dots\\
t_{l+m\ts 1}(u)&\dots&t_{l+m\ts l+m}(u)&t_{l+m\ts j}(u)\\
t_{i1}(u)&\dots&t_{i\ts l+m}(u)&\boxed{t_{ij}(u)}
\end{vmatrix}
\een
which holds for all $l+m+1\leqslant i,j\leqslant (l+m+1)\pr$
due to the Sylvester theorem for quasideterminants
\cite{gr:dm}, \cite{kl:mi}.
\epf

For subsets $\{a_1,\dots,a_k\}$ and $\{b_1,\dots,b_k\}$ of $\{1,\dots,N\}$
introduce $A$-type quantum minors by the formula
\ben
{t\ts}^{a_1\dots\ts a_k}_{b_1\dots\ts b_k}(u)=
\sum_{ p\in \Sym_k} \sgn p\cdot t_{a_{ p(1)}b_1}(u)\dots
t_{a_{ p(k)}b_k}(u-k+1).
\een
These are formal series in $u^{-1}$ with coefficients in $\X(\g_{N})$.

\bpr\label{prop:quasiqua}
For all $m+1\leqslant i,j\leqslant (m+1)\pr$ we have the identity
\ben
\left|\begin{matrix}
t_{11}(u)&\dots&t_{1m}(u)&t_{1j}(u)\\
\dots&\dots&\dots&\dots\\
t_{m1}(u)&\dots&t_{mm}(u)&t_{mj}(u)\\
t_{i1}(u)&\dots&t_{im}(u)&\boxed{t_{ij}(u)}
\end{matrix}\right|={t\ts}^{1\dots\ts m}_{1\dots\ts m}(u+m)^{-1}\cdot
{t\ts}^{1\dots\ts m\ts i}_{1\dots\ts m\ts j}(u+m).
\een
\epr

\bpf
This identity holds for the Yangian $\Y(\gl_N)$; see, e.g., \cite[Section~1.11]{m:yc}.
The commutation relations between the generator series $t_{ab}(u)$ of $\X(\g_{N})$
occurring in this identity are the same as for the Yangian $\Y(\gl_N)$, with a possible
exception of the commutators $[t_{ij}(u),t_{ij}(v)]$, and in addition in the $B$ case,
the commutators of the series in the last row or last column of the quasideterminant.
However, since the quasideterminant depends linearly on such generators,
the identity does not depend on such
commutators. Hence it holds for $\X(\g_{N})$ as well.
\epf

The following is a counterpart of the corresponding property of the Yangian
for $\gl_N$; see, e.g., \cite{bk:pp}.

\bco\label{cor:commu}
We have the relations
\ben
\big[t_{ab}(u),\psi_m(t_{ij}(v))\big]=0
\een
for all $1\leqslant a,b\leqslant m$ and $m+1\leqslant i,j\leqslant (m+1)\pr$.
\eco

\bpf
By Proposition~\ref{prop:quasiqua} we only need to verify that $t_{ab}(u)$
commutes with the quantum minors. This follows by the same argument as for the proof
of the proposition.
\epf

\section{Gauss decomposition}
\label{sec:gd}

As we pointed out in the Introduction, the Gauss decomposition \eqref{gd}
will play a key role in constructing the Drinfeld
generators. We will assume that $T(u)$ is the generator matrix
for the extended Yangian $\X(\g_N)$ (that is, we ignore relation \eqref{unitaryint})
and recall the well-known formulas for the entries
of the matrices $F(u)$, $H(u)$ and $E(u)$ which occur in \eqref{gd};
see, e.g., \cite[Sec.~1.11]{m:yc}. We have
\beql{hmqua}
h_i(u)=\begin{vmatrix} t_{1\tss 1}(u)&\dots&t_{1\ts i-1}(u)&t_{1\tss i}(u)\\
                          \vdots&\ddots&\vdots&\vdots\\
                         t_{i-1\ts 1}(u)&\dots&t_{i-1\ts i-1}(u)&t_{i-1\ts i}(u)\\
                         t_{i\tss 1}(u)&\dots&t_{i\ts i-1}(u)&\boxed{t_{i\tss i}(u)}\\
           \end{vmatrix},\qquad i=1,\dots,N,
\eeq
whereas
\beql{eijmlqua}
e_{ij}(u)=h_i(u)^{-1}\ts\begin{vmatrix} t_{1\tss 1}(u)&\dots&t_{1\ts i-1}(u)&t_{1\ts j}(u)\\
                          \vdots&\ddots&\vdots&\vdots\\
                         t_{i-1\ts 1}(u)&\dots&t_{i-1\ts i-1}(u)&t_{i-1\ts j}(u)\\
                         t_{i\tss 1}(u)&\dots&t_{i\ts i-1}(u)&\boxed{t_{i\tss j}(u)}\\
           \end{vmatrix}
\eeq
and
\beql{fijlmqua}
f_{ji}(u)=\begin{vmatrix} t_{1\tss 1}(u)&\dots&t_{1\ts i-1}(u)&t_{1\tss i}(u)\\
                          \vdots&\ddots&\vdots&\vdots\\
                         t_{i-1\ts 1}(u)&\dots&t_{i-1\ts i-1}(u)&t_{i-1\ts i}(u)\\
                         t_{j\ts 1}(u)&\dots&t_{j\ts i-1}(u)&\boxed{t_{j\tss i}(u)}\\
           \end{vmatrix}\ts h_i(u)^{-1}
\eeq
for $1\leqslant i<j\leqslant N$. Obviously, the algebra $\X(\g_N)$
is generated by the coefficients of the series $f_{ji}(u)$, $e_{ij}(u)$ and $h_i(u)$
which we will refer to as the {\em Gaussian generators}.

Suppose that $0\leqslant m< n$ if $N=2n$ and $0\leqslant m\leqslant n$ if $N=2n+1$.
We will use the superscript $[m]$ to indicate
square submatrices corresponding to rows and columns labelled by
$m+1,m+2,\dots,(m+1)'$. In particular, we set
\ben
F^{[m]}(u)=\begin{bmatrix}
1&0&\dots&0\ts\\
f_{m+2\ts m+1}(u)&1&\dots&0\\
\vdots&\ddots&\ddots&\vdots\\
f_{(m+1)'\tss m+1}(u)&\dots&f_{(m+1)'\ts (m+2)'}(u)&1
\end{bmatrix},
\een
\ben
E^{[m]}(u)=\begin{bmatrix} 1&e_{m+1\tss m+2}(u)&\ldots&e_{m+1\tss (m+1)'}(u)\\
                        0&1&\ddots &\vdots\\
                         \vdots&\vdots&\ddots&e_{(m+2)'\tss(m+1)'}(u)\\
                         0&0&\ldots&1\\
           \end{bmatrix}
\een
and $H^{[m]}(u)=\diag\ts\big[h_{m+1}(u),\dots,h_{(m+1)'}(u)\big]$. Furthermore,
introduce the product of these matrices by
\ben
T^{[m]}(u)=F^{[m]}(u)\tss H^{[m]}(u)\tss E^{[m]}(u).
\een
Accordingly, the entries of $T^{[m]}(u)$ will be denoted by $t^{[m]}_{ij}(u)$
with $m+1\leqslant i,j\leqslant (m+1)'$.

\bpr\label{prop:gauss-consist}
The series $t^{[m]}_{ij}(u)$ coincides with the image of the generator series $t_{ij}(u)$
of the extended Yangian $\X(\g_{N-2m})$
under the embedding \eqref{redu},
\ben
t^{[m]}_{ij}(u)=\psi_m\big(t_{ij}(u)\big),\qquad m+1\leqslant i,j\leqslant (m+1)'.
\een
\epr

\bpf
Set $s_{ij}(u)=\psi_m\big(t_{ij}(u)\big)$.
Since the Gauss decomposition \eqref{gd}
uniquely determines the matrices $F(u)$, $H(u)$ and $E(u)$,
it suffices to verify that such matrices obtained by the Gauss decomposition
of the matrix
\ben
S(u)=\big[s_{ij}(u)\big],\qquad m+1\leqslant i,j\leqslant (m+1)',
\een
coincide with $F^{[m]}(u)$, $H^{[m]}(u)$ and $E^{[m]}(u)$, respectively.
Let $S(u)=\wt F(u)\tss \wt H(u)\tss \wt E(u)$ be the Gauss decomposition
of $S(u)$. By formulas \eqref{hmqua}, \eqref{eijmlqua} and \eqref{fijlmqua},
the entries of the matrices $\wt F(u)$, $\wt H(u)$ and $\wt E(u)$
are found as the quasideterminants of certain submatrices of $S(u)$.
However, as we pointed out in the proof of Proposition~\ref{prop:consist},
such quasideterminants coincide with the corresponding quasideterminants of
submatrices of the matrix $T(u)$.
\epf

We record the following as an immediate consequence of Proposition~\ref{prop:gauss-consist}.

\bco\label{cor:guass-embed}
The subalgebra $\X^{[m]}(\g_N)$ generated by the coefficients of all
series $t^{[m]}_{ij}(u)$ with $m+1\leqslant i,j\leqslant (m+1)'$
is isomorphic to the extended Yangian $\X(\g_{N-2m})$. In particular,
we have the relation
\beql{subRTT}
R^{[m]}_{12}(u-v)\ts T^{[m]}_1(u)\ts T^{[m]}_2(v)=T^{[m]}_2(v)\ts T^{[m]}_1(u)\ts R^{[m]}_{12}(u-v),
\eeq
where $R^{[m]}(u)$ is the $R$-matrix associated with
$\X(\g_{N-2m})$.
Moreover,
\beql{subunitary}
T^{[m]\tss\prime}(u+\ka^{[m]})T^{[m]}(u)=T^{[m]}(u)T^{[m]\tss\prime}(u+\ka^{[m]})=z^{[m]}_N(u)\tss 1,
\eeq
with $\ka^{[m]}=\ka-m$, for a certain series
$z^{[m]}_N(u)$ whose coefficients
generate the center of the subalgebra $\X^{[m]}(\g_N)$.
\qed
\eco

Introduce the coefficients of the series defined in
\eqref{hmqua}, \eqref{eijmlqua} and \eqref{fijlmqua} by
\ben
e_{ij}(u)=\sum_{r=1}^{\infty} e_{ij}^{(r)}\tss u^{-r},\qquad
f_{ji}(u)=\sum_{r=1}^{\infty} f_{ji}^{(r)}\tss u^{-r},\qquad
h_i(u)=1+\sum_{r=1}^{\infty} h_i^{(r)}\tss u^{-r}.
\een
Furthermore, define the series by
\beql{defkn}
k_{i}(u)=h_i(u)^{-1}h_{i+1}(u),\qquad
e_{i}(u)=e_{i\ts i+1}(u),
\qquad f_{i}(u)=f_{i+1\ts i}(u)
\eeq
for $i=1,\dots, n-1$, and set
\beql{defenfn}
e_n(u)=\begin{cases}
e_{n\ts n+1}(u)
\quad&\text{for}\quad \oa_{2n+1}\\
e_{n\ts n+1}(u)
\quad&\text{for}\quad \spa_{2n}\\
e_{n-1\ts n+1}(u)
\quad&\text{for}\quad \oa_{2n}
\end{cases},
\qquad
f_n(u)=\begin{cases}
f_{n+1\ts n}(u)
\quad&\text{for}\quad \oa_{2n+1}\\
f_{n+1\ts n}(u)
\quad&\text{for}\quad \spa_{2n}\\
f_{n+1\ts n-1}(u)
\quad&\text{for}\quad \oa_{2n}
\end{cases}
\eeq
and
\beql{defknu}
k_n(u)=\begin{cases}
h_n(u)^{-1}\ts h_{n+1}(u)
\qquad&\text{for}\quad \oa_{2n+1}\\[0.2em]
2\tss h_n(u)^{-1}\ts h_{n+1}(u)
\qquad&\text{for}\quad \spa_{2n}\\[0.2em]
h_{n-1}(u)^{-1}\ts h_{n+1}(u)
\qquad&\text{for}\quad \oa_{2n}.
\end{cases}
\eeq

\ble\label{lem:emjmtkl}
For the parameter $m$ chosen as above, suppose that the indices
$j,k,l$ satisfy $m+1\leqslant j,k,l\leqslant (m+1)'$ and $j\neq l'$.
Then the following relations hold in the extended Yangian
$\X(\g_N)$,
\beql{emjmtkl}
\big[e_{m\tss j}(u), {t^{\ts[m]}_{k\tss l}(v)}\big]
=\frac{1}{u-v}{t^{\ts[m]}_{k\tss j}(v)}\big(e_{m\tss l}(v)-e_{m\tss l}(u)\big),
\eeq
\beql{fjmmtkl}
\big[f_{j\tss m}(u), {t^{\ts[m]}_{k\tss l}(v)}\big]
=\frac{1}{u-v}\big(f_{k\tss m}(u)-f_{k\tss m}(v)\big)\ts{t^{\ts[m]}_{j\tss l}(v)}.
\eeq
\ele

\bpf
It is sufficient to verify the relations for $m=1$; the general case
will follow by the application of the homomorphism $\psi_m$;
see Proposition~\ref{prop:gauss-consist}.
Both relations follow by similar arguments so we
only verify \eqref{emjmtkl}.
Since
\ben
h_1(v)=t_{11}(v),\qquad f_{k\tss 1}(v)=t_{k1}(v)t_{11}(v)^{-1},\qquad
e_{1\tss l}(v)=t_{11}(v)^{-1}t_{1\tss l}(v),
\een
we can write
\ben
t^{\ts[1]}_{kl}(v)=t_{kl}(v)-t_{k1}(v)t_{11}(v)^{-1}t_{1\tss l}(v)=t_{kl}(v)-f_{k1}(v)h_{1}(v)e_{1l}(v).
\een
The defining relations
\eqref{defrel} imply
\ben
[t_{1\tss j}(u), t_{k\tss l}(v)]=\frac{1}{u-v}\big(t_{kj}(u)t_{1l}(v)-t_{kj}(v)t_{1l}(u)\big)
\een
and so
\begin{multline}
\non
\big[t_{1\tss j}(u), {t^{\ts[1]}_{k\tss l}(v)}\big]
+\big[t_{1j}(u), f_{k\tss 1}(v)h_{1}(v)e_{1\tss l}(v)\big]
=\frac{1}{u-v}\big({t^{\ts[1]}_{kj}(u)}t_{1l}(v)-{t^{\ts[1]}_{kj}(v)}t_{1l}(u)\big)\\[0.4em]
{}+\frac{1}{u-v}\big(f_{k1}(u)h_{1}(u)e_{1j}(u)t_{1l}(v)-f_{k1}(v)h_{1}(v)e_{1j}(v)t_{1l}(u)\big).
\end{multline}
The second commutator
on the left hand side can be transformed as
\begin{multline}
\non
\big[t_{1j}(u), f_{k\tss 1}(v)h_{1}(v)e_{1\tss l}(v)\big]
=\big[t_{1j}(u), t_{k1}(v)\big]\ts e_{1\tss l}(v)
+t_{k1}(v)\big[t_{1j}(u),e_{1\tss l}(v)\big]\\[0.4em]
{}=\big[t_{1j}(u), t_{k1}(v)\big]e_{1\tss l}(v)
+t_{k1}(v)\big[t_{1j}(u),t_{11}(v)^{-1}t_{1\tss l}(v)\big],
\end{multline}
which equals
\begin{multline}
\non
\big[t_{1j}(u), t_{k1}(v)\big]e_{1\tss l}(v)
+t_{k1}(v)\big[t_{1j}(u),t_{11}(v)^{-1}\big]t_{1\tss l}(v)
+f_{k1}(v)\big[t_{1j}(u),t_{1\tss l}(v)\big]\\[0.3em]
{}=\big[t_{1j}(u), t_{k1}(v)\big]e_{1\tss l}(v)
-f_{k1}(v)\big[t_{1j}(u),t_{11}(v)\big]e_{1\tss l}(v)+f_{k1}(v)\big[t_{1j}(u),t_{1\tss l}(v)\big].
\end{multline}
Hence, calculating the commutators by \eqref{defrel}, we come to the relation
\begin{multline}
\non
\big[t_{1j}(u), f_{k\tss 1}(v)h_{1}(v)e_{1\tss l}(v)\big]
=\frac{1}{u-v}\big(t_{kj}(u)h_{1}(v)e_{1\tss l}(v)-t_{kj}(v)h_{1}(u)e_{1\tss l}(v)\big)\\[0.4em]
-\frac{1}{u-v}\big(f_{k1}(v)h_{1}(u)e_{1j}(u)
h_{1}(v)e_{1\tss l}(v)-f_{k1}(v)h_{1}(v)e_{1j}(v)h_{1}(u)e_{1\tss l}(v)\big)\\[0.4em]
+\frac{1}{u-v}\big(f_{k1}(v)t_{1j}(u)
h_{1}(v)e_{1\tss l}(v)-f_{k1}(v)t_{1j}(v)h_{1}(u)e_{1\tss l}(u)\big).
\end{multline}
This gives
\ben
\big[t_{1\tss j}(u), {t^{\ts[1]}_{k\tss l}(v)}\big]
=\frac{1}{u-v}\big({t^{\ts[1]}_{kj}(v)}h_{1}(u)e_{1\tss l}(v)
-{t^{\ts[1]}_{kj}(v)}h_{1}(u)e_{1\tss l}(u)\big).
\een
By Corollary \ref{cor:commu}, $t_{11}(u)$ commutes with $t^{\ts[1]}_{kj}(v)$
and so \eqref{emjmtkl} with $m=1$ follows.
\epf

\section{Drinfeld presentation of extended Yangian}
\label{sec:dpey}

Here we will give a Drinfeld presentation for the extended
Yangian $\X(\g_N)$ analogous to that of the Yangian $\Y(\gl_N)$
\cite{bk:pp}; see also \cite[Sec.~3.1]{m:yc}. Isomorphisms between classical Lie algebras
in low ranks lead to corresponding Yangian isomorphisms; see \cite{amr:rp}
and \cite{jl:ib}. We begin by reviewing them in the context of Drinfeld presentations.

\subsection{Low rank isomorphisms}
\label{subsec:lri}

We will follow the notation of \cite[Sec.~3.1]{m:yc} for the Gauss decomposition
of the generator matrix of the Yangian $\Y(\gl_N)$, but use the corresponding
capital letters $H_i(u)$, $E_{ij}(u)$ and $F_{ji}(u)$
for the entries of the respective matrices occurring in the $\Y(\gl_N)$ counterpart
of \eqref{gd}. The next lemmas are implied by
the results of \cite[Sec.~4]{amr:rp}.

\ble\label{lem:lowrank C}
In terms of the Gaussian generators, the isomorphism $\X(\spa_2)\to \Y(\gl_2)$
has the form
\begin{alignat*}{2}
h_{1}(u)&\mapsto H_{1}(u/2),\qquad\quad  & e_{12}(u)&\mapsto E_{12}(u/2),
\\[0.2em]
h_{2}(u)&\mapsto H_{2}(u/2),\qquad\quad & f_{21}(u)&\mapsto F_{21}(u/2).
\end{alignat*}
\ele

\ble\label{lem:lowrank B}
In terms of the Gaussian generators, the isomorphism $\X(\oa_3)\to \Y(\gl_2)$
has the form
\begin{alignat*}{3}
h_{1}(u)&\mapsto H_{1}(2u)H_{1}(2u+1),\qquad
&e_{12}(u)&\mapsto \sqrt{2}\tss E_{12}(2u+1),\qquad
&f_{21}(u)&\mapsto \sqrt{2}\tss F_{21}(2u+1),\\[0.3em]
h_{2}(u)&\mapsto H_{1}(2u)H_{2}(2u+1),\qquad
&e_{23}(u)&\mapsto -\sqrt{2}\tss E_{12}(2u),\qquad
&f_{32}(u)&\mapsto -\sqrt{2}\tss F_{21}(2u),\\[0.3em]
h_3(u)&\mapsto H_{2}(2u)H_{2}(2u+1),\qquad
&e_{13}(u)&\mapsto -E_{12}(2u+1)^2,\qquad
&f_{31}(u)&\mapsto -F_{21}(2u+1)^2.
\end{alignat*}
\ele

\ble\label{lem:lowrank D}
 Using the
embedding $\X(\oa_4)\hra \Y(\gl_2)\ot \Y(\gl_2)$,
the correspondence in terms of the Gaussian generators between $\X(\oa_4)$ and $\Y(\gl_2)\ot \Y(\gl_2)$ is given by
\begin{alignat*}{2}
h_{1}(u)&\mapsto H_{1}(u)H'_{1}(u),\qquad\quad h_{2'}(u)&&\mapsto H_{2}(u)H'_{1}(u),
\\[0.2em]
h_{2}(u)&\mapsto H_{1}(u)H'_{2}(u),\qquad\quad
h_{1'}(u)&&\mapsto H_{2}(u)H'_{2}(u),
\end{alignat*}
together with
\begin{alignat*}{2}
e_{12}(u)&\mapsto E'_{12}(u),\qquad\quad & e_{12'}(u)&\mapsto E_{12}(u),
\\[0.2em]
e_{11'}(u)&\mapsto -E_{12}(u)E'_{12}(u),\qquad\quad &e_{22'}(u)&\mapsto 0,\\[0.2em]
e_{21'}(u)&\mapsto -E_{12}(u),\qquad\quad &e_{2'1'}(u)&\mapsto -E'_{12}(u),
\end{alignat*}
and
\begin{alignat*}{2}
f_{21}(u)&\mapsto F'_{21}(u),\qquad\quad & f_{2'1}(u)&\mapsto F_{21}(u),
\\[0.2em]
f_{1'1}(u)&\mapsto -F_{21}(u)F'_{21}(u),\qquad\quad &f_{2'2}(u)&\mapsto 0,\\[0.2em]
f_{1'2}(u)&\mapsto -F_{21}(u),\qquad\quad &f_{1'2'}(u)&\mapsto -F'_{21}(u),
\end{alignat*}
where $H'_1(u)$, $H'_2(u)$, $E'_{12}(u)$ and $F'_{21}(u)$
denote the Gaussian generators of the second copy of $\Y(\gl_2)$
in the tensor product.
\ele

It will be convenient to use a uniform root notation for all three cases, so we will
assume that the simple roots of $\g_N$ are $\al_1,\dots,\al_n$ with
\beql{ali}
\al_i=\ep_i-\ep_{i+1},\qquad i=1,\dots,n-1,
\eeq
and
\beql{aln}
\al_n=\begin{cases}
\ep_n
\qquad&\text{for}\quad \oa_{2n+1}\\
2\tss \ep_n
\qquad&\text{for}\quad \spa_{2n}\\
\ep_{n-1}+\ep_n
\qquad&\text{for}\quad \oa_{2n},
\end{cases}
\eeq
where $\ep_1,\dots,\ep_n$ is an orthonormal basis of an Euclidian space
with the bilinear form $(.\ts,.)$.

\bpr\label{pro:rellowrankBCD}
We have the relations in $\X(\g_N)$,
\ben
\bal
\big[h_n(u),e_{n}(v)\big]&=-(\ep_n, \al_n)\frac{h_{n}(u)\big(e_{n}(u)-e_{n}(v)\big)}{u-v},\\[0.3em]
\big[h_n(u),f_{n}(v)\big]&=(\ep_n, \al_n)\frac{\big(f_{n}(u)-f_{n}(v)\big)h_{n}(u)}{u-v},\\[0.3em]
\big[e_{n}(u),e_{n}(v)\big]&=
\frac{(\al_n,\al_n)}{2}\frac{\big(e_{n}(u)-e_{n}(v)\big)^2}{u-v},\\[0.3em]
\big[f_{n}(u),f_{n}(v)\big]&=
-\frac{(\al_n,\al_n)}{2}\frac{\big(f_{n}(u)-f_{n}(v)\big)^2}{u-v},\\[0.3em]
\big[e_{n}(u),f_{n}(v)\big]&=\frac{k_{n}(u)-k_{n}(v)}{u-v}.
\eal
\een
Moreover, for $\g_N=\oa_{2n+1}$ we have
\begin{align}\label{Bhn+1en}
[h_{n+1}(u),e_{n}(v)]&=\frac{1}{2(u-v)}h_{n+1}(u)\big(e_{n}(u)-e_{n}(v)\big)\\
{}&-\frac{1}{2(u-v-1)}\big(e_{n}(u-1)-e_{n}(v)\big)h_{n+1}(u)
\non
\end{align}
and
\begin{align}\label{Bhn+1fn}
[h_{n+1}(u),f_{n}(v)]={}&-\frac{1}{2(u-v)}h_{n+1}(u)\big(f_{n}(u)-f_{n}(v)\big)\\
{}&+\frac{1}{2(u-v-1)}\big(f_{n}(u-1)-f_{n}(v)\big)h_{n+1}(u),
\non
\end{align}
whereas for $\g_N=\spa_{2n}$ and $\oa_{2n}$ we have
\ben
\big[h_{n+1}(u),e_{n}(v)\big]=(\ep_n, \al_n)\tss\frac{h_{n+1}(u)\big(e_{n}(u)-e_{n}(v)\big)}{u-v}
\een
and
\ben
\big[h_{n+1}(u),f_{n}(v)\big]=-(\ep_n, \al_n)\tss \frac{\big(f_{n}(u)-f_{n}(v)\big)h_{n+1}(u)}{u-v}.
\een
\epr

\bpf
By Corollary~\ref{cor:guass-embed}, the subalgebra $\X^{[n-1]}(\g_N)$ of
$\X(\g_N)$ is isomorphic to $\X(\oa_3)$ and $\X(\spa_2)$
in types $B$ and $C$, respectively, while
the subalgebra of $\X^{[n-2]}(\oa_{2n})$
is isomorphic to $\X(\oa_4)$. Hence the relations are implied by
Lemmas~\ref{lem:lowrank C}, \ref{lem:lowrank B}, \ref{lem:lowrank D},
and the Drinfeld presentation of the Yangian $\Y(\gl_2)$; see \cite[Sec.~3.1]{m:yc}.
For instance, to verify the first relation in type $B$, use Lemma~\ref{lem:lowrank B}
to get
\ben
\bal
\big[h_1(u)&,e_1(v)\big]=\big[H_{1}(2u)H_{1}(2u+1),\sqrt{2}\tss E_{12}(2v+1)\big]\\[0.3em]
{}&=\big[H_{1}(2u),\sqrt{2}\tss E_{12}(2v+1)\big]H_{1}(2u+1)+
H_{1}(2u)\big[H_{1}(2u+1),\sqrt{2}\tss E_{12}(2v+1)\big].
\eal
\een
Applying the commutator formula
\ben
(u-v)\big[H_1(u),E_{12}(v)\big]=-H_{1}(u)\big(E_{12}(u)-E_{12}(v)\big)
\een
from \cite[Lemma~3.1.1]{m:yc}, we bring the right hand side to the required form.
\epf

\subsection{Type {\sl A} relations}
\label{subsec:tar}

Since the extended Yangian $\X(\g_{N})$ contains a subalgebra isomorphic
to the Yangian $\Y(\gl_n)$, some relations between the Gaussian generators of
$\X(\g_{N})$ can be obtained from those of the Drinfeld presentation of $\Y(\gl_n)$.
We record them in the next proposition, where we use
the generating functions
\beql{efcirc}
e^\circ_i(u)=\sum_{r=2}^{\infty}e_i^{(r)}u^{-r}\Fand
f^\circ_i(u)=\sum_{r=2}^{\infty}f_i^{(r)}u^{-r}.
\eeq

\bpr\label{pro:Arelations}
The following relations hold in
$\X(\g_{N})$, with the conditions on the indices $1\leqslant i,j\leqslant n-1$
and $1\leqslant k,l \leqslant n$,
\begin{align}
\label{hihj}
\big[h_k(u),h_l(v)\big]&=0, \\
\label{eifj}
\big[e_i(u),f_j(v)\big]&=\delta_{i\tss j}\ts\frac{k_i(u)-k_i(v)}{u-v},\\
\label{hkej}
\big[h_k(u),e_j(v)\big]&=-(\ep_k,\al_j)\ts
\frac{h_k(u)\tss\big(e_j(u)-e_j(v)\big)}{u-v}, \\[0.4em]
\label{hkfj}
\big[h_k(u),f_j(v)\big]&=(\ep_k,\al_j)\ts
\frac{\big(f_j(u)-f_j(v)\big)\tss h_k(u)}{u-v}.
\end{align}
Moreover, for $1\leqslant i\leqslant n-2$ we have
\ben
\bal
u\big[e^{\circ}_i(u),e_{i+1}(v)\big]-v\big[e_i(u),e^{\circ}_{i+1}(v)\big]
&=e_{i}(u)\tss e_{i+1}(v),\\[0.4em]
u\big[f^{\circ}_i(u),f_{i+1}(v)\big]-v\big[f_i(u),f^{\circ}_{i+1}(v)\big]
&=-f_{i+1}(v)\tss f_{i}(u),
\eal
\een
and for $1\leqslant i\leqslant n-1$ we have
\ben
\bal
&\big[e_i(u),e_{i}(v)\big]=\frac{(e_{i}(u)-e_{i}(v))^2}{u-v}
\Fand\\[0.4em]
&\big[f_i(u),f_{i}(v)\big]=-\frac{(f_{i}(u)-f_{i}(v))^2}{u-v}.
\eal
\een
For $1\leqslant i,j\leqslant n-1$ we have
\ben
\big[e_i(u),e_j(v)\big]=0\Fand
\big[f_i(u),f_j(v)\big]=0, \qquad \text{if} \quad(\al_i,\al_j)=0,
\een
whereas
\ben
\bal
\big[e_i(u),[e_i(v),e_j(w)]\big]+\big[e_i(v),[e_i(u),e_j(w)]\big]&=0\Fand\\[0.4em]
\big[f_i(u),[f_i(v),f_j(w)]\big]+\big[f_i(v),[f_i(u),f_j(w)]\big]&=0\qquad \text{if} \quad |i-j|=1.
\eal
\een
\epr

\bpf
The coefficients of the series $t_{ij}(u)$ with $i,j\in\{1,\dots,n\}$
generate a subalgebra of $\X(\g_{N})$ isomorphic to the Yangian $\Y(\gl_n)$.
Hence, the upper left $n\times n$ submatrices of the
matrices $F(u)$, $H(u)$ and $E(u)$
defined by
the Gauss decomposition \eqref{gd} are
given by the same formulas as the corresponding
elements of $\Y(\gl_n)$. Therefore they satisfy the relations as
described in \cite[Section~5]{bk:pp}; see also \cite[Section~3.1]{m:yc}.
\epf

We point out two useful consequences of
\eqref{hkej} and \eqref{hkfj}:
\beql{usefulf1}
h_{i}(u)e_{i}(u)=e_{i}(u-1)h_{i}(u)\Fand
h_{i}(u)f_{i}(u-1)=f_{i}(u)h_{i}(u),
\eeq
which hold for all $i=1,\ldots, n-1$.

Another subalgebra of $\X(\g_N)$ isomorphic to the Yangian $\Y(\gl_n)$
is generated by the coefficients of the series $t_{ij}(u)$ with
$n'\leqslant i,j\leqslant 1'$. This corresponds to the
lower right $n\times n$ submatrix of $T(u)$.
It is known that the map
\beql{autoins}
\vs:T(u)\mapsto T(-u)^{-1}
\eeq
defines an automorphism of $\X(\g_N)$; see \cite[Sec.~2]{amr:rp}.
Hence, the subalgebra of $\X(\g_N)$ generated by
the coefficients of the series $\vs\big(t_{ij}(u)\big)$ with $n'\leqslant i,j\leqslant 1'$
is isomorphic to the Yangian $\Y(\gl_n)$.
On the other hand, by inverting the matrices in the Gauss decomposition \eqref{gd}
we get
\beql{gdinv}
T(u)^{-1}=E(u)^{-1}\tss H(u)^{-1}\tss F(u)^{-1}.
\eeq
Since the upper triangular matrix $E(u)^{-1}$ appears on the left and
the lower triangular matrix $F(u)^{-1}$ appears on the right,
the lower right $n\times n$ submatrix of $T(u)^{-1}$ will only involve
the corresponding submatrices of $E(u)^{-1}$, $H(u)^{-1}$ and $F(u)^{-1}$.
Regarding this lower right submatrix of $T(-u)^{-1}$ as the generator matrix for
$\Y(\gl_n)$, apply now the automorphism of the Yangian $\Y(\gl_n)$
defined by the same formula \eqref{autoins}. We can thus conclude that the product
of the matrices
\ben
\begin{bmatrix}
1&0&\dots&0\ts\\
f_{(n-1)'\ts n'}(u)&1&\dots&0\\
\vdots&\vdots&\ddots&\vdots\\
f_{1'\tss n'}(u)&\dots&\dots&1
\end{bmatrix}
\begin{bmatrix}
h_{n'}(u)&0&\dots&0\ts\\
0&\ddots&&\vdots\\
\vdots&&\ddots&\vdots\\
0&\dots&\dots&h_{1'}(u)
\end{bmatrix}
\begin{bmatrix}
\ts1&e_{n'\ts (n-1)'}(u)&\dots&e_{n'\tss 1'}(u)\ts\\
\ts0&1&\dots&\dots\\
\vdots&\vdots&\ddots&\vdots\\
0&0&\dots&1
\end{bmatrix}
\een
yields the Gauss decomposition of the generator matrix for $\Y(\gl_n)$.
So we have derived the following.

\bpr\label{prop:Areldual}
All relations in $\X(\g_{N})$ given in Proposition~\ref{pro:Arelations}
remain valid under the replacement of the indices of all series
by $i\mapsto (n-i+1)'$ for $1\leqslant i\leqslant n$.
\qed
\epr

\subsection{Central elements in terms of Gaussian generators}
\label{subsec:cegd}

We note some relations implied by Corollary~\ref{cor:guass-embed}
and low rank isomorphisms pointed out in Section~\ref{subsec:lri}.
Write relation \eqref{zcenter} in the form
\beql{ztinv}
T\pr(u+\ka)=z_N(u)\tss T(u)^{-1}.
\eeq
By the Gauss decomposition \eqref{gd} we have $t_{11}(u)=h_1(u)$
so that using \eqref{gdinv}
and taking the $(N,N)$-entry on both sides of \eqref{ztinv} we get
\beql{gaussTrans}
h_1(u+\ka)=z_N(u)\tss h_{1'}(u)^{-1}.
\eeq
Taking $N=2n+1$ and
applying this relation to the subalgebra $\X^{[n-1]}(\oa_N)\cong\X(\oa_3)$
(see Corollary~\ref{cor:guass-embed}), we get
\ben
z^{[n-1]}_N(u)=h_n(u+\ka-n+1)\tss h_{n'}(u).
\een
Lemma~\ref{lem:lowrank B} allows us to bring this to the form
\beql{Z1B}
z^{[n-1]}_N(u)=
h_{n}\big(u+1/2\big)\tss h_{n}\big(u-1/2\big)^{-1}\tss h_{n+1}(u)\tss h_{n+1}\big(u-1/2\big).
\eeq

Similarly, the application of \eqref{gaussTrans} to the subalgebra
$\X^{[n-1]}(\spa_N)\cong\X(\spa_2)$ with $N=2n$ gives
\beql{Z1C}
z^{[n-1]}_N(u)=h_{n}(u+\ka-n+1)\tss h_{n+1}(u)=h_{n}(u+2)\tss h_{n+1}(u).
\eeq

If $\g_N=\oa_N$ with $N=2n$ we apply \eqref{gaussTrans} to the subalgebra
$\X^{[n-2]}(\oa_N)\cong\X(\oa_4)$ to get
\ben
z^{[n-2]}_N(u)=h_{n-1}(u+{\ka-n+2})\tss h_{(n-1)'}(u)=h_{n-1}(u+1)\tss h_{(n-1)'}(u).
\een
We have $h_{n-1}(u)\tss h_{(n-1)'}(u)=h_n(u)\tss h_{n'}(u)$ by Lemma~\ref{lem:lowrank D}
and so
\beql{Z2}
z^{[n-2]}_N(u)=h_{n-1}(u+1)\tss h_{n-1}(u)^{-1}\tss h_n(u)\tss h_{n'}(u).
\eeq

We will need some symmetry properties for the entries
of the matrices in the Gauss decomposition \eqref{gd}.

\bpr\label{prop:eiei'}
The following relations hold in $\X(\g_{N})$,
\beql{eiei'}
e_{(i+1)'\ts i'}(u)=-e_{i}(u+\ka-i)\Fand
f_{i'\ts (i+1)'}(u)=-f_{i}(u+\ka-i)
\eeq
for $i=1,\dots,n-1$.
Furthermore, if $\g_N=\oa_{2n+1}$ then we have
\ben
e_{n+1\ts n+2}(u)=-e_{n}\big(u-1/2\big)\Fand f_{n+2\ts n+1}(u)=-f_{n}\big(u-1/2\big),
\een
and if $\g_N=\oa_{2n}$ then
\ben
e_{n\ts (n-1)'}(u)=-e_{n}(u)\Fand f_{(n-1)'\ts n}(u)=-f_{n}(u).
\een
\epr

\bpf
Suppose that $1\leqslant i\leqslant n-1$.
By Corollary~\ref{cor:guass-embed}, we have the following
consequence of relation \eqref{subunitary},
\beql{ztinvi}
T^{[i-1]}(u)^{-1}\tss z^{[i-1]}_N(u)=T^{[i-1]\tss\prime}(u+\ka^{[i-1]}).
\eeq
As with the above derivation of \eqref{gaussTrans},
take the $(i\pr,i\pr)$-entry on both sides of \eqref{ztinvi} to get
\beql{hiprin}
h_{i'}(u)^{-1}\tss z^{[i-1]}_N(u)=h_{i}(u+\ka-i+1).
\eeq
Similarly, by taking the $\big((i+1)\pr,i\pr\big)$-entry in \eqref{ztinvi} we obtain
\ben
-e_{(i+1)'\tss i'}(u)\tss h_{i'}(u)^{-1}\tss z^{[i-1]}_N(u)=h_{i}(u+\ka-i+1)\tss e_{i}(u+\ka-i+1).
\een
Together with \eqref{hiprin} this gives
\ben
-e_{(i+1)'\tss i'}(u)=h_{i}(u+\ka-i+1)\tss e_{i}(u+\ka-i+1)\tss h_{i}(u+\ka-i+1)^{-1}.
\een
The first relation in \eqref{eiei'} now follows from
\eqref{usefulf1} and a similar argument verifies the
second.
The additional relations in types $B$ and $D$
follow from Lemmas~\ref{lem:lowrank B} and \ref{lem:lowrank D}, respectively.
\epf

We are now in a position to give explicit formulas for the series $z_N(u)$
in terms of the Gaussian generators $h_i(u)$. Recall that by Proposition~\ref{pro:Arelations}
the coefficients of the series $h_i(u)$ pairwise commute for $i=1,\dots,n$; see \eqref{hihj}.

\bth\label{thm:Center}
We have the identities in $\X(\g_{N})${\rm :}
\ben
z_N(u)=\prod_{i=1}^n h_i(u+\ka-i)^{-1}\tss \prod_{i=1}^n h_i(u+\ka-i+1)\cdot
h_{n+1}(u)\tss h_{n+1}\big(u-1/2\big)
\een
for $\g_N=\oa_{2n+1}$,
\ben
z_N(u)=\prod_{i=1}^{n-1} h_i(u+\ka-i)^{-1}\tss \prod_{i=1}^n h_i(u+\ka-i+1)\cdot
h_{n+1}(u)
\een
for $\g_N=\oa_{2n}$ and $\g_N=\spa_{2n}$.
\eth

\bpf
Take the $(2\pr,2\pr)$ entries on both sides of \eqref{ztinv}.
Expressing the entries of the matrices $T\pr(u+\ka)$ and $T(u)^{-1}$
in terms of the Gaussian generators from \eqref{gd} and \eqref{gdinv},
we get
\ben
h_{2}(u+\ka)+f_{1}(u+\ka)\tss h_{1}(u+\ka)\tss e_{1}(u+\ka)=
z_N(u)\big(h_{2'}(u)^{-1}+e_{2'\tss 1'}(u)\tss h_{1'}(u)^{-1}f_{1'\tss 2'}(u)\big).
\een
Since $z_N(u)$ is central in $\X(\g_{N})$, using \eqref{gaussTrans} we can rewrite
this as
\ben
h_{2'}(u)^{-1}z_N(u)=h_{2}(u+\ka)+f_{1}(u+\ka)\tss h_{1}(u+\ka)\tss e_{1}(u+\ka)
-e_{2'\tss 1'}(u)\tss h_{1}(u+\ka)\tss f_{1'\tss 2'}(u).
\een
Now apply \eqref{eiei'} to get
\begin{multline}
\non
h_{2'}(u)^{-1}z_N(u)=h_{2}(u+\ka)
+f_{2\tss 1}(u+\ka)\tss h_{1}(u+\ka)\tss e_{1\tss 2}(u+\ka)\\[0.3em]
{}-e_{1\tss 2}(u+\ka-1)\tss h_{1}(u+\ka)\tss f_{2\tss 1}(u+\ka-1).
\end{multline}
Due to \eqref{usefulf1}, this simplifies to
\ben
h_{2'}(u)^{-1}z_N(u)=h_{2}(u+\ka)+\big[f_{2\tss 1}(u+\ka),e_{1\tss 2}(u+\ka-1)\big]\tss h_{1}(u+\ka).
\een
Calculating the commutator by \eqref{eifj}, we bring this relation to the form
\ben
h_{2'}(u)^{-1}z_N(u)=h_{1}(u+\ka-1)^{-1}\tss h_1(u+\ka)\tss h_{2}(u+\ka-1).
\een
Finally, use \eqref{hiprin} with $i=2$ to get the recurrence formula
\ben
z_N(u)=h_1(u+\ka-1)^{-1}\tss h_1(u+\ka)\ts z^{[1]}_N(u).
\een
By Corollary~\ref{cor:guass-embed},
the desired identities now follow from the respective
base cases \eqref{Z1B}, \eqref{Z1C} and \eqref{Z2}.
\epf

\bre\label{rem:compa}
The expansions provided by Theorem~\ref{thm:Center} are analogous to
the multiplicative formula for the central series $z(u)$ for the Yangian $\Y(\gl_N)$
implied by the quantum Liouville formula and a decomposition of the quantum determinant
\cite[Theorem~1.9.5 and Corollary~1.11.8]{m:yc}.
\qed
\ere

\subsection{Relations for Gaussian generators}
\label{subsec:rgd}

In addition to the type $A$ relations in $\X(\g_N)$ described in Section~\ref{subsec:tar},
we will now derive some root system specific relations for each of the types $B$, $C$ and $D$.

\bpr\label{prop:relhn+1en-1}
We have the relations in $\X(\g_N)${\rm :}
\ben
\big[h_{n+1}(u),e_{n-1}(v)\big]=0\Fand \big[h_{n+1}(u),f_{n-1}(v)\big]=0
\een
for $\g_N=\oa_{2n+1}$,
\ben
\bal
\big[h_{n+1}(u),e_{n-1}(v)\big]&=\frac{h_{n+1}(u)\big(e_{n-1}(v)-e_{n-1}(u+2)\big)}{u-v+2}\Fand\\[0.4em]
\big[h_{n+1}(u),f_{n-1}(v)\big]&=-\frac{\big(f_{n-1}(v)-f_{n-1}(u+2)\big)\ts h_{n+1}(u)}{u-v+2}
\eal
\een
for $\g_N=\spa_{2n}$, and
\ben
\bal
\big[h_{n+1}(u),e_{n-1}(v)\big]&=\frac{h_{n+1}(u)\ts\big(e_{n-1}(v)-e_{n-1}(u)\big)}{u-v}\Fand\\[0.4em]
\big[h_{n+1}(u),f_{n-1}(v)\big]&=-\frac{\big(f_{n-1}(v)-f_{n-1}(u)\big)\ts h_{n+1}(u)}{u-v}
\eal
\een
for $\g_N=\oa_{2n}$.
\epr

\bpf
First take $\g_N=\oa_{2n+1}$.
Corollary~\ref{cor:commu} implies that
$h_{n+1}(u)$ commutes with each element of the subalgebra
generated by the $t_{ij}(u)$ with $1\leqslant i,j \leqslant n$
and so the relations follow. Now let $\g_N=\spa_{2n}$ or $\g_N=\oa_{2n}$.
By Corollary~\ref{cor:guass-embed}, the subalgebra
$\X^{[n-2]}(\g_N)$ of $\X(\g_N)$ is isomorphic to $\X(\g_4)$.
Applying Proposition~\ref{prop:Areldual} to this subalgebra,
we get
\ben
\big[h_{n+1}(u),e_{n+1\ts n+2}(v)\big]
=\frac{h_{n+1}(u)\big(e_{n+1\ts n+2}(v)-e_{n+1\ts n+2}(u)\big)}{u-v}
\een
and
\ben
\big[h_{n+1}(u),f_{n+2\ts n+1}(v)\big]
=-\frac{\big(f_{n+2\ts n+1}(v)-f_{n+2\ts n+1}(u)\big)\tss h_{n+1}(u)}{u-v}.
\een
It remains to apply
Proposition~\ref{prop:eiei'}.
\epf

In the next proposition we use the root notation \eqref{ali} and \eqref{aln}.

\bpr\label{pro:relenei}
For $i=1,\dots,n$ in the algebra $\X(\g_{N})$ we have
\begin{align}
\label{hihn+1}
\big[h_i(u),h_{n+1}(v)\big]&=0,\\
\label{hien}
\big[h_i(u),e_{n}(v)\big]&=-(\ep_i, \al_n)
\tss\frac{h_{i}(u)\big(e_{n}(u)-e_{n}(v)\big)}{u-v},\\
\label{hifn}
\big[h_i(u),f_{n}(v)\big]&=(\ep_i, \al_n)
\tss\frac{\big(f_{n}(u)-f_{n}(v)\big)h_{i}(u)}{u-v},
\end{align}
and
\beql{eifn}
\big[e_{i}(u),f_{n}(v)\big]=\big[e_{n}(u),f_{i}(v)\big]=
\delta_{i\tss n}\frac{k_{n}(u)-k_{n}(v)}{u-v}.
\eeq
Moreover, if $i<n-1$ then
\beql{hn+1ei}
\big[h_{n+1}(u),e_{i}(v)\big]=0\Fand
\big[h_{n+1}(u),f_{i}(v)\big]=0.
\eeq
\epr

\bpf
Note that relations \eqref{hien} and \eqref{hifn} for $i=n$ were already
verified in Proposition~\ref{pro:rellowrankBCD}.
Now suppose that $\g_N=\oa_{2n}$. By the defining relations
\eqref{defrel}, the subalgebra generated by the coefficients of the series $t_{ij}(u)$
with $i,j$ running over the set
$J=\{1,\dots,n-1, n+1\}$
is isomorphic to $\Y(\gl_n)$.
Furthermore, Lemma~\ref{lem:lowrank D} implies that
\ben
e_{n\ts n+1}(u)=f_{n+1\ts n}(u)=0.
\een
Hence, we have the Gauss decomposition
\beql{gdj}
T_J(u)=F_J(u)\ts H_J(u)\ts E_J(u),
\eeq
where the subscript $J$ indicates the submatrices in \eqref{gd} whose rows and columns
are labelled by the elements of the set $J$.  Therefore, the Gaussian
generators which occur as the entries of the matrices
$F_J(u)$, $H_J(u)$ and $E_J(u)$ satisfy the type $A$ relations as described
in Proposition~\ref{pro:Arelations}. This completes the proof for
type $D$.

Now let $\g_N=\oa_{2n+1}$ or $\g_N=\spa_{2n}$. Almost all of the relations
\eqref{hien} and \eqref{hifn} with $i<n$, as well as \eqref{eifn}
and \eqref{hn+1ei} with $i<n-1$ follow from Corollary~\ref{cor:commu}.
For instance, \eqref{hien} and \eqref{hifn} are immediate from
the observation that all elements of
the subalgebra of $\X(\g_N)$ generated by $t_{ij}(u)$ with $i,j=1,\dots,n-1$
commute with the subalgebra $\X^{[n-1]}(\g_N)$. In addition, we use
Lemma~\ref{lem:lowrank C} to see that $[h_{n}(u),h_{n+1}(v)]=0$ in type $C$.

Furthermore, the case $i=n$ of \eqref{eifn} was already
pointed out in Proposition~\ref{pro:rellowrankBCD}.
To verify the remaining cases $i=n-1$ of \eqref{eifn}, apply
Lemma~\ref{lem:emjmtkl} with $m=n-1$. Since
$t^{\ts[n-1]}_{n+1\ts n}(v)=f_{n}(v)h_{n}(v)$,
relation \eqref{emjmtkl} gives
\ben
\big[e_{n-1}(u),f_{n}(v)h_{n}(v)\big]
=\frac{1}{u-v}\tss f_{n}(v)h_{n}(v)\tss \big(e_{n-1}(v)-e_{n-1}(u)\big).
\een
On the other hand,
\ben
\big[e_{n-1}(u),f_{n}(v)h_{n}(v)\big]=\big[e_{n-1}(u),f_{n}(v)\big]\tss h_{n}(v)
+f_{n}(v)\tss \big[e_{n-1}(u),h_{n}(v)\big],
\een
whereas
\beql{enhnuv}
\big[e_{n-1}(u),h_{n}(v)\big]=\frac{1}{u-v}\tss h_{n}(v)\tss \big(e_{n-1}(v)-e_{n-1}(u)\big)
\eeq
by Proposition~\ref{pro:Arelations}.
This gives $\big[e_{n-1}(u),f_n(v)\big]=0$. The other case of \eqref{eifn}
is verified in the same way.
\epf

\ble\label{Lemma:en-1en}
In the algebra $\X(\oa_{2n+1})$ we have
\begin{align}\label{Ben-1en}
&\big[e_{n-1}(u), e_{n}(v)\big]
=\frac{e_{n-1\tss n+1}(v)-e_{n-1\tss n+1}(u)
-e_{n-1}(v)e_{n}(v)+e_{n-1}(u)e_{n}(v)}{u-v},\\
\label{Bfn-1fn}
&\big[f_{n-1}(u), f_{n}(v)\big]
=\frac{f_{n+1\tss n-1}(u)-f_{n+1\tss n-1}(v)
-f_{n-1}(u)f_{n}(v)+f_{n-1}(v)f_{n}(v)}{u-v}.
\end{align}
In the algebra $\X(\oa_{2n})$ we have
\beql{Den-1en}
\big[e_{n-1}(u), e_{n}(v)\big]=0,
\qquad\quad \big[f_{n-1}(u), f_{n}(v)\big]=0,
\eeq
and
\begin{align}\label{Den-2en}
&\big[e_{n-2}(u), e_{n}(v)\big]
=\frac{e_{n-2\tss n+1}(v)-e_{n-2\tss n+1}(u)
-e_{n-2}(v)e_{n}(v)+e_{n-2}(u)e_{n}(v)}{u-v},\\
\label{Dfn-2fn}
&\big[f_{n-2}(u), f_{n}(v)\big]
=\frac{f_{n+1\tss n-2}(u)-f_{n+1\tss n-2}(v)
-f_{n-2}(u)f_{n}(v)+f_{n-2}(v)f_{n}(v)}{u-v}.
\end{align}
In the algebra $\X(\spa_{2n})$ we have
\begin{align}\label{Cen-1en}
\big[e_{n-1}(u), e_{n}(v)\big]&
=\frac{2\tss\big(e_{n-1\tss n+1}(v)-e_{n-1\ts n+1}(u)
-e_{n-1}(v)e_{n}(v)+e_{n-1}(u)e_{n}(v)\big)}{u-v},\\
\label{Cfn-1fn}
\big[f_{n-1}(u), f_{n}(v)\big]&
=\frac{2\tss\big(f_{n+1\tss n-1}(u)
-f_{n+1\ts n-1}(v)-f_{n-1}(u)f_{n}(v)+f_{n-1}(v)f_{n}(v)\big)}{u-v}.
\end{align}
\ele

\bpf
By Lemma~\ref{lem:emjmtkl}, in $\X(\oa_{2n+1})$ we have
\ben
\big[e_{n-1\tss n}(u), {t^{\ts[n-1]}_{n\tss n+1}}(v)\big]
=\frac{1}{u-v}{t^{\ts[n-1]}_{n\tss n}(v)}\big(e_{n-1\tss n+1}(v)-e_{n-1\tss n+1}(u)\big).
\een
Using the Gauss decomposition for $T^{\tss[n-1]}(u)$, we can write the left hand
side as
\ben
\big[e_{n-1\tss n}(u), h_n(v)\tss e_{n\tss n+1}(v)\big]
=\big[e_{n-1\tss n}(u), h_n(v)\big]\tss e_{n\tss n+1}(v)
+h_{n}(v)\tss \big[e_{n-1\tss n}(u),e_{n\tss n+1}(v)\big].
\een
Now observe that ${t^{\ts[n-1]}_{n\tss n}(v)}=h_{n}(v)$ so that
\eqref{Ben-1en} follows by the application of \eqref{enhnuv}.
A similar argument yields \eqref{Bfn-1fn}.

Now turn to the case $\g_N=\oa_{2n}$.
Relation \eqref{Den-1en} follows from Lemma~\ref{lem:lowrank D}.
As we pointed out in the proof of Proposition~\ref{pro:relenei},
the Gaussian
generators which occur as the entries of the matrices
$F_J(u)$, $H_J(u)$ and $E_J(u)$ in \eqref{gdj}
satisfy the type $A$ relations as described
in Proposition~\ref{pro:Arelations}. Therefore,
\eqref{Den-2en} and \eqref{Dfn-2fn}
follow from the Drinfeld presentation of the Yangian $\Y(\gl_n)$;
see~\cite[Lemma~3.1.2]{m:yc}.

To prove \eqref{Cen-1en} and \eqref{Cfn-1fn}, note that
by Corollary~\ref{cor:guass-embed}, the subalgebra
$\X^{[n-2]}(\spa_{2n})$ is isomorphic to $\X(\spa_4)$.
Hence, we may take $n=2$ and will work with this case throughout the rest of the argument.

The defining relations \eqref{defrel} for $\X(\spa_4)$ give
\begin{align}\label{t12t23}
\big[t_{12}(u),t_{23}(v)\big]&=\frac{1}{u-v}\Big(t_{22}(u)t_{13}(v)-t_{22}(v)t_{13}(u)\Big)\\[0.3em]
&+\frac{1}{u-v-3}\Big(t_{24}(v)t_{11}(u)+t_{23}(v)t_{12}(u)-t_{22}(v)t_{13}(u)-t_{21}(v)t_{14}(u)\Big).
\non
\end{align}
Applying the Gauss decomposition \eqref{gd},
for the left hand side of \eqref{t12t23} we can write
\begin{multline}
\big[t_{12}(u),t_{23}(v)\big]
=\big[h_1(u)e_{12}(u),h_2(v)e_{23}(v)+f_{21}(v)h_{1}(v)e_{13}(v)\big]=\\[0.4em]
=h_1(u)\big[e_{12}(u),h_2(v)e_{23}(v)\big]+\big[h_1(u),h_2(v)e_{23}(v)\big]e_{12}(u)\\[0.4em]
+\big[h_1(u)e_{12}(u),f_{21}(v)h_1(v)\big]e_{13}(v)
+f_{21}(v)h_1(v)\big[h_{1}(u)e_{12}(u),e_{13}(v)\big].
\non
\end{multline}
Corollary~\eqref{cor:commu} implies that
$\big[h_1(u),h_2(v)e_{23}(v)\big]=0$. Furthermore, by \eqref{defrel}
\begin{multline}
\big[h_1(u)e_{12}(u),f_{21}(v)h_{1}(v)\big]=\big[t_{12}(u),t_{21}(v)\big]
=\frac{1}{u-v}\big(t_{22}(u)t_{11}(v)-t_{22}(v)t_{11}(u)\big)\\[0.4em]
=\frac{1}{u-v}\Big(h_2(u)h_1(v)-h_2(v)h_{1}(u)
+f_{21}(u)h_1(u)e_{12}(u)h_1(v)-f_{21}(v)h_1(v)e_{12}(v)h_1(u)\Big).
\non
\end{multline}
Similarly, we have
\ben
\big[h_{1}(u)e_{12}(u),e_{13}(v)\big]
=\big[t_{12}(u),t_{11}(v)^{-1}t_{13}(v)\big],
\een
and \eqref{defrel} gives
\ben
\big[t_{12}(u),t_{11}(v)^{-1}\big]
=\frac{1}{u-v}t_{11}(v)^{-1}\Big(t_{12}(v)t_{11}(u)-t_{12}(u)t_{11}(v)\Big)t_{11}(v)^{-1}
\een
together with
\ben
\bal
\big[t_{12}(u),t_{13}(v)\big]&=
\frac{1}{u-v}\big(t_{12}(u)t_{13}(v)-t_{12}(v)t_{13}(u)\big)\\
{}&+\frac{1}{u-v-3}\big(t_{14}(v)t_{11}(u)+t_{13}(v)t_{12}(u)
-t_{12}(v)t_{13}(u)-t_{11}(v)t_{14}(u)\big).
\eal
\een
Writing the resulting expression back in terms of the Gaussian generators
and applying \eqref{hkej} with $k=2$ and $j=1$
we find that the left hand side of \eqref{t12t23} equals
\ben
\begin{aligned}
&h_1(u)h_2(v)[e_{12}(u),e_{23}(v)]+\frac{1}{u-v}
\Big(h_1(u)h_2(v)e_{12}(v)e_{23}(v)-h_1(u)h_2(v)e_{12}(u)e_{23}(v)\Big)\\[0.4em]
&+\frac{1}{u-v}\Big(h_2(u)h_{1}(v)e_{13}(v)-h_2(v)h_{1}(u)e_{13}(v)\\
&\qquad\qquad\qquad\qquad\qquad
 {}+f_{21}(u)h_1(u)e_{12}(u)h_1(v)e_{13}(v)-f_{21}(v)h_{1}(v)e_{12}(v)h_1(u)e_{13}(u)\Big)+\\[0.4em]
&\frac{1}{u-v-3}\ts f_{21}(v)h_1(v)\Big(\hspace{-1pt}e_{14}(v)h_{1}(u)+e_{13}(v)h_1(u)e_{12}(u)
-e_{12}(v)h_1(u)e_{13}(u)-h_{1}(u)e_{14}(u)\hspace{-1pt}\Big).
\end{aligned}
\een
As a next step, write the right hand side of \eqref{t12t23} in terms of
the Gaussian generators. Cancelling common terms on both sides,
we bring the relation to the form
\ben
\bal
h_1(u)h_2(v)\tss &\big[e_{12}(u),e_{23}(v)\big]\\[0.4em]
&=\frac{1}{u-v}\ts h_1(u)h_{2}(v)
\Big(e_{13}(v)-e_{13}(u)-e_{12}(v)e_{23}(v)+e_{12}(u)e_{23}(v)\Big)\\[0.4em]
&+\frac{1}{u-v-3}\ts h_2(v)\Big(e_{24}(v)h_{1}(u)+h_{1}(u)e_{23}(v)e_{12}(u)-h_{1}(u)e_{13}(u)\Big),
\eal
\een
where we also used the property $\big[h_{1}(u),e_{23}(v)\big]=0$ implied
by Corollary~\ref{cor:commu}.
Since the series $h_1(u)$ and $h_2(v)$ are invertible, we thus get
\begin{align}
\label{e12e23}
\big[e_{12}(u),e_{23}(v)\big]&=\frac{1}{u-v}
\Big(e_{13}(v)-e_{13}(u)-e_{12}(v)e_{23}(v)+e_{12}(u)e_{23}(v)\Big)\\[0.4em]
&+\frac{1}{u-v-3}\Big(h_1(u)^{-1}e_{24}(v)h_{1}(u)+e_{23}(v)e_{12}(u)-e_{13}(u)\Big).
\non
\end{align}

Now we need an expression for the commutator $\big[h_{1}(u),e_{24}(v)\big]$
which is obtained by a calculation similar to the above derivation of \eqref{e12e23}.
Namely, we begin with the following analogue of \eqref{t12t23},
\begin{align}
\label{t11t24}
\big[t_{11}(u),t_{24}(v)\big]&=\frac{1}{u-v}\Big(t_{21}(u)t_{14}(v)-t_{21}(v)t_{14}(u)\Big)\\[0.3em]
&+\frac{1}{u-v-3}\Big(t_{24}(v)t_{11}(u)+t_{23}(v)t_{12}(u)-t_{22}(v)t_{13}(u)-t_{21}(v)t_{14}(u)\Big).
\non
\end{align}
Then write the left hand side in terms of the Gaussian generators so that it equals
\ben
h_2(v)\big[h_{1}(u),e_{24}(v)\big]+\big[h_1(u),f_{21}(v)\big]
h_{1}(v)e_{14}(v)+f_{21}(v)\big[h_{1}(u),h_1(v)e_{14}(v)\big].
\een
Furthermore, use \eqref{hkfj} with $k=j=1$
and expand
\ben
\big[h_1(u),h_1(v)e_{14}(v)\big]=\big[t_{11}(u),t_{14}(v)\big]
\een
by the defining relations \eqref{defrel}.
Now writing the resulting expressions on both sides of \eqref{t11t24}
in terms of the Gaussian generators and simplifying as with the
derivation of \eqref{e12e23},
we come to the desired commutation relation
\ben
\big[h_{1}(u),e_{24}(v)\big]=\frac{1}{u-v-3}
\Big(e_{24}(v)h_1(u)+h_1(u)e_{23}(v)e_{12}(u)-h_1(u)e_{13}(u)\Big).
\een
Applying it to \eqref{e12e23}, we can transform the latter as
\ben
\bal
\big[e_{12}(u),e_{23}(v)\big]&=\frac{1}{u-v}
\Big(e_{13}(v)-e_{13}(u)-e_{12}(v)e_{23}(v)+e_{12}(u)e_{23}(v)\Big)\\[0.3em]
&+\frac{1}{u-v-2}\Big(e_{24}(v)+e_{23}(v)e_{12}(u)-e_{13}(u)\Big).
\eal
\een
By rearranging the terms, write it in an equivalent form,
\ben
\bal
\big[e_{12}(u),e_{23}(v)\big]&=\frac{2}{u-v}
\Big(e_{13}(v)-e_{13}(u)-e_{12}(v)e_{23}(v)+e_{12}(u)e_{23}(v)\Big)\\[0.3em]
&+\frac{1}{u-v-1}\Big(e_{24}(v)+e_{12}(v)e_{23}(v)-e_{13}(v)\Big).
\eal
\een
Finally,
multiply both sides by $u-v-1$ and set $u=v+1$ to see that the second summand vanishes.
This yields \eqref{Cen-1en}.
Relation \eqref{Cfn-1fn} is verified by a similar argument.
\epf

In the next proposition we use notation \eqref{efcirc}
for generating functions of elements of the extended Yangian $\X(\g_{N})$.

\bpr\label{pro:inrelation}
Suppose that $1\leqslant i\leqslant n-1$.
If $(\al_i,\al_n)=0$ then
\beql{eien=0}
\big[e_i(u),e_n(v)\big]=0\Fand
\big[f_i(u),f_n(v)\big]=0.
\eeq
If $(\al_i,\al_n)\ne 0$ then
\begin{align}
\label{eien}
u\big[e^{\circ}_i(u),e_{n}(v)\big]-v\big[e_i(u),e^{\circ}_{n}(v)\big]
&=-(\al_i,\al_{n})\tss e_{i}(u)\tss e_{n}(v)\Fand\\[0.4em]
\label{fifn}
u\big[f^{\circ}_i(u),f_{n}(v)\big]-v\big[f_i(u),f^{\circ}_{n}(v)\big]
&=(\al_i,\al_{n})\tss f_{n}(v)\tss f_{i}(u).
\end{align}
\epr

\bpf
If $i\leqslant n-2$ in types $B$ and $C$ or $i\leqslant n-3$ in type $D$, then
\eqref{eien=0} is implied by Corollary~\ref{cor:commu}. If $i=n-1$ in type $D$,
then \eqref{eien=0} follows from Lemma~\ref{lem:lowrank D} with $e_1(u)=e_{12}(u)$
and $e_2(u)=e_{12'}(u)$ by taking into account
Corollary~\ref{cor:guass-embed}.
To get \eqref{eien} and \eqref{fifn}, consider the expressions
for $(u-v)\big[e_i(u),e_{n}(v)\big]$ and $(u-v)\big[f_i(u),f_{n}(v)\big]$
provided by Lemma \ref{Lemma:en-1en} and
take the coefficients of $u^{-r}v^{-s}$ for $r,s\geqslant 1$.
\epf

Note that relations \eqref{eien} and \eqref{fifn} hold in the case $(\al_i,\al_n)=0$ as well;
they are implied by \eqref{eien=0}.

\subsection{Theorem on the Drinfeld presentation}
\label{subsec:tdp}

We will now prove the theorem on the Drinfeld presentation for $\X(\g_N)$.
We use notation \eqref{defkn}, \eqref{defenfn}, \eqref{defknu} and \eqref{efcirc}
for the generating series and the root notation \eqref{ali}
and \eqref{aln}.

\bth\label{thm:dp}
The extended Yangian $\X(\g_{N})$ is generated by
the coefficients of the series
$h_i(u)$ with $i=1,\dots,n+1$, and
$e_i(u)$ and $f_i(u)$ with $i=1,\dots, n$,
subject only to the following sets of relations, where the indices
take all admissible values unless specified otherwise.
We have
\begin{align}
\label{CompleteBhihj}
\big[h_i(u),h_j(v)\big]&=0, \\
\label{CompleteBeifj}
\big[e_i(u),f_j(v)\big]&=\delta_{i\tss j}\ts\frac{k_i(u)-k_i(v)}{u-v}.
\end{align}
For $i\leqslant n$ we have
\begin{align}
\label{CompleteBhiej}
\big[h_i(u),e_j(v)\big]&=-(\ep_i,\al_j)\ts
\frac{h_i(u)\tss\big(e_j(u)-e_j(v)\big)}{u-v},\\[0.4em]
\label{CompleteBhifj}
\big[h_i(u),f_j(v)\big]&=(\ep_i,\al_j)\ts
\frac{\big(f_j(u)-f_j(v)\big)\tss h_i(u)}{u-v}.
\end{align}
For $j\leqslant n-2$ we have
\beql{CompleteBhn+1ej}
\big[h_{n+1}(u),e_j(v)\big]=0,\qquad
\big[h_{n+1}(u),f_j(v)\big]=0.
\eeq
For $\g_N=\oa_{2n+1}$ we have
\begin{align}\label{Bhn+1enthm}
[h_{n+1}(u),e_{n}(v)]&=\frac{1}{2(u-v)}h_{n+1}(u)\big(e_{n}(u)-e_{n}(v)\big)\\
{}&-\frac{1}{2(u-v-1)}\big(e_{n}(u-1)-e_{n}(v)\big)h_{n+1}(u)
\non
\end{align}
and
\begin{align}\label{Bhn+1fnthm}
[h_{n+1}(u),f_{n}(v)]={}&-\frac{1}{2(u-v)}h_{n+1}(u)\big(f_{n}(u)-f_{n}(v)\big)\\
{}&+\frac{1}{2(u-v-1)}\big(f_{n}(u-1)-f_{n}(v)\big)h_{n+1}(u),
\non
\end{align}
whereas for $\g_N=\spa_{2n}$ and $\oa_{2n}$ we have
\beql{Chn+1en}
\big[h_{n+1}(u),e_{n}(v)\big]=(\ep_n, \al_n)\tss\frac{h_{n+1}(u)\big(e_{n}(u)-e_{n}(v)\big)}{u-v}
\eeq
and
\beql{Chn+1fn}
\big[h_{n+1}(u),f_{n}(v)\big]=-(\ep_n, \al_n)\tss \frac{\big(f_{n}(u)-f_{n}(v)\big)h_{n+1}(u)}{u-v}.
\eeq
Moreover, for $\g_N=\oa_{2n+1}$ we have
\beql{completehn+1en-1b}
\big[h_{n+1}(u),e_{n-1}(v)\big]=0\Fand \big[h_{n+1}(u),f_{n-1}(v)\big]=0,
\eeq
for $\g_N=\spa_{2n}$ we have
\begin{align}\label{completehn+1en-1c}
\big[h_{n+1}(u),e_{n-1}(v)\big]&=\frac{h_{n+1}(u)\big(e_{n-1}(v)-e_{n-1}(u+2)\big)}{u-v+2}\Fand\\[0.4em]
\big[h_{n+1}(u),f_{n-1}(v)\big]&=-\frac{\big(f_{n-1}(v)-f_{n-1}(u+2)\big)\ts h_{n+1}(u)}{u-v+2}
\non
\end{align}
and for $\g_N=\oa_{2n}$ we have
\begin{align}\label{completehn+1en-1d}
\big[h_{n+1}(u),e_{n-1}(v)\big]&=\frac{h_{n+1}(u)\ts\big(e_{n-1}(v)-e_{n-1}(u)\big)}{u-v}\Fand\\[0.4em]
\big[h_{n+1}(u),f_{n-1}(v)\big]&=-\frac{\big(f_{n-1}(v)-f_{n-1}(u)\big)\ts h_{n+1}(u)}{u-v}.
\non
\end{align}
In all three cases we have
\begin{align}
\label{eiei}
&\big[e_i(u),e_{i}(v)\big]=\frac{(\al_i,\al_{i})}{2}\frac{(e_{i}(u)-e_{i}(v))^2}{u-v}
\Fand\\[0.4em]
\label{fifi}
&\big[f_i(u),f_{i}(v)\big]=-\frac{(\al_i,\al_{i})}{2}\frac{(f_{i}(u)-f_{i}(v))^2}{u-v}.
\end{align}
Furthermore,
\beql{eiej=0}
\big[e_i(u),e_j(v)\big]=0\Fand
\big[f_i(u),f_j(v)\big]=0 \qquad \text{if}\quad(\al_i,\al_j)=0,
\eeq
whereas for $i\ne j$ we have
\begin{align}
\label{eiej}
u\big[e^{\circ}_i(u),e_{j}(v)\big]-v\big[e_i(u),e^{\circ}_{j}(v)\big]
&=-(\al_i,\al_{j})\tss e_{i}(u)\tss e_{j}(v)\Fand\\[0.4em]
\label{fifj}
u\big[f^{\circ}_i(u),f_{j}(v)\big]-v\big[f_i(u),f^{\circ}_{j}(v)\big]
&=(\al_i,\al_{j})\tss f_{j}(v)\tss f_{i}(u).
\end{align}
Finally, for $i\ne j$ we have
the Serre relations
\begin{align}
\label{CompleteESerre}
\sum_{p\in\Sym_m}\big[e_{i}(u_{p(1)}),
\big[e_{i}(u_{p(2)}),\dots,\big[e_{i}(u_{p(m)}),e_{j}(v)\big]\dots\big]\big]&=0\qquad\text{and}\\
\label{CompleteFSerre}
\sum_{p\in\Sym_m}\big[f_{i}(u_{p(1)}),
\big[f_{i}(u_{p(2)}),\dots,\big[f_{i}(u_{p(m)}),f_{j}(v)\big]\dots\big]\big]&=0,
\end{align}
where $m=1-a_{ij}$.
\eth

\bpf
Apart from the Serre relations where $i$ or $j$ takes the value $n$,
all the relations are satisfied in the algebra $\X(\g_{N})$ due to
Propositions~\ref{pro:rellowrankBCD}, \ref{pro:Arelations}, \ref{prop:relhn+1en-1},
\ref{pro:relenei} and \ref{pro:inrelation}. To derive \eqref{CompleteESerre} and
\eqref{CompleteFSerre} we use a theorem of Levendorski\v{\i}~\cite{l:gd}
which provides a simplified presentation of the Drinfeld Yangian $\Y^{D}(\g)$;
see also \cite{gu:co}. In particular, the theorem implies that
the Serre relations \eqref{xisym} in the definition of $\Y^{D}(\g)$ (see Section~\ref{sec:int})
can be replaced by their level zero case $r_1=\dots=r_m=s=0$.
As we demonstrate in Section~\ref{sec:isom} below,
the defining relations of $\Y^{D}(\g_N)$ are satisfied by the
elements of the extended Yangian $\X(\g_{N})$ which are defined by
the respective coefficients of the series $\ka_i(u)$, $\xi_i^{\pm}(u)$
in Section~\ref{sec:int}. This includes the level zero case
of \eqref{xisym} which is implied by the embedding $\U(\g_N)\hra \X(\g_{N})$;
see \cite[Proposition~3.11]{amr:rp}. Hence, by \cite{l:gd},
relations \eqref{xisym} hold for the coefficients of the series $\xi_i^{\pm}(u)$.
However, the series
$\xi_i^-(u)$ and $\xi_i^+(u)$ coincide with $e_i(u)$ and $f_i(u)$,
respectively, up to a shift of the variable $u$ by a constant.
This shift does not affect the Serre relations, and so
we can conclude that \eqref{CompleteESerre} and \eqref{CompleteFSerre} hold as well
for all $i\ne j$.

Now consider the algebra $\wh \X(\g_{N})$
with generators and relations as in the statement of the theorem.
The above argument shows that there is a homomorphism
\beql{surjhom}
\wh \X(\g_{N})\to\X(\g_{N})
\eeq
which takes
the generators $h_{i}^{(r)}$,
$e_{i}^{(r)}$ and $f_{i}^{(r)}$ of $\wh \X(\g_{N})$ to the elements
of $\X(\g_{N})$ denoted by the same symbols. We need to demonstrate that
this homomorphism is surjective and injective.
To prove the surjectivity we need a lemma.

\ble\label{lem:eijr}
For all $1\leqslant i<j\leqslant n$ in the algebra
$\X(\oa_{2n+1})$ we have
\begin{alignat}{2}
e_{i\ts j+1}^{(r)}&=\big[e_{i\tss j}^{(r)},\ts e_j^{(1)}\big],\qquad
&f_{j+1\ts i}^{(r)}&=\big[f_j^{(1)},\ts f_{j\tss i}^{(r)}\big],
\non\\[0.3em]
e_{i\ts j'}^{(r)}&=-\big[e_{i\tss j'-1}^{(r)},\ts e_{j}^{(1)}\big],\qquad
&f_{j'\ts i}^{(r)}&=-\big[f_{j}^{(1)},\ts f_{j'-1\tss i}^{(r)}\big].
\non
\end{alignat}
For all $1\leqslant i<j\leqslant n-1$ in the algebra
$\X(\spa_{2n})$ we have
\begin{alignat}{2}
e_{i\ts j+1}^{(r)}&=\big[e_{i\tss j}^{(r)},\ts e_j^{(1)}\big],\qquad
&f_{j+1\ts i}^{(r)}&=\big[f_j^{(1)},\ts f_{j\tss i}^{(r)}\big],
\non\\[0.3em]
e_{i\ts j'}^{(r)}&=-\big[e_{i\tss j'-1}^{(r)},\ts e_{j}^{(1)}\big],\qquad
&f_{j'\ts i}^{(r)}&=-\big[f_{j}^{(1)},\ts f_{j'-1\tss i}^{(r)}\big].
\non
\end{alignat}
Moreover, for $1\leqslant i\leqslant n-1$ we have
\begin{alignat}{2}
e_{i\ts n'}^{(r)}&=\frac{1}{2}\ts\big[e_{i\tss n'-1}^{(r)},\ts e_{n}^{(1)}\big],\qquad
&f_{n'\ts i}^{(r)}&=\frac{1}{2}\ts\big[f_{n}^{(1)},\ts f_{n'-1\tss i}^{(r)}\big],
\non\\[0.3em]
e_{i\ts i'}^{(r)}&=-\big[e_{i\tss i'-1}^{(r)},\ts e_{i}^{(1)}\big]-\sum_{k=1}^{r-1}e_{i}^{(k)}e_{i\tss i'-1}^{(r-k)},\qquad
&f_{i'\ts i}^{(r)}&=-\big[f_{i}^{(1)},\ts f_{i'-1\tss i}^{(r)}\big]-\sum_{k=1}^{r-1}f_{i'-1 \tss i}^{(r-k)}f_{i}^{(k)}.
\non
\end{alignat}
For all $1\leqslant i<j\leqslant n-1$ in the algebra
$\X(\oa_{2n})$ we have
\begin{alignat}{2}
e_{i\ts j+1}^{(r)}&=\big[e_{i\tss j}^{(r)},\ts e_j^{(1)}\big],\qquad
&f_{j+1\ts i}^{(r)}&=\big[f_j^{(1)},\ts f_{j\tss i}^{(r)}\big],
\non\\[0.3em]
e_{i\ts j'}^{(r)}&=-\big[e_{i\tss j'-1}^{(r)},\ts e_{j}^{(1)}\big],\qquad
&f_{j'\ts i}^{(r)}&=-\big[f_{j}^{(1)},\ts f_{j'-1\tss i}^{(r)}\big],
\non
\end{alignat}
and for $1\leqslant i\leqslant n-2$ we have
\ben
e_{i\ts n'}^{(r)}=\big[e_{i\tss n-1}^{(r)},\ts e_{n}^{(1)}\big],\qquad
f_{n'\ts i}^{(r)}=\big[f_{n}^{(1)},\ts f_{n-1\tss i}^{(r)}\big].
\een
\ele

\bpf{}\hspace{-0.5em}\footnote{We thank Oleksandr Tsymbaliuk
for pointing out mistakes in Lemma~\ref{lem:eijr} and some calculations
of this section in the previous version of the paper.}
All relations follow easily from the Gauss decomposition \eqref{gd}
and defining relations \eqref{defrel}. To illustrate, consider the
case $\g_N=\spa_{2n}$. By taking the coefficients of $v^{-1}$
on both sides of \eqref{defrel}, for $1<j\leqslant n-1$ we get
$[t_{1j}(u),t_{j\tss j+1}^{(1)}]=t_{1\tss j+1}(u)$.
Writing this in terms of the Gaussian generators we come to
the relation $h_{1}(u)[e_{1j}(u),e_{j}^{(1)}]=h_{1}(u)e_{1\tss j+1}(u)$,
which gives
$[e_{1j}(u),e_{j}^{(1)}]=e_{1\tss j+1}(u)$ so that
$e_{1\ts j+1}^{(r)}=[e_{1\tss j}^{(r)},\ts e_j^{(1)}]$.

By a similar argument, for $1\leqslant j\leqslant n-1$ we find that
$e_{1\ts j'}^{(r)}=[e_{1\tss j'-1}^{(r)},\ts e_{(j+1)'\tss j'}^{(1)}]$.
Now apply Proposition~\ref{prop:eiei'} to write this as
$e_{1\ts j'}^{(r)}=-[e_{1\tss j'-1}^{(r)},\ts e_{j}^{(1)}]$.
Similarly,
\ben
\bal{}
[t_{1,2'}(u),t_{2'\tss 1'}^{(1)}]&=[h_1(u)e_{1,2'}(u),e_{2',1'}^{(1)}]
=-[h_1(u)e_{1,2'}(u),e_{1}^{(1)}]\\
&=h_{1}(u)[e_{1}^{(1)},e_{1,2'}(u)]-[h_{1}(u),e_{1}^{(1)}]e_{1,2'}(u)
=h_{1}(u)e_{1\tss 1'}(u).
\eal
\een
Since
$[h_1(u), e_1^{(1)}]=h_1(u)e_{1}(u)$,
we get
\ben
[e_{1}^{(1)},e_{1,2'}(u)]-e_{1}(u)e_{1,2'}(u)
=e_{1\tss 1'}(u),
\een
which is equivalent to $e_{1\ts 1'}^{(r)}=-\big[e_{1\tss 1'-1}^{(r)},\ts e_{1}^{(1)}\big]-\sum_{k=1}^{r-1}e_{1}^{(k)}e_{1\tss 1'-1}^{(r-k)}$.

The remaining cases with $i=1$ are treated in the same way.
The extension to arbitrary values of $i$ follows by
the application of Corollary~\ref{cor:guass-embed}.
\epf

By Lemma~\ref{lem:eijr}, all elements
$e_{ij}^{(r)}$ and $f_{ji}^{(r)}$ with $r\geqslant 1$
and the conditions $i<j$ and $i<j\pr$ in the orthogonal case, and
$i<j$ and $i\leqslant j\pr$ in the symplectic case, belong to the subalgebra $\wt\X(\gl_N)$
of $\X(\g_{N})$ generated by the coefficients of the series
$h_i(u)$ with
$i=1,\dots,n+1$, and
$e_i(u)$, $f_i(u)$ with $i=1,\dots,n$. Hence,
the Gauss decomposition \eqref{gd} implies that all coefficients
of the series $t_{ij}(u)$ with the same respective conditions
on the indices $i$ and $j$ also belong to the subalgebra $\wt\X(\gl_N)$.
Furthermore, by Theorem~\ref{thm:Center}, all coefficients of the series $z_N(u)$
are also in $\wt\X(\gl_N)$. Finally, taking the coefficients of $u^{-r}$
for $r=1,2,\dots$ in \eqref{zcenter} and using induction on $r$, we conclude that
the coefficients of all series $t_{ij}(u)$ belong to $\wt\X(\gl_N)$
so that $\wt\X(\gl_N)=\X(\gl_N)$.
This proves that the homomorphism
\eqref{surjhom} is surjective.

In the rest of the proof we will
show that this homomorphism is injective.
We will follow the argument of \cite{bk:pp} dealing with type $A$,
and adapt it to the orthogonal and symplectic Lie algebras.
As a first step, observe that the set of monomials in
the generators
$h_{i}^{(r)}$
with $i=1,\dots,n+1$ and $r\geqslant 1$,
and $e_{ij}^{(r)}$ and $f_{ji}^{(r)}$ with $r\geqslant 1$
and the conditions $i<j$ and $i<j\pr$ in the orthogonal case, and
$i<j$ and $i\leqslant j\pr$ in the symplectic case,
taken in some fixed order, is linearly independent in the extended Yangian
$\X(\g_{N})$.
Indeed, under the isomorphism \eqref{isomgrpol},
the images of the elements $e_{ij}^{(r)}$
and $f_{ji}^{(r)}$ in the $(r-1)$-th component of the graded
algebra $\gr \X(\g_{N})$ respectively correspond to
$F_{ij}\tss x^{r-1}$ and $F_{ji}\tss x^{r-1}$. Similarly,
the image of $h_{i}^{(r)}$ correspond to
$F_{ii}\tss x^{r-1}+\ze_r/2$ for $i=1,\dots, n$, while
for the image of $h^{(r)}_{n+1}$ we have
\ben
\bar h^{(r)}_{n+1}\mapsto
\begin{cases}
\ze_r/2
\qquad&\text{for}\quad \oa_{2n+1}\\[0.2em]
-F_{n\tss n}\tss x^{r-1}+\ze_r/2
\qquad&\text{for}\quad \spa_{2n}\\[0.2em]
-F_{n-1\tss n-1}\tss x^{r-1}-F_{n\tss n}\tss x^{r-1}+\ze_r/2
\qquad&\text{for}\quad \oa_{2n},
\end{cases}
\een
which follows from \eqref{znim} and Theorem~\ref{thm:Center}.
Hence the claim is implied by
the Poincar\'e--Birkhoff--Witt theorem for $\U(\g_{N}[x])$.

Define elements $e_{ij}^{(r)}$ and $f_{ji}^{(r)}$ of $\wh \X(\g_{N})$
inductively as follows. For $\g_N=\oa_{2n+1}$ set $e_{i\ts i+1}^{(r)}=e_{i}^{(r)}$ and
$f_{i+1\ts i}^{(r)}=f_{i}^{(r)}$, and
\begin{alignat}{2}\label{Beijpone}
e_{i\ts j+1}^{(r)}&=[e_{i\tss j}^{(r)},\ts e_j^{(1)}],\qquad
&f_{j+1\ts i}^{(r)}&=[f_j^{(1)},\ts f_{j\tss i}^{(r)}],\\
e_{i\ts j'}^{(r)}&=-[e_{i\tss j'-1}^{(r)},\ts e_{j}^{(1)}],\qquad
&f_{j'\ts i}^{(r)}&=-[f_{j}^{(1)},\ts f_{j'-1\tss i}^{(r)}],
\non
\end{alignat}
for $1\leqslant i<j\leqslant n$. For $\g_N=\spa_{2n}$
set $e_{i\ts i+1}^{(r)}=e_{i}^{(r)}$ and
$f_{i+1\ts i}^{(r)}=f_{i}^{(r)}$, and
\begin{alignat}{2}\label{Ceijpone}
e_{i\ts j+1}^{(r)}&=[e_{i\tss j}^{(r)},\ts e_j^{(1)}],\qquad
&f_{j+1\ts i}^{(r)}&=[f_j^{(1)},\ts f_{j\tss i}^{(r)}],\\
e_{i\ts j'}^{(r)}&=-[e_{i\tss j'-1}^{(r)},\ts e_{j}^{(1)}],\qquad
&f_{j'\ts i}^{(r)}&=-[f_{j}^{(1)},\ts f_{j'-1\tss i}^{(r)}],
\non
\end{alignat}
for $1\leqslant i< j\leqslant n-1$. Furthermore, set
$e_{n\ts n'}^{(r)}=e_n^{(r)}$ and $f_{n'\ts n}^{(r)}=f_n^{(r)}$, and
\begin{alignat}{2}\label{Ceijponen}
&e_{i\ts n'}^{(r)}=\frac{1}{2}\ts[e_{i\tss n'-1}^{(r)},\ts e_{n}^{(1)}],\qquad
f_{n'\ts i}^{(r)}=\frac{1}{2}\ts[f_{n}^{(1)},\ts f_{n'-1\tss i}^{(r)}],\\
e_{i\ts i'}^{(r)}&=-\big[e_{i\tss i'-1}^{(r)},\ts e_{i}^{(1)}\big]-\sum_{k=1}^{r-1}e_{i}^{(k)}e_{i\tss i'-1}^{(r-k)},\quad
f_{i'\ts i}^{(r)}=-\big[f_{i}^{(1)},\ts f_{i'-1\tss i}^{(r)}\big]-\sum_{k=1}^{r-1}f_{i'-1 \tss i}^{(r-k)}f_{i}^{(k)}.
\non
\end{alignat}
for $1\leqslant i\leqslant n-1$. For $\g_N=\oa_{2n}$
set $e_{i\ts i+1}^{(r)}=e_{i}^{(r)}$ and
$f_{i+1\ts i}^{(r)}=f_{i}^{(r)}$, and
\begin{alignat}{2}\label{Deijpone}
e_{i\ts j+1}^{(r)}&=[e_{i\tss j}^{(r)},\ts e_j^{(1)}],\qquad
&f_{j+1\ts i}^{(r)}&=[f_j^{(1)},\ts f_{j\tss i}^{(r)}],\\
e_{i\ts j'}^{(r)}&=-[e_{i\tss j'-1}^{(r)},\ts e_{j}^{(1)}],\qquad
&f_{j'\ts i}^{(r)}&=-[f_{j}^{(1)},\ts f_{j'-1\tss i}^{(r)}],
\non
\end{alignat}
for $1\leqslant i< j\leqslant n-1$. Furthermore, set
$e_{n-1\ts n'}^{(r)}=e_n^{(r)}$ and $f_{n'\ts n-1}^{(r)}=f_n^{(r)}$, and
\beql{Deijponen}
e_{i\ts n'}^{(r)}=[e_{i\tss n-1}^{(r)},\ts e_{n}^{(1)}],\qquad
f_{n'\ts i}^{(r)}=[f_{n}^{(1)},\ts f_{n-1\tss i}^{(r)}]
\eeq
for $1\leqslant i\leqslant n-2$.

By Lemma~\ref{lem:eijr}, these definitions
are consistent with those of the elements of the algebra $\X(\g_{N})$
in the sense that the images of the elements
$e_{ij}^{(r)}$ and $f_{ji}^{(r)}$ of the algebra
$\wh \X(\g_{N})$ under the homomorphism \eqref{surjhom}
coincide with the elements of $\X(\g_{N})$ with the same name.

The injectivity property of the homomorphism \eqref{surjhom} will follow if we prove that
the algebra $\wh \X(\g_{N})$ is spanned by
monomials in
$h_{i}^{(r)}$, $e_{ij}^{(r)}$ and $f_{ji}^{(r)}$
taken in some fixed order.
Denote by $\wh \Ec$, $\wh \Fc$
and $\wh \Hc$ the subalgebras of $\wh \X(\g_{N})$ respectively
generated by all elements
of the form $e_{i}^{(r)}$, $f_{i}^{(r)}$ and $h_{i}^{(r)}$.
Define an ascending filtration
on $\wh \Ec$ by setting $\deg e_{i}^{(r)}=r-1$.
Denote by $\gr\wh \Ec$ the corresponding graded algebra.
Let $\eb_{ij}^{\tss(r)}$ be the image of $e_{ij}^{(r)}$ in the
$(r-1)$-th component of the graded algebra $\gr\wh \Ec$.
Extend the range of subscripts of
$\eb_{ij}^{\tss(r)}$ to all values $1\leqslant i<j\leqslant 1'$
by using the skew-symmetry conditions
\beql{skew-symmetry}
\eb_{i\tss j}^{\tss(r)}=-\theta_{ij}\eb_{j'\tss i'}^{\tss(r)}.
\eeq
The desired spanning property of
the monomials in the $e_{ij}^{(r)}$
clearly follows from the relations
\beql{TypeBeijeklbar}
[\eb_{i\tss j}^{\tss(r)},\eb_{k\tss l}^{\tss(s)}]=
\de_{k\tss j}\ts\eb_{i\tss l}^{\tss(r+s-1)}-\de_{i\tss l}\ts\eb_{kj}^{\tss(r+s-1)}
-\theta_{ij}\tss \de_{k\tss i'}\ts\eb_{j'\tss l}^{\tss(r+s-1)}
+\theta_{ij}\tss \de_{j'\tss l}\ts\eb_{k\tss i'}^{\tss(r+s-1)}.
\eeq
We will be verifying these
relations separately for each of the three cases.

\paragraph{Type $B_n$.}
If $i,j,k,l\in\{1,\dots,n\}$, then
\eqref{TypeBeijeklbar} are essentially type $A$ relations and they
were already verified in \cite{bk:pp}.
We will often use these particular cases
of \eqref{TypeBeijeklbar} in the arguments below.

By relation \eqref{eiej} (in the algebra $\wh \X(\oa_{2n+1})$)
we have $[\eb^{\tss(r)}_{i\ts i+1},\eb^{\tss(s)}_{n\ts n+1}]
=\delta_{i+1,n}\tss \eb^{\tss(r+s-1)}_{i\ts n+1}$.
For
$i<j<n$, using the definition \eqref{Beijpone} we obtain
\begin{multline}
\non
[\eb^{\tss(r)}_{i\ts j},\eb^{\tss(s)}_{n\ts n+1}]=
\big[[\eb^{\tss(r)}_{i\ts j-1},\eb^{\tss(1)}_{j-1\ts j}],\eb^{\tss(s)}_{n\ts n+1}\big]\\[0.3em]
{}=\big[\eb^{\tss(r)}_{i\ts j-1},[\eb^{\tss(1)}_{j-1\ts j},\eb^{(s)}_{n\tss n+1}]\big]
+\big[\eb^{\tss(1)}_{j-1\ts j},[\eb^{(s)}_{n\tss n+1},\eb^{\tss(r)}_{i\ts j-1}]\big]
=\big[\eb^{\tss(1)}_{j-1\ts j},[\eb^{(s)}_{n\tss n+1},\eb^{\tss(r)}_{i\ts j-1}]\big].
\end{multline}
Hence, an obvious induction gives $[\eb^{\tss(r)}_{i\ts j},\eb^{\tss(s)}_{n\ts n+1}]=0$
for $i<j<n$.
Now we verify $[\eb^{\tss(r)}_{i\ts n},\eb^{\tss(s)}_{n\ts n+1}]=\eb^{\tss(r+s-1)}_{i\ts n+1}.$
Using \eqref{Beijpone}, we get
\begin{multline}
\non
[\eb^{\tss(r)}_{i\ts n},\eb^{\tss(s)}_{n\ts n+1}]=
\big[[\eb^{\tss(r)}_{i\ts n-1},\eb^{\tss(1)}_{n-1\ts n}],\eb^{\tss(s)}_{n\ts n+1}\big]\\[0.3em]
=\big[\eb^{\tss(r)}_{i\ts n-1},[\eb^{\tss(1)}_{n-1\ts n},\eb^{(s)}_{n\tss n+1}]\big]
+\big[\eb^{\tss(1)}_{n-1\ts n},[\eb^{(s)}_{n\tss n+1},\eb^{\tss(r)}_{i\ts n-1}]\big]
=[\eb^{\tss(r)}_{i\ts n-1},\eb^{\tss(s)}_{n-1\ts n+1}],
\end{multline}
where the last equality holds by $[\eb^{\tss(r)}_{i\ts n-1},\eb^{\tss(s)}_{n\ts n+1}]=0$
and $[\eb^{\tss(s)}_{n-1\ts n},\eb^{\tss(1)}_{n\ts n+1}]=\bar{e}^{(s)}_{n-1\ts n+1}$.
Furthermore, $\eb^{\tss(s)}_{n-1\ts n+1}
=[\eb^{\tss(s)}_{n-1\ts n},\eb^{\tss(1)}_{n\ts n+1}]$ by \eqref{Beijpone},
and so
\begin{multline}
\non
[\eb^{\tss(r)}_{i\ts n},\eb^{\tss(s)}_{n\ts n+1}]=
\big[\eb^{\tss(r)}_{i\ts n-1},[\eb^{\tss(s)}_{n-1\ts n},\eb^{\tss(1)}_{n\ts n+1}]\big]\\[0.3em]
=\big[\eb^{\tss(s)}_{n-1\ts n},[\eb^{\tss(r)}_{i\ts n-1},\eb^{\tss(1)}_{n\ts n+1}]\big]+
\big[\eb^{\tss(1)}_{n\ts n+1},[\eb^{\tss(s)}_{n-1\ts n},\eb^{\tss(r)}_{i\ts n-1}]\big]
=-[\eb^{\tss(1)}_{n\ts n+1},\eb^{\tss(r+s-1)}_{i\ts n}]=\eb^{\tss(r+s-1)}_{i\ts n+1}.
\end{multline}
Thus, we have verified that $[\eb^{\tss(r)}_{i\ts j},\eb^{\tss(s)}_{n\ts n+1}]
=\delta_{jn}\tss\eb^{\tss(r+s-1)}_{i\ts n+1}$ for $1\leqslant i<j\leqslant n$.

Next we will check
\beql{eijekn+1}
[\eb^{\tss(r)}_{i\ts j},\eb^{\tss(s)}_{k\ts n+1}]=\delta_{jk}\tss
\eb^{\tss(r+s-1)}_{i\ts n+1}
\eeq
for $1\leqslant i<j\leqslant n$ and $1\leqslant k< n$.
Suppose first that $1\leqslant i<j<n$. We have
\begin{multline}
\non
[\eb^{\tss(r)}_{i\ts j},\eb^{\tss(s)}_{k\ts n+1}]=
\big[\eb^{\tss(r)}_{i\ts j},[\eb^{\tss(s)}_{k\ts n},\eb^{\tss(1)}_{n\ts n+1}]\big]\\[0.3em]
=\big[\eb^{\tss(s)}_{k\ts n},[\eb^{\tss(r)}_{i\ts j},\eb^{\tss(1)}_{n\ts n+1}]\big]+
\big[\eb^{\tss(1)}_{n\ts n+1},[\eb^{\tss(s)}_{k\ts n},\eb^{\tss(r)}_{i\ts j}]\big]
=-\delta_{kj}\tss [\eb^{\tss(1)}_{n\ts n+1},\eb^{\tss(r+s-1)}_{i\ts n}]
=\delta_{kj}\tss\eb^{\tss(r+s-1)}_{i\ts n+1}.
\end{multline}
Now let $1\leqslant i<j=n$.
Note first that
\begin{multline}
\non
[\eb^{\tss(r)}_{n-1\ts n},\eb^{\tss(s)}_{n-1\ts n+1}]=
\big[\eb^{\tss(r)}_{n-1\ts n},[\eb^{\tss(1)}_{n-1\ts n},\eb^{\tss(s)}_{n\ts n+1}]\big]
=\big[\eb^{\tss(1)}_{n-1\ts n},[\eb^{\tss(r)}_{n-1\ts n},\eb^{\tss(s)}_{n\ts n+1}]\big]\\[0.3em]
{}+\big[\eb^{\tss(s)}_{n\ts n+1},[\eb^{\tss(1)}_{n-1\ts n},\eb^{\tss(r)}_{n-1\ts n}]\big]
=\big[\eb^{\tss(1)}_{n-1\ts n},[\eb^{\tss(r)}_{n-1\ts n},\eb^{\tss(s)}_{n\ts n+1}]\big]
=\big[\eb^{\tss(1)}_{n-1\ts n},[\eb^{\tss(1)}_{n-1\ts n},\eb^{\tss(r+s)}_{n\ts n+1}]\big]=0,
\end{multline}
where the last equality holds by the Serre relations \eqref{CompleteESerre}
with $a_{n-1\tss n}=-1$. As a next step, verify
$[\eb^{\tss(r)}_{n-2\ts n},\eb^{\tss(s)}_{n-1\ts n+1}]=0$. Indeed, we have
\begin{multline}
\non
[\eb^{\tss(r)}_{n-2\ts n},\eb^{\tss(s)}_{n-1\ts n+1}]=
\big[[\eb^{\tss(r)}_{n-2\ts n-1},\eb^{\tss(1)}_{n-1\ts n}],
[\eb^{\tss(1)}_{n-1\ts n},\eb^{\tss(s)}_{n\ts n+1}]\big]\\[0.3em]
=\big[\eb^{\tss(1)}_{n-1\ts n},[[\eb^{\tss(r)}_{n-2\ts n-1},
\eb^{\tss(1)}_{n-1\ts n}],\eb^{\tss(s)}_{n\ts n+1}]\big]
+\big[\eb^{\tss(s)}_{n\ts n+1},[\eb^{\tss(1)}_{n-1\ts n},
[\eb^{\tss(r)}_{n-2\ts n-1},\eb^{\tss(1)}_{n-1\ts n}]]\big].
\end{multline}
Since the second term is zero, this equals
\ben
\big[\eb^{\tss(1)}_{n-1\ts n},[\eb^{\tss(1)}_{n-1\ts n},
[\eb^{\tss(s)}_{n\ts n+1},\eb^{\tss(r)}_{n-2\ts n-1}]]\big]
+\big[\eb^{\tss(1)}_{n-1\ts n},[\eb^{\tss(r)}_{n-2\ts n-1},
[\eb^{\tss(1)}_{n-1\ts n},\eb^{\tss(s)}_{n\ts n+1}]]\big],
\een
where the first term is zero, so the expression equals
\ben
\big[[\eb^{\tss(1)}_{n-1\ts n},\eb^{\tss(s)}_{n\ts n+1}],
[\eb^{\tss(r)}_{n-2\ts n-1},\eb^{\tss(1)}_{n-1\ts n}]\big]+
\big[\eb^{\tss(r)}_{n-2\ts n-1},[\eb^{\tss(1)}_{n-1\ts n},
[\eb^{\tss(1)}_{n-1\ts n},\eb^{\tss(s)}_{n\ts n+1}]]\big]
=[\eb^{\tss(s)}_{n-1\ts n+1},\eb^{\tss(r)}_{n-2\ts n}]=0.
\een
Furthermore, for $k<n-1$ write the commutator $[\eb^{\tss(r)}_{n-1\ts n},\eb^{\tss(s)}_{k\ts n+1}]$
as
\begin{multline}
\non
\big[\eb^{\tss(r)}_{n-1\ts n},[\eb^{\tss(1)}_{k\ts k+1},\eb^{\tss(s)}_{k+1\ts n+1}]\big]
=\big[\eb^{\tss(1)}_{k\ts k+1},[\eb^{\tss(r)}_{n-1\ts n},\eb^{\tss(s)}_{k+1\ts n+1}]\big]
+\big[\eb^{\tss(s)}_{k+1\ts n+1},[\eb^{\tss(1)}_{k\ts k+1},\eb^{\tss(r)}_{n-1\ts n}]\big]\\[0.3em]
=\big[\eb^{\tss(1)}_{k\ts k+1},[\eb^{\tss(r)}_{n-1\ts n},\eb^{\tss(s)}_{k+1\ts n+1}]\big]
+\delta_{k\ts n-2}\tss[\eb^{\tss(s)}_{n-1\ts n+1},\eb^{\tss(r)}_{n-2\ts n}]
=\big[\eb^{\tss(1)}_{k\ts k+1},[\eb^{\tss(r)}_{n-1\ts n},\eb^{\tss(s)}_{k+1\ts n+1}]\big].
\end{multline}
Hence, by induction,
the relation $[\eb^{\tss(r)}_{n-1\ts n},\eb^{\tss(s)}_{k\ts n+1}]=0$
holds for $k\leqslant n-1$.
Finally,
\ben
[\eb^{\tss(r)}_{i\ts n},\eb^{\tss(s)}_{k\ts n+1}]=
\big[[\eb^{\tss(r)}_{i\ts n-1},\eb^{\tss(1)}_{n-1\ts n}],\eb^{\tss(s)}_{k\ts n+1}\big]
=\big[\eb^{\tss(r)}_{i\ts n-1},[\eb^{\tss(1)}_{n-1\ts n},\eb^{\tss(s)}_{k\ts n+1}]\big]
+\big[\eb^{\tss(1)}_{n-1\ts n},[\eb^{\tss(s)}_{k\ts n+1},\eb^{\tss(r)}_{i\ts n-1}]\big]=0,
\een
which completes the verification of \eqref{eijekn+1}.

For the next case of \eqref{TypeBeijeklbar} to verify we take
the relations
\beql{bar:ein+1ekn+1}
[\eb^{\tss(r)}_{i\ts n+1},\eb^{\tss(s)}_{k\ts n+1}]
=\eb^{\tss(r+s-1)}_{k\ts i'}.
\eeq
We may assume that $i\geqslant k$.
If $i=k=n$ then the relations follow from \eqref{eiei}.
In its turn, this implies \eqref{bar:ein+1ekn+1} for $i=n$ and arbitrary $k$
by an argument similar to the one used in the proof of \eqref{eijekn+1}.
Furthermore,
a reverse induction on $i$, beginning with $i=n$ shows
that \eqref{bar:ein+1ekn+1} holds
for $i>k$. For the remaining case
of \eqref{bar:ein+1ekn+1} with $i=k<n$ use induction on $i$ starting with $i=n$.
For $i<n$ we have
\ben
\bal[]
[\eb^{\tss(r)}_{i\tss n+1},\eb^{\tss(s)}_{i\ts n+1}]&=
\big[[\eb^{\tss(1)}_{i\ts i+1},\eb^{\tss(r)}_{i+1\ts n+1}],
[\eb^{\tss(1)}_{i\ts i+1},\eb^{\tss(s)}_{i+1\ts n+1}]\big]\\[0.3em]
{}&=\big[[\eb^{\tss(1)}_{i\ts i+1},[\eb^{\tss(1)}_{i\ts i+1},\eb^{\tss(s)}_{i+1\ts n+1}]],
\eb^{\tss(r)}_{i+1\ts n+1}\big]
+\big[\eb^{\tss(1)}_{i\ts i+1},[\eb^{\tss(r)}_{i+1\ts n+1},
[\eb^{\tss(1)}_{i\ts i+1},\eb^{\tss(s)}_{i+1\ts n+1}]]\big].
\eal
\een
The first term vanishes, while by the induction hypothesis, the second term equals
\ben
\big[\eb^{\tss(1)}_{i\ts i+1},[[\eb^{\tss(r)}_{i+1\ts n+1},
\eb^{\tss(1)}_{i\ts i+1}],\eb^{\tss(s)}_{i+1\ts n+1}]\big]
=\big[[\eb^{\tss(r)}_{i+1\ts n+1},\eb^{\tss(1)}_{i\ts i+1}],
[\eb^{\tss(1)}_{i\ts i+1},\eb^{\tss(s)}_{i+1\ts n+1}]\big]
=-[\eb^{\tss(r)}_{i\tss n+1},\eb^{\tss(s)}_{i\ts n+1}]
\een
so that $[\eb^{\tss(r)}_{i\tss n+1},\eb^{\tss(s)}_{i\ts n+1}]=0$,
completing the proof of \eqref{bar:ein+1ekn+1}. As its consequence,
we derive a more general relation
\beql{Bijkl's}
[\eb^{\tss(r)}_{i\tss j},\eb^{\tss(s)}_{k\ts l'}]=
\de_{k\tss j}\ts\eb_{i\ts l'}^{\tss(r+s-1)}+\de_{j\tss l}\ts\eb_{k\ts i'}^{\tss(r+s-1)}
\eeq
which holds for all $1\leqslant i<j\leqslant n$ and $1\leqslant k<l\leqslant n$.
Indeed, using \eqref{bar:ein+1ekn+1} we get
\begin{multline}
\non
[\eb^{\tss(r)}_{i\tss j},\eb^{\tss(s)}_{k\ts l'}]=
\big[\eb^{\tss(r)}_{i\tss j},[\eb^{\tss(1)}_{l\tss n+1},
\eb^{\tss(s)}_{k\ts n+1}]\big]=\big[\eb^{\tss(1)}_{l\tss n+1},[\eb^{\tss(r)}_{i\tss j},
\eb^{\tss(s)}_{k\ts n+1}]\big]\\[0.3em]
{}+\big[\eb^{\tss(s)}_{k\ts n+1},[\eb^{\tss(1)}_{l\tss n+1},
\eb^{\tss(r)}_{i\tss j}]\big]
=[\eb^{\tss(1)}_{l\tss n+1},\de_{k\tss j}\ts\eb^{\tss(r+s-1)}_{i\ts n+1}]
+[\de_{j\tss l}\ts\eb^{\tss(r)}_{i\tss n+1},\eb^{\tss(s)}_{k\ts n+1}]
\end{multline}
which equals the right hand side of \eqref{Bijkl's}.
Our next goal is to verify the relations
\beql{Berinskl}
[\eb^{\tss(r)}_{i\tss n+1},\eb^{\tss(s)}_{k\ts l'}]=0
\eeq
for all admissible $i,k,l$ with $k<l\leqslant n$. We begin
with the particular case
\beql{enn+1kn'}
[\eb^{\tss(r)}_{n\tss n+1},\eb^{\tss(s)}_{k\ts n'}]=0
\eeq
which we check by a reverse induction on $k$.
Using the previously checked cases of \eqref{TypeBeijeklbar},
for $k=n-1$ we obtain
\ben
\bal[]
[\eb^{\tss(r)}_{n\tss n+1},\eb^{\tss(s)}_{n-1\ts n'}]
&=\big[\eb^{\tss(r)}_{n\ts n+1},[\eb^{\tss(1)}_{n\ts n+1},
\eb^{\tss(s)}_{n-1\ts n+1}]\big]\\[0.3em]
{}&=\big[\eb^{\tss(r)}_{n\ts n+1},[\eb^{\tss(1)}_{n\ts n+1},
[\eb^{\tss(s)}_{n-1\ts n},\eb^{\tss(1)}_{n\ts n+1}]]\big]
=-\big[\eb^{\tss(r)}_{n\ts n+1},[\eb^{\tss(1)}_{n\ts n+1},
[\eb^{\tss(1)}_{n\ts n+1},\eb^{\tss(s)}_{n-1\ts n}]]\big].
\eal
\een
Now apply the Serre relations \eqref{CompleteESerre} with $a_{n\ts n-1}=-2$
to write this commutator as
\ben
\big[\eb^{\tss(1)}_{n\ts n+1},[\eb^{\tss(r)}_{n\ts n+1},
[\eb^{\tss(1)}_{n\ts n+1},\eb^{\tss(s)}_{n-1\ts n}]]\big]
+\big[\eb^{\tss(1)}_{n\ts n+1},[\eb^{\tss(1)}_{n\ts n+1},
[\eb^{\tss(r)}_{n\ts n+1},\eb^{\tss(s)}_{n-1\ts n}]]\big].
\een
Hence,
\begin{multline}
\non
[\eb^{\tss(r)}_{n\tss n+1},\eb^{\tss(s)}_{n-1\ts n'}]
=2\big[\eb^{\tss(1)}_{n\ts n+1},[\eb^{\tss(r)}_{n\ts n+1},
[\eb^{\tss(1)}_{n\ts n+1},\eb^{\tss(s)}_{n-1\ts n}]]\big]
+\big[\eb^{\tss(1)}_{n\ts n+1},[\eb^{\tss(s)}_{n-1\ts n},
[\eb^{\tss(r)}_{n\ts n+1},\eb^{\tss(1)}_{n\ts n+1}]]\big]\\[0.3em]
=-2\big[\eb^{\tss(1)}_{n\ts n+1},[\eb^{\tss(r)}_{n\ts n+1},
\eb^{\tss(s)}_{n-1\ts n+1}]\big]=-2[\eb^{\tss(r)}_{n\ts n+1},\eb^{\tss(s)}_{n-1\ts n'}]=0
\end{multline}
establishing the induction base. The induction step is a straightforward
application of \eqref{eijekn+1}, which completes the proof of \eqref{enn+1kn'}.
Now assume that $i=n$ in \eqref{Berinskl} to show that
\beql{enn+1ekl'}
[\eb^{\tss(r)}_{n\ts n+1},\eb^{\tss(s)}_{k\ts l'}]=0
\eeq
for all $k<l\leqslant n$. By writing
\ben
[\eb^{\tss(r)}_{n\ts n+1},\eb^{\tss(s)}_{k\ts l'}]
=\big[\eb^{\tss(r)}_{n\ts n+1},[\eb^{\tss(1)}_{l\ts n+1},\eb^{\tss(s)}_{k\ts n+1}]\big]
\een
and relying on the already verified cases of \eqref{TypeBeijeklbar},
we reduce checking of \eqref{enn+1ekl'} to the particular case $r=1$
where we may also assume $l\leqslant n-1$. Proceed by
\ben
\big[\eb^{\tss(1)}_{n\ts n+1},\eb^{\tss(s)}_{k\ts l'}\big]
=\big[\eb^{\tss(1)}_{n\ts n+1},[\eb^{\tss(1)}_{l\ts n},\eb^{\tss(s)}_{k\ts n'}]\big]
=\big[\eb^{\tss(1)}_{l\ts n},[\eb^{\tss(1)}_{n\ts n+1},\eb^{\tss(s)}_{k\ts n'}]\big]+
\big[\eb^{\tss(s)}_{k\ts n'},[\eb^{\tss(1)}_{l\ts n},\eb^{\tss(1)}_{n\ts n+1}]\big]
\een
and use \eqref{enn+1kn'} to see that the expression is zero.
Now \eqref{Berinskl} follows from \eqref{Bijkl's} and \eqref{enn+1ekl'}:
\ben
[\eb^{\tss(r)}_{i\ts n+1},\eb^{\tss(s)}_{k\ts l'}]=
\big[[\eb^{\tss(r)}_{i\ts n},\eb^{\tss(1)}_{n\ts n+1}],\eb^{\tss(s)}_{k\ts l'}\big]
=\big[\eb^{\tss(1)}_{n\ts n+1},[\eb^{\tss(s)}_{k\ts l'},\eb^{\tss(r)}_{i\ts n}]\big]+
\big[\eb^{\tss(r)}_{i\ts n},[\eb^{\tss(1)}_{n\ts n+1},\eb^{\tss(s)}_{k\ts l'}]\big]=0.
\een
The proof of \eqref{TypeBeijeklbar} will be completed by checking that
$[\eb^{\tss(r)}_{i\ts j'},\eb^{\tss(s)}_{k\ts l'}]=0$
for all admissible values with $i<j$ and $k<l$. This relation follows from
\eqref{Berinskl} by
\ben
[\eb^{\tss(r)}_{i\ts j'},\eb^{\tss(s)}_{k\ts l'}]
=\big[[\eb^{\tss(r)}_{j\ts n+1},\eb^{\tss(1)}_{i\ts n+1}],\eb^{\tss(s)}_{k\ts l'}\big]=
\big[\eb^{\tss(r)}_{j\ts n+1},[\eb^{\tss(1)}_{i\ts n+1},\eb^{\tss(s)}_{k\ts l'}]\big]+
\big[\eb^{\tss(1)}_{i\ts n+1},[\eb^{\tss(s)}_{k\ts l'},\eb^{\tss(r)}_{j\ts n+1}]\big]=0.
\een

\paragraph{Type $C_n$.} We will now
verify \eqref{TypeBeijeklbar} for $\g_N=\spa_{2n}$, where the arguments
are quite similar to type $B_n$. We will outline the sequence of steps.
By \eqref{eiej=0}, we have $[\eb^{\tss(r)}_{i\ts i+1},\eb^{\tss(s)}_{n\ts n'}\big]=0$ for $i<n-1$.
This implies
\beq\label{typeC eijenn'}
[\eb^{\tss(r)}_{i\ts j},\eb^{\tss(s)}_{n\ts n'}\big]=0
\eeq
for $i<j<n$ by an easy induction.
Using \eqref{eiej} and the
definition \eqref{Ceijponen}, we derive
\beq\label{typeC en-1nenn'}
[\eb_{n-1\ts n}^{\tss(r)},\eb_{n\ts n'}^{\tss(s)}]
=\eb_{n-1\ts n'}^{\tss(r+s-1)}+\eb_{n\tss (n-1)'}^{\tss(r+s-1)}
\eeq
which then implies a more general relation
\beql{typeC einenn'}
\big[\eb^{\tss(r)}_{i\ts n},\eb^{\tss(s)}_{n\ts n'}\big]
=2\tss\eb^{\tss(r+s-1)}_{i\ts n'}=\eb^{\tss(r+s-1)}_{i\ts n'}+\eb^{\tss(r+s-1)}_{n\ts i'}
\eeq
for $i<n-1$. Using \eqref{Ceijponen} we extend it further to
\beql{TypeC eijekn'}
\big[\eb^{\tss(r)}_{i\ts j},\eb^{\tss(s)}_{k\ts n'}\big]=\delta_{kj}\eb^{\tss(r+s-1)}_{i\ts n'},
\eeq
for all $i<j<n$ and $k<n$, and then apply
\eqref{typeC en-1nenn'}
to check
\beql{TypeC en-1nekn'}
\big[\eb^{\tss(r)}_{n-1\ts n},\eb^{\tss(s)}_{k\ts n'}\big]=\eb^{\tss(r+s-1)}_{k\ts (n-1)'},
\eeq
for $k<n$. Next, we verify
\beql{TypeC einekn'}
\big[\eb^{\tss(r)}_{i\ts n},\eb^{\tss(s)}_{k\ts n'}\big]=\eb^{\tss(r+s-1)}_{k\ts i'},
\eeq
for $k\leqslant n-1$ and $i<n-1$.
We use the reverse induction on $i$, beginning with $i=n-1$ to show first
that this relation holds
for $i\geqslant k$.
The induction base is \eqref{TypeC en-1nekn'}, while for $i<n-1$ write
$\big[\eb^{\tss(r)}_{i\ts n},\eb^{\tss(s)}_{k\ts n'}\big]
=\big[[\eb^{\tss(1)}_{i\ts i+1},\eb^{\tss(r)}_{i+1\ts n}],\eb^{\tss(s)}_{k\ts n'}\big]$
then proceed by using the definition \eqref{Ceijpone} and \eqref{typeC eijenn'}.
The case \eqref{TypeC einekn'} with $i<k$ now follows
by an application of \eqref{Ceijponen}.
Furthermore, assuming now that $i<j\leqslant n$ and $k\leqslant l<n$
we get from \eqref{TypeC einekn'}:
\ben
[\eb^{\tss(r)}_{i\tss j},\eb^{\tss(s)}_{k\ts l'}]=
\big[\eb^{\tss(r)}_{i\tss j},[\eb^{\tss(1)}_{l\tss n},
\eb^{\tss(s)}_{k\ts n'}]\big]=[\de_{j\tss l}\ts\eb^{\tss(r)}_{i\tss n},\eb^{\tss(s)}_{k\ts n'}]
+[\eb^{\tss(1)}_{l\tss n},\de_{k\tss j}\ts\eb^{\tss(r+s-1)}_{i\ts n'}]
\een
thus proving
\beql{TypeC eijekl'}
[\eb^{\tss(r)}_{i\ts j},\eb^{\tss(s)}_{k\ts l'}]
=\delta_{kj}\tss\eb^{\tss(r+s-1)}_{i\ts l'}+\delta_{jl}\tss\eb^{\tss(r+s-1)}_{k\ts i'}.
\eeq

Our next goal is to verify
\beql{TypeC einekl'}
[\eb^{\tss(r)}_{i\tss n},\eb^{\tss(s)}_{k\ts l'}]=0
\eeq
for all admissible $i,k,l$ with $k\leqslant l\leqslant n-1$.
We begin with the particular case
\beql{TypeC en-1nen-1n-1'}
[\eb^{\tss(r)}_{n-1\tss n},\eb^{\tss(s)}_{n-1\ts (n-1)'}]=0.
\eeq
By \eqref{TypeC einekn'}, we have
\ben
\bal[]
[\eb^{\tss(r)}_{n-1\tss n},\eb^{\tss(s)}_{n-1\ts (n-1)'}]&=
\big[\eb^{\tss(r)}_{n-1\tss n},[\eb^{\tss(1)}_{n-1\tss n},
\eb^{\tss(s)}_{n-1\ts n'}]\big]=\frac{1}{2}\ts\big[\eb^{\tss(r)}_{n-1\tss n},[\eb^{\tss(1)}_{n-1\tss n},
[\eb^{\tss(1)}_{n-1\ts n},\eb^{\tss(s)}_{n\ts n'}]]\big]\\[0.3em]
&=-\frac{1}{2}\ts\big[\eb^{\tss(1)}_{n-1\tss n},[\eb^{\tss(r)}_{n-1\tss n},
[\eb^{\tss(1)}_{n-1\ts n},\eb^{\tss(s)}_{n\ts n'}]]\big]
-\frac{1}{2}\ts\big[\eb^{\tss(1)}_{n-1\tss n},[\eb^{\tss(1)}_{n-1\tss n},
[\eb^{\tss(r)}_{n-1\ts n},\eb^{\tss(s)}_{n\ts n'}]]\big],
\eal
\een
where we also used the Serre relation \eqref{CompleteESerre} with $a_{n-1\ts n}=-2$.
Hence
\ben
[\eb^{\tss(r)}_{n-1\tss n},\eb^{\tss(s)}_{n-1\ts (n-1)'}]
=-\big[\eb^{\tss(1)}_{n-1\tss n},[\eb^{\tss(r)}_{n-1\tss n},
[\eb^{\tss(1)}_{n-1\ts n},\eb^{\tss(s)}_{n\ts n'}]]\big]
=-2\tss [\eb^{\tss(r)}_{n-1\tss n},\eb^{\tss(s)}_{n-1\ts (n-1)'}]
\een
and so \eqref{TypeC en-1nen-1n-1'} follows. By a reverse induction on $k$
we then derive
\beql{TypeC en-1nekn-1'}
[\eb^{\tss(r)}_{n-1\tss n},\eb^{\tss(s)}_{k\ts (n-1)'}]=0
\eeq
for $k<n-1$. This relations extends to all values $k\leqslant l\leqslant n-1$,
\beql{TypeC en-1nekl'}
[\eb^{\tss(r)}_{n-1\tss n},\eb^{\tss(s)}_{k\ts l'}]=0.
\eeq
The verification of \eqref{TypeC en-1nekl'} reduces to the case $r=1$
with the use of \eqref{TypeC eijekl'}, while in this particular case
it follows from \eqref{TypeC einekn'}, \eqref{TypeC eijekl'} and \eqref{TypeC en-1nekn-1'}.
Furthermore, by \eqref{TypeC eijekl'} and \eqref{TypeC en-1nekl'},
\ben
[\eb^{\tss(r)}_{i\ts n},\eb^{\tss(s)}_{k\ts l'}]
=\big[[\eb^{\tss(r)}_{i\ts n-1},\eb^{\tss(1)}_{n-1\ts n}],\eb^{\tss(s)}_{k\ts l'}\big]
=\big[\eb^{\tss(r)}_{i\ts n-1},[\eb^{\tss(1)}_{n-1\ts n},\eb^{\tss(s)}_{k\ts l'}]\big]+
\big[\eb^{\tss(1)}_{n-1\ts n},[\eb^{\tss(s)}_{k\ts l'},\eb^{\tss(r)}_{i\ts n-1}]\big]=0
\een
thus completing the proof of \eqref{TypeC einekl'}. As a next case of \eqref{TypeBeijeklbar}, we prove
\beql{TypeC ein'ekn'}
[\eb^{\tss(r)}_{i\ts n'},\eb^{\tss(s)}_{k\ts n'}]=0.
\eeq
If $i=k=n$ then this follows from
\eqref{eiei}. If $i=n-1$ and $k=n$ write
\ben
\bal
&[\eb^{\tss(r)}_{n-1\ts n'},\eb^{\tss(s)}_{n\ts n'}]
=\frac{1}{2}\ts\big[[\eb^{\tss(1)}_{n-1\ts n},\eb^{\tss(r)}_{n\ts n'}],\eb^{\tss(s)}_{n\ts n'}\big]
=\frac{1}{2}\ts\big[\eb^{\tss(s)}_{n\ts n'},[\eb^{\tss(r)}_{n\ts n'},\eb^{\tss(1)}_{n-1\ts n}]\big].
\eal
\een
Applying the Serre relation \eqref{CompleteESerre} with $a_{n\ts n-1}=-1$ we get
\ben
[\eb^{\tss(r)}_{n-1\ts n'},\eb^{\tss(s)}_{n\ts n'}]
=-\frac{1}{2}\ts\big[\eb^{\tss(r)}_{n\ts n'},[\eb^{\tss(s)}_{n\ts n'},\eb^{\tss(1)}_{n-1\ts n}]\big]
=-\frac{1}{2}\ts\big[\eb^{\tss(s)}_{n\ts n'},[\eb^{\tss(r)}_{n\ts n'},\eb^{\tss(1)}_{n-1\ts n}]\big]
=\big[\eb^{\tss(s)}_{n\ts n'},\eb^{\tss(r)}_{n-1\ts n'}\big]
\een
and so \eqref{TypeC ein'ekn'} with $i=n-1$ and $k=n$ follows.
Together with \eqref{typeC eijenn'} this implies
\beql{TypeC ein'enn'}
[\eb^{\tss(r)}_{i\ts n'},\eb^{\tss(s)}_{n\ts n'}]=0
\eeq
for all $i$. For the remaining values of the indices, \eqref{TypeC ein'ekn'}
follows by writing
\ben
[\eb^{\tss(r)}_{i\ts n'},\eb^{\tss(s)}_{k\ts n'}]
=\frac{1}{2}\ts\big[[\eb^{\tss(r)}_{i\ts n},\eb^{\tss(1)}_{n\ts n'}],\eb^{\tss(s)}_{k\ts n'}]\big]
\een
and applying
\eqref{typeC einenn'}, \eqref{TypeC einekn'} and \eqref{TypeC ein'enn'}.
Furthermore, we will verify
\beql{TypeC ein'ekl'}
[\eb^{\tss(r)}_{i\ts n'},\eb^{\tss(s)}_{k\ts l'}]=0
\eeq
for $k\leqslant l\leqslant n-1$.
By the definition of $e^{\tss(r)}_{i\ts n'}$ in \eqref{Ceijponen}, we have
\ben
\bal
\big[\eb^{(r)}_{i\ts n'},\eb^{(s)}_{k\ts l'}\big]
=\frac{1}{2}\big[[\eb_{in}^{(r)},\eb_{nn'}^{(1)}],\eb_{kl'}^{(s)}\big]
&=\frac{1}{2}\big[\eb_{in}^{(r)},[\eb_{nn'}^{(1)},\eb_{kl'}^{(s)}]\big]
+\frac{1}{2}\big[\eb_{nn'}^{(1)},[\eb_{kl'}^{(s)},\eb_{in}^{(r)}]\big]\\
&=\frac{1}{2}\big[\eb_{in}^{(r)},[\eb_{nn'}^{(1)},\eb_{kl'}^{(s)}]\big]
\eal
\een
 where the last equality holds by  \eqref{TypeC einekl'}.
 Next we will show
 \ben
 [\eb_{nn'}^{(1)},\eb_{kl'}^{(s)}]=0.
 \een
 By equations \eqref{TypeC en-1nekn'} and \eqref{TypeC einekn'}, we have $\eb_{kl'}^{(s)}=[\eb_{ln}^{(1)},\eb_{kn'}^{(s)}]$.
 Thus,
 \ben
 \bal
 \big[\eb_{nn'}^{(1)},\eb_{kl'}^{(s)}\big]=\big[\eb_{nn'}^{(1)},[\eb_{ln}^{(1)},\eb_{kn'}^{(s)}]\big]
 &=\big[\eb_{ln'}^{(1)},[\eb_{nn'}^{(1)},\eb_{kn'}^{(s)}]\big]
 +\big[\eb_{kn'}^{(s)},[\eb_{ln}^{(1)},\eb_{nn'}^{(1)}]\big]\\
 &=\frac{1}{2}\big[\eb_{kn'}^{(s)},\eb_{ln'}^{(1)}\big]
 \eal
 \een
 where the last equality holds by \eqref{TypeC ein'ekn'}
 and \eqref{Ceijponen}.
 Using \eqref{TypeC ein'ekn'} again, we get
 \ben
 \big[\eb_{nn'}^{(1)},\eb_{kl'}^{(s)}\big]=0
 \een
so that \eqref{TypeC ein'ekl'} holds.

The last remaining case of \eqref{TypeBeijeklbar} is
$[\eb^{\tss(r)}_{i\ts j'},\eb^{\tss(s)}_{k\ts l'}]=0$ which
holds for all $i\leqslant j\leqslant n-1$ and $k\leqslant l\leqslant n-1$.
Indeed, we have
\ben
[\eb^{\tss(r)}_{i\ts j'},\eb^{\tss(s)}_{k\ts l'}]
=\big[[\eb^{\tss(1)}_{j\ts n},\eb^{\tss(r)}_{i\ts n'}],\eb^{\tss(s)}_{k\ts l'}\big]
=\big[\eb^{\tss(1)}_{j\ts n},[\eb^{\tss(r)}_{i\ts n'},\eb^{\tss(s)}_{k\ts l'}]\big]+
\big[\eb^{\tss(r)}_{i\ts n'},[\eb^{\tss(s)}_{k\ts l'},\eb^{\tss(1)}_{j\ts n}]\big]
\een
which is zero by
\eqref{TypeC einekl'} and \eqref{TypeC ein'ekl'}.

\paragraph{Type $D_n$.} Now let $\g_N=\oa_{2n}$.
The arguments follow the same pattern as for types $B_n$ and $C_n$ above.
We will assume throughout the rest of the proof that the indices
$i,j,k,l$ run over the set $\{1,\dots,n\}$.
By \eqref{eiej=0} we have $[\eb^{\tss(r)}_{n-1\ts n},\eb^{\tss(s)}_{n-1\ts n'}]=0$.
Therefore, for
$i<n-1$, using the definition \eqref{Deijpone} we obtain
\ben
[\eb^{\tss(r)}_{i\ts n},\eb^{\tss(s)}_{n-1\ts n'}]=
\big[[\eb^{\tss(r)}_{i\ts n-1},\eb^{\tss(1)}_{n-1\ts n}],\eb^{\tss(s)}_{n-1\ts n'}\big]
=[\eb^{\tss(r+s-1)}_{i\ts n'},\eb^{\tss(1)}_{n-1\ts n}]=-\eb^{\tss(r+s-1)}_{i\ts (n-1)'}.
\een
Hence, for $k<n-1$ we derive
\beql{eknpr}
[\eb^{\tss(r)}_{n-1\ts n},\eb^{\tss(s)}_{k\ts n'}]=
\big[\eb^{\tss(r)}_{n-1\ts n},[\eb^{\tss(s)}_{k\ts n-1},\eb^{\tss(1)}_{n-1\ts n'}]\big]
=-[\eb^{\tss(r+s-1)}_{k\tss n},\eb^{\tss(1)}_{n-1\ts n'}]=\eb^{\tss(r+s-1)}_{k\ts (n-1)'}.
\eeq
Next we verify that
\beql{erins}
[\eb^{\tss(r)}_{i\tss n},\eb^{\tss(s)}_{k\ts n'}]=\eb^{\tss(r+s-1)}_{k\ts i'}
\eeq
by the reverse induction on $i$, beginning with $i=n-1$. Show first
that this relation holds
for $i>k$ by taking
\eqref{eknpr} as the induction base and using \eqref{Deijpone}.
To verify \eqref{erins} for $i=k$ write
\ben
[\eb^{\tss(r)}_{i\tss n},\eb^{\tss(s)}_{i\ts n'}]=
\big[[\eb^{\tss(1)}_{i\ts i+1},\eb^{\tss(r)}_{i+1\ts n}],
[\eb^{\tss(1)}_{i\ts i+1},\eb^{\tss(s)}_{i+1\ts n'}]\big]
\een
and proceed by induction.
As a next step, we point out the relations
\beql{iipo}
\eb_{i\ts j'}^{\tss(r)}=[\eb_{i\ts i+1}^{\tss(1)},\ts \eb_{i+1\ts j'}^{\tss(r)}]
\eeq
which hold for $j>i+1$.
Indeed, they follow by taking consecutive brackets of both sides
of the first relation in \eqref{eknpr}
with the elements
$\eb^{\tss(1)}_{n-1\ts n}$, $\eb^{\tss(1)}_{n-2\ts n-1},\dots,\eb^{\tss(1)}_{j\ts j+1}$
with the use of \eqref{eknpr}.
Now we check \eqref{erins} for $i\leqslant k$ by induction on $k-i$ taking
the case $i=k$ considered above as the induction base.
Suppose that $i<k$ and write
\beql{eirnrs}
[\eb^{\tss(r)}_{i\tss n},\eb^{\tss(s)}_{k\ts n'}]=
\big[[\eb^{\tss(1)}_{i\ts i+1},\eb^{\tss(r)}_{i+1\ts n}],\eb^{\tss(s)}_{k\ts n'}\big].
\eeq
If $i+1<k$, then by the induction hypothesis and \eqref{iipo} this equals
\ben
-[\eb^{\tss(1)}_{i\ts i+1},\eb^{\tss(r+s-1)}_{i+1\ts k'}]=-\eb^{\tss(r+s-1)}_{i\ts k'}.
\een
If $i+1=k$, then \eqref{eirnrs} equals
\ben
\big[[\eb^{\tss(1)}_{i\ts i+1},\eb^{\tss(s)}_{i+1\ts n'}],\eb^{\tss(r)}_{i+1\ts n}\big]
=[\eb^{\tss(s)}_{i\ts n'},\eb^{\tss(r)}_{i+1\ts n}]
=-[\eb^{\tss(r)}_{i+1\ts n},\eb^{\tss(s)}_{i\ts n'}]
=-\eb^{\tss(r+s-1)}_{i\ts (i+1)'}
\een
by \eqref{erins} with $i>k$ verified above. Thus, \eqref{erins} holds
for all admissible values of $i$ and $k$.

By employing \eqref{erins} we can derive a more general relation: for $j<n$ we have
\beql{erinjs}
[\eb^{\tss(r)}_{i\tss j},\eb^{\tss(s)}_{k\ts l'}]=
\de_{k\tss j}\ts\eb_{i\ts l'}^{\tss(r+s-1)}+\de_{j\tss l}\ts\eb_{k\ts i'}^{\tss(r+s-1)}.
\eeq
Indeed,
\ben
[\eb^{\tss(r)}_{i\tss j},\eb^{\tss(s)}_{k\ts l'}]=
\big[\eb^{\tss(r)}_{i\tss j},[\eb^{\tss(1)}_{l\tss n},
\eb^{\tss(s)}_{k\ts n'}]\big]=[\de_{j\tss l}\ts\eb^{\tss(r)}_{i\tss n},\eb^{\tss(s)}_{k\ts n'}]
+[\eb^{\tss(1)}_{l\tss n},\de_{k\tss j}\ts\eb^{\tss(r+s-1)}_{i\ts n'}]
\een
which gives \eqref{erinjs}.
Our next goal is to verify the relations
\beql{erinskl}
[\eb^{\tss(r)}_{i\tss n},\eb^{\tss(s)}_{k\ts l'}]=0
\eeq
for all admissible $i,k,l$ with $k<l\leqslant n-1$. We begin
with the particular case
\beql{ernns}
[\eb^{\tss(r)}_{n-1\tss n},\eb^{\tss(s)}_{k\ts (n-1)'}]=0
\eeq
which we check by reverse induction on $k$. For $k=n-2$ by
\eqref{erins} we have
\ben
[\eb^{\tss(r)}_{n-1\tss n},\eb^{\tss(s)}_{n-2\ts (n-1)'}]
=\big[\eb^{\tss(r)}_{n-1\ts n},[\eb^{\tss(1)}_{n-1\ts n},
[\eb^{\tss(1)}_{n-2\ts n-1},\eb^{\tss(s)}_{n-1\ts n'}]]\big]
=-\big[\eb^{\tss(1)}_{n-1\ts n},[\eb^{\tss(r)}_{n-2\ts n},\eb^{\tss(s)}_{n-1\ts n'}]\big]=0.
\een
For $k<n-2$ use \eqref{iipo} to write
\ben
[\eb^{\tss(r)}_{n-1\ts n},\eb^{\tss(s)}_{k\ts (n-1)'}]
=\big[\eb^{\tss(r)}_{n-1\ts n},
[\eb^{\tss(1)}_{k\ts k+1},\eb^{\tss(s)}_{k+1\ts (n-1)'}]\big]
\een
which is zero by the induction hypothesis. This verifies \eqref{ernns}.
We will now extend \eqref{ernns} to all values $k<l\leqslant n-1$:
\beql{ernnstwo}
[\eb^{\tss(r)}_{n-1\ts n},\eb^{\tss(s)}_{k\ts l'}]=0.
\eeq
First, verify it for $r=1$. Assume that $l\leqslant n-2$.
By using \eqref{erinjs} with $j=l=n-1$ and
\eqref{ernns} we get
\ben
\bal[]
[\eb^{\tss(1)}_{n-1\ts n},\eb^{\tss(s)}_{k\ts l'}]
&=\big[\eb^{\tss(1)}_{n-1\ts n},[\eb^{\tss(1)}_{l\ts n-1},\eb^{\tss(s)}_{k\ts (n-1)'}]\big]
=-[\eb^{\tss(1)}_{l\ts n},\eb^{\tss(s)}_{k\ts (n-1)'}]\\
{}&=-\big[\eb^{\tss(1)}_{l\ts n},[\eb^{\tss(1)}_{n-1\ts n},\eb^{\tss(s)}_{k\ts n'}]\big]=
-[\eb^{\tss(1)}_{n-1\ts n},\eb^{\tss(s)}_{k\ts l'}],
\eal
\een
where we also used \eqref{erins}. Hence \eqref{ernnstwo}
holds for $r=1$. Furthermore, by a similar transformation we find
\ben
\bal[]
[\eb^{\tss(r)}_{n-1\ts n},\eb^{\tss(s)}_{k\ts l'}]
&=\big[\eb^{\tss(r)}_{n-1\ts n},[\eb^{\tss(1)}_{l\ts n},\eb^{\tss(s)}_{k\ts n'}]\big]
=[\eb^{\tss(1)}_{l\ts n},\eb^{\tss(r+s-1)}_{k\ts (n-1)'}]\\
{}&=\big[\eb^{\tss(1)}_{l\ts n},[\eb^{\tss(1)}_{n-1\ts n},\eb^{\tss(r+s-1)}_{k\ts n'}]\big]
=[\eb^{\tss(1)}_{n-1\ts n},\eb^{\tss(r+s-1)}_{k\ts l'}]=0,
\eal
\een
thus proving \eqref{ernnstwo}.
Consider \eqref{erinjs} with $j=n-1$ and $k<l\leqslant n-1$ and
take the bracket of both sides with $\eb^{\tss(1)}_{n-1\ts n}$. This
finally gives \eqref{erinskl}.

To complete checking \eqref{TypeBeijeklbar}, it remains to show that
\beql{eijprklprxxx}
[\eb^{\tss(r)}_{i\ts j'},\eb^{\tss(s)}_{k\ts l'}]=0
\eeq
for all admissible values with $i<j$ and $k<l$. To this end, observe
that the mapping $\si:\wh \X(\oa_{2n})\to \wh \X(\oa_{2n})$
which swaps generators of each pair $(e^{\tss(r)}_{n-1},e^{\tss(r)}_{n})$,
$(f^{\tss(r)}_{n-1},f^{\tss(r)}_{n})$ and $(h^{\tss(r)}_{n},h^{\tss(r)}_{n'})$,
and leaves all other generators
unchanged,
defines an involutive automorphism of $\wh \X(\oa_{2n})$. The subalgebra
$\wh\Ec$ is invariant under $\si$ and we have the induced automorphism $\bar\si$
of the associated graded algebra $\gr\wh\Ec$. The definitions \eqref{Deijpone}
and \eqref{Deijponen}
imply that $\bar\si$ acts by
\ben
\eb^{\tss(r)}_{i\tss n}\mapsto \eb^{\tss(r)}_{i\tss n'},\qquad
\eb^{\tss(r)}_{i\tss n'}\mapsto \eb^{\tss(r)}_{i\tss n}
\een
for $i<n$, and leaves each element
$\eb^{\tss(r)}_{i\tss j}$ and $\eb^{\tss(r)}_{i\tss j'}$ unchanged for $i<j\leqslant n-1$.
Hence, applying $\bar\si$ to \eqref{erinskl} we get
\beql{erinprskl}
[\eb^{\tss(r)}_{i\tss n'},\eb^{\tss(s)}_{k\tss l'}]=0
\eeq
for all admissible $i,k,l$ with $k<l\leqslant n-1$. Finally, if $i<j\leqslant n-1$
and $k<l\leqslant n-1$, then using \eqref{erinskl} and \eqref{erinprskl} we get
\ben
[\eb^{\tss(r)}_{i\tss j'},\eb^{\tss(s)}_{k\tss l'}]
=\big[[\eb^{\tss(1)}_{j\tss n},\eb^{\tss(r)}_{i\tss n'}],\eb^{\tss(s)}_{k\tss l'}\big]=0
\een
thus verifying the remaining cases of \eqref{TypeBeijeklbar} for type $D_n$.

\medskip

To complete the proof of the theorem, note that in all three types
relations \eqref{TypeBeijeklbar} imply that the graded algebra $\gr\wh \Ec$ is spanned
by the set of monomials in the elements $\eb_{ij}^{\tss(r)}$
taken in some fixed order. Hence the algebra $\wh \Ec$ is spanned
by the corresponding monomials in the elements $e_{ij}^{\tss(r)}$.
It is immediate from
the defining relation of $\wh \X(\g_{N})$
that the mapping
\ben
h_i^{(r)}\mapsto h_i^{(r)},\qquad e_i^{(r)}\mapsto f_i^{(r)},\qquad
f_i^{(r)}\mapsto e_i^{(r)}
\een
defines an anti-automorphism of $\wh \X(\g_{N})$. By applying
this anti-automorphism, we deduce that the ordered monomials
in the elements $f_{ji}^{(r)}$ span the subalgebra $\wh\Fc$.
Note also that the ordered monomials in $h_i^{(r)}$
span $\wh\Hc$. Furthermore,
by the defining relations of $\wh \X(\g_{N})$, the multiplication map
\ben
\wh\Fc\ot\wh\Hc\ot\wh \Ec\to \wh \X(\g_{N})
\een
is surjective. Thus, ordering the elements
$h_{i}^{(r)}$, $e_{ij}^{(r)}$ and $f_{ji}^{(r)}$ in such a way that
the elements of $\wh\Fc$ precede the elements
of $\wh\Hc$, and the latter precede the elements of $\wh \Ec$,
we can conclude that the ordered monomials in these elements
span $\wh \X(\g_{N})$. This
proves that \eqref{surjhom} is an isomorphism.
\epf

We point out another version of
the Poincar\'e--Birkhoff--Witt theorem for $\X(\g_{N})$
which is implied by the proof of Theorem~\ref{thm:dp}.
Denote by $\Ec$, $\Fc$
and $\Hc$ the subalgebras of $\X(\g_{N})$ respectively
generated by all elements
of the form $e_{i}^{(r)}$, $f_{i}^{(r)}$ and $h_{i}^{(r)}$.
Consider the generators
$h_{i}^{(r)}$
with $i=1,\dots,n+1$ and $r\geqslant 1$,
and $e_{ij}^{(r)}$ and $f_{ji}^{(r)}$ with $r\geqslant 1$
and the conditions $i<j$ and $i<j\pr$ in the orthogonal case, and
$i<j$ and $i\leqslant j\pr$ in the symplectic case.
Order the elements
$h_{i}^{(r)}$, $e_{ij}^{(r)}$ and $f_{ji}^{(r)}$ in such a way that
the elements of $\Fc$ precede the elements
of $\Hc$, and the latter precede the elements of $\Ec$.

\bco\label{cor:pbwdp}
The set of all ordered monomials in the elements
$h_{i}^{(r)}$, $e_{ij}^{(r)}$ and $f_{ji}^{(r)}$
with the respective conditions on the indices
forms a basis of $\X(\g_{N})$.
\qed
\eco

\section{Isomorphism theorem for the Drinfeld Yangian}
\label{sec:isom}

We will now prove the Main Theorem as stated in the Introduction.
By the $R$-matrix definition of the Yangian $\Y(\g_N)=\Y^{R}(\g_N)$ in Section~\ref{sec:nd},
this is the subalgebra of $\X(\g_N)$, whose elements are stable under all automorphisms
\eqref{muf}. It is clear from the definition
of the series $\kappa^{}_i(u)$ and $\xi_i^{\pm}(u)$
with $i=1,\dots,n$, and the explicit formulas
\eqref{hmqua}, \eqref{eijmlqua} and \eqref{fijlmqua}, that all the
coefficients $\kappa^{}_{i\tss r}$ and
$\xi_{i\tss r}^{\pm}$ defined in
\eqref{kaxicoeff} belong to the subalgebra $\Y(\g_N)$.

\bpr\label{prop:YgN}
The subalgebra $\Y(\g_N)$ of $\X(\g_N)$ is generated by
the elements $\kappa^{}_{i\tss r}$, $\xi_{i\tss r}^{+}$ and
$\xi_{i\tss r}^{-}$ with
$i=1,\dots,n$ and $r=0,1,\dots$.
\epr

\bpf
Due to the tensor decomposition \eqref{tensordecom} of $\X(\g_N)$, it suffices to check that
these elements
together with the coefficients $z^{(r)}_N$ of the series $z_N(u)$
given in \eqref{zn} generate the
algebra $\X(\g_N)$.
Using the definition of the series $\ka^{}_i(u)$ and formulas \eqref{defkn}
it is straightforward to express $h_{1}(u)h_{n+1}(u)^{-1}$ as a product of the
series of the form $\ka^{}_i(u)^{-1}$ with some shifts of $u$ by constants.
On the other hand, Theorem~\ref{thm:Center} implies that
$h_{1}(u+\kappa)h_{n+1}(u)$ equals $z_N(u)$
times the same kind of product of the shifted series $\ka^{}_i(u)^{-1}$.
Therefore, all coefficients of $h_1(u)$ and hence all coefficients
of the series $h_i(u)$ with $i=1,\dots,n+1$ belong to the subalgebra
of $\X(\g_N)$ generated by the elements given in the proposition.
Furthermore, for each $i$, the elements $e_i^{(r)}$ and $f_i^{(r)}$ are found as
linear combinations of the $\xi_{i\tss s}^{-}$ and
$\xi_{i\tss s}^{+}$, respectively.
By Theorem~\ref{thm:dp}, the coefficients of the series
$h_i(u)$ with $i=1,\dots,n+1$, and $e_i(u)$ and $f_i(u)$ with
$i=1,\dots,n$ generate the algebra
$\X(\g_N)$ thus completing the proof.
\epf

Now we will verify that the generators
$\kappa^{}_{i\tss r}$, $\xi_{i\tss r}^{+}$ and
$\xi_{i\tss r}^{-}$ of the subalgebra $\Y(\g_N)$
provided by Proposition~\ref{prop:YgN} satisfy the
defining relations of the Drinfeld Yangian $\Y^{D}(\g_N)$
as given in the Introduction.
We will do this
in terms of the equivalent generating series relations for
$\kappa^{}_i(u)$ and $\xi_i^{\pm}(u)$. We use the notation $\{a,b\}=ab+ba$.

\bpr\label{pro:relinYgN}
The following relations hold in $\Y(\g_N)$:
\begin{align}
\label{kikj}
\big[\kappa_i(u),\kappa_j(v)\big]&=0,\\
\label{xpixmj}
\big[\xi_{i}^{+}(u),\xi_{j}^{-}(v)\big]&=-\de_{i\tss j}\ts\frac{\kappa_i(u)-\kappa_i(v)}{u-v},\\
\label{kixpj}
\big[\kappa_i(u),\xi^{\pm}_j(v)\big]&={}\mp\frac{1}{2}\ts\big(\al_i,\al_j\big)\ts
\frac{\big\{\kappa_i(u),\xi^{\pm}_j(u)-\xi^{\pm}_j(v)\big\}}{u-v},\\
\label{xpixpj}
\big[\xi^{\pm}_i(u),\xi^{\pm}_{j}(v)\big]
+\tss\big[\xi^{\pm}_j(u),\xi^{\pm}_{i}(v)\big]&=
{}\mp\frac{1}{2}\ts(\al_i,\al_j)
\frac{\big\{\xi^{\pm}_i(u)-\xi^{\pm}_i(v),
\xi^{\pm}_j(u)-\xi^{\pm}_j(v)\big\}}{u-v},
\end{align}
\beql{Serrexipm}
\sum_{p\in\Sym_m}\big[\xi_i^{\pm}(u_{p(1)}),
\big[\xi_i^{\pm}(u_{p(2)}),\dots,
\big[\xi_i^{\pm}(u_{p(m)}),\xi_j^{\pm}(v)\big]\dots\big]\big]=0,
\eeq
where the last relation holds for
all $i\ne j$, and we set $m=1-a_{ij}$.
\epr

\bpf
The proof amounts to writing the relations of Theorem~\ref{thm:dp} in terms
the series $\kappa^{}_i(u)$ and $\xi_i^{\pm}(u)$. In particular,
\eqref{kikj} and \eqref{xpixmj} are immediate from \eqref{hihj}
and \eqref{eifj}, respectively. Moreover, the Serre relations
\eqref{Serrexipm} are implied by \eqref{CompleteESerre}
and \eqref{CompleteFSerre}. In the case where $i,j\leqslant n-1$,
relations \eqref{kixpj} and \eqref{xpixpj} hold due to the corresponding
results in type $A$ as shown in Section~\ref{subsec:tar}; see
\cite{bk:pp} and \cite[Sec.~3.1]{m:yc} for
the calculations. The remaining cases of \eqref{kixpj} and \eqref{xpixpj} are dealt
with in a way quite similar to type $A$, so we will only point out
a few necessary modifications specific for types $B$, $C$ and $D$.

\paragraph{Type $B_n$.}
To verify the $i=j=n$ case of \eqref{kixpj} for $\xi^{-}_n(v)$ write
\ben
\big[h_{n}(u)^{-1}h_{n+1}(u),e_n(v)\big]=
h_n(u)^{-1}\big[h_{n+1}(u),e_n(v)\big]+\big[h_n(u)^{-1},e_n(v)\big]h_{n+1}(u).
\een
Calculating the commutators by
\eqref{Bhn+1en} and \eqref{CompleteBhiej}, we get
\begin{multline}
\non
(u-v)\big[h_{n}(u)^{-1}h_{n+1}(u),e_n(v)\big]
=\frac{1}{2}h_{n}(u)^{-1}h_{n+1}(u)\big(e_n(u)-e_n(v)\big)\\[0.4em]
+\big(e_n(u)-e_n(v)\big)h_{n}(u)^{-1}h_{n+1}(u)
-\frac{u-v}{2(u-v-1)}h_{n}(u)^{-1}\big(e_{n}(u-1)-e_n(v)\big)h_{n+1}(u).
\end{multline}
Relation \eqref{CompleteBhiej} gives
\ben
h_n(u)^{-1}e_{n}(u-1)=e_{n}(u)h_{n}(u)^{-1}
\een
and
\ben
h_n(u)^{-1}e_{n}(v)=\frac{1}{u-v}\ts e_{n}(u)h_{n}(u)^{-1}+\frac{u-v-1}{u-v}\ts e_n(v)h_n(u)^{-1}.
\een
Therefore, we derive
\ben
(u-v)[h_{n}(u)^{-1}h_{n+1}(u),e_n(v)]=\frac{1}{2}\ts\big\{h_{n}(u)^{-1}h_{n+1}(u),e_n(u)-e_n(v)\big\}.
\een
This yields \eqref{kixpj} with $i=j=n$ by writing the relation in terms of
the series $\kappa_n(u)$ and $\xi_n^{-}(v)$. The remaining cases of \eqref{kixpj}
follow by similar arguments.

Now choose the minus signs in \eqref{xpixpj} and verify it for
$i=n-1$ and $j=n$.
By \eqref{Ben-1en} we have
\begin{align}
\label{en-1en}
(u-v+1/2)\big[e_{n-1}(u),e_{n}(v-1/2)\big]
&=e_{n-1\tss n+1}(v-1/2)-e_{n-1\tss n+1}(u)\\[0.4em]
{}&-e_{n-1}(v-1/2)e_n(v-1/2)+e_{n-1}(u)e_n(v-1/2).
\non
\end{align}
By swapping $u$ and $v$ we also get
\begin{align}
\label{enen-1}
(u-v-1/2)\big[e_{n}(u-1/2),e_{n-1}(v)\big]
&=e_{n-1\tss n+1}(u-1/2)e_{n-1\tss n+1}(v)\\[0.4em]
{}&-e_{n-1}(u-1/2)e_n(u-1/2)+e_{n-1}(v)e_n(u-1/2)
\non
\end{align}
which we can also write in the form
\begin{align}
\label{en-1en2}
e_{n-1}(v)&e_{n}(u-1/2)=\frac{u-v-1/2}{u-v+1/2}\ts e_{n}(u-1/2)e_{n-1}(v)\\
&-\frac{1}{u-v+1/2}\Big(e_{n-1\ts n+1}(u-1/2)-e_{n-1\ts n+1}(v)-e_{n-1}(u-1/2)e_n(u-1/2)\Big).
\non
\end{align}
Now \eqref{enen-1} and \eqref{en-1en2} give
\begin{align}
\label{enen-12}
(u-v+1/2)\big[e_{n}(u-1/2),e_{n-1}(v)\big]&=
e_{n-1\tss n+1}(u-1/2)-e_{n-1\tss n+1}(v)\\[0.4em]
{}&-e_{n-1}(u-1/2)e_n(u-1/2)+e_n(u-1/2)e_{n-1}(v).
\non
\end{align}
By \eqref{en-1en} and \eqref{enen-12},
we have
\ben
\bal
&(u-v+1/2)\Big(\big[e_{n-1}(u),e_{n}(v-1/2)\big]+\big[e_{n}(u-1/2),e_{n-1}(v)\big]\Big)\\
&{}=e_{n-1\tss n+1}(u-1/2)-e_{n-1\tss n+1}(v)
-e_{n-1}(u-1/2)e_n(u-1/2)+e_n(u-1/2)e_{n-1}(v)\\[0.3em]
&+e_{n-1\tss n+1}(v-1/2)-e_{n-1\tss n+1}(u)-e_{n-1}(v-1/2)e_n(v-1/2)+e_{n-1}(u)e_n(v-1/2).\\
\eal
\een
Setting
$v=u$ in \eqref{en-1en2} we get
\ben
e_{n-1}(u-1/2)e_n(u-1/2)=\frac{1}{2}\ts\big\{e_{n-1}(u),e_{n}(u-1/2)\big\}
+e_{n-1\ts n+1}(u-1/2)-e_{n-1\ts n+1}(u).
\een
Using also this relation with $u$ replaced by $v$, we thus come to
\ben
\bal
&(u-v+1/2)\Big(\big[e_{n-1}(u),e_{n}(v-1/2)\big]+\big[e_{n}(u-1/2),e_{n-1}(v)\big]\Big)\\
&{}=e_n(u-1/2)e_{n-1}(v)+e_{n-1}(u)e_n(v-1/2)\\[0.3em]
&{}-\frac{1}{2}\ts\big\{e_{n-1}(u),e_{n}(u-1/2)\big\}
-\frac{1}{2}\ts\big\{e_{n-1}(v),e_{n}(v-1/2)\big\}
\eal
\een
which is equivalent to \eqref{xpixpj} for $i=n-1$ and $j=n$.
The case with the plus signs and all other remaining cases follow by similar
calculations.

\paragraph{Type $C_n$.}
Since $[h_i(u),e_n(v)]=0$ and $[h_i(u),f_n(v)]=0$
for $i\leqslant n-1$ by \eqref{CompleteBhiej} and \eqref{CompleteBhifj},
relation \eqref{kixpj} holds for $i<n-1$ and $j=n$. Furthermore,
\eqref{CompleteBhiej} gives
\beql{Chnen}
(u-v+1)\big[h_n(u),e_n(v-1)\big]=-2h_n(u)\big(e_n(u)-e_{n}(v-1)\big)
\eeq
which implies
\ben
\bal
(u-v+1)\big[h_{n-1}(u)^{-1}h_{n}(u),e_{n}(v-1)\big]&
=(u-v+1)h_{n-1}(u)^{-1}\big[h_{n}(u),e_{n}(v-1)\big]\\[0.4em]
&=-2h_{n-1}(u)^{-1}h_{n}(u)\big(e_n(u)-e_{n}(v-1)\big).
\eal
\een
Taking $v=u$ in \eqref{Chnen} we get
\ben
-2\tss h_{n-1}(u)^{-1}h_{n}(u)e_n(u)=-\big\{h_{n-1}(u)^{-1}h_{n}(u), e_{n}(u-1)\big\}.
\een
Hence,
\ben
\bal
(u-v+1)\big[h_{n-1}(u)^{-1}h_{n}(u),e_{n}(v-1)\big]
=(u-v+1)h_{n-1}(u)^{-1}\big[h_{n}(u),e_{n}(v-1)\big]&\\[0.4em]
{}=2h_{n-1}(u)^{-1}h_{n}(u)e_{n}(v-1)-\big\{h_{n-1}(u)^{-1}h_{n}(u), e_{n}(u-1)\big\}&.
\eal
\een
This gives
\ben
(u-v)\big[h_{n-1}(u)^{-1}h_{n}(u),e_{n}(v-1)\big]=
-\big\{h_{n-1}(u)^{-1}h_{n}(u), e_{n}(u-1)-e_{n}(v-1)\big\}
\een
which implies \eqref{kixpj} for $i=n-1$ and $j=n$ for the series $\xi^{-}_n(v)$.
The remaining cases of \eqref{kixpj} follow by similar arguments. In particular,
the case $i=j=n$ for the same series is straightforward from
\eqref{CompleteBhiej} and \eqref{Chn+1en},
whereas the case $i=n$ and $j=n-1$ is derived from
\eqref{CompleteBhiej} and \eqref{completehn+1en-1c}.

Relations \eqref{eiei} and \eqref{eiej=0} imply the respective cases of \eqref{xpixpj}.
The derivation of \eqref{xpixpj} for
$i=n-1$ and $j=n$ relies on
\eqref{Cen-1en} and \eqref{Cfn-1fn}
and follows the same pattern as for type $B_n$ above.

\paragraph{Type $D_n$.}
Applying \eqref{CompleteBhiej} and \eqref{Chn+1en} we get
\ben
\bal
\big[h_{n-1}(u)^{-1}h_{n+1}(u),e_n(v)\big]
&=\frac{1}{u-v}h_{n-1}(u)^{-1}h_{n+1}(u)\big(e_n(u)-e_n(v)\big)\\
&+\frac{1}{u-v}\big(e_n(u)-e_n(v)\big)h_{n-1}(u)^{-1}h_{n+1}(u)
\eal
\een
so that \eqref{kixpj} holds for $i=j=n$ for the series $\xi^{-}_n(v)$.
The remaining cases of \eqref{kixpj} follow with the use of \eqref{CompleteBhiej}
by similar calculations. Relations \eqref{eiei}, \eqref{fifi} and \eqref{eiej=0}
imply the corresponding cases of \eqref{xpixpj}. The remaining case of \eqref{xpixpj}
requiring a longer calculation is $i=n-2$ and $j=n$.
It is performed with the use of \eqref{Den-2en}
and \eqref{Dfn-2fn}
in the same way as for the case $i=n-1$ and $j=n$
in type $B_n$ above.
\epf

Propositions~\ref{prop:YgN} and \ref{pro:relinYgN} imply that the mapping
$\Y^D(\g_{N})\to\Y(\g_N)$ considered in the Main Theorem
is a surjective homomorphism. Its injectivity follows from the decomposition
\beql{decom}
\Y(\g_N)= \Ec\otimes (\Y(\g_N)\cap\Hc)\otimes \Fc
\eeq
and the corresponding arguments of the proof of Theorem \ref{thm:dp}.
This completes the proof of the Main Theorem.

\centerline{\bf Acknowledgments}

We thank the support of South China University of Technology and State Administration of Foreign Experts Affairs during the work.
%the State High-end Project for Foreign Experts during the work.
Jing acknowledges the National Natural Science Foundation of China grant no. 11531004 and Simons Foundation grant no.
523868, Liu acknowledges the National Natural Science Foundation of China grant nos. 11531004 and 11701182, and
Molev acknowledges the support of Australian Research Council.

%\newpage
\bigskip\bigskip

\small

\noindent
N.J.:\qquad\qquad\qquad\qquad\\
Department of Mathematics\\
North Carolina State University, Raleigh, NC 27695, USA\\
jing@math.ncsu.edu

\vspace{3 mm}

\noindent
N.J. \& M.L.:\newline
School of Mathematical Sciences\\
South China University of Technology\\
Guangzhou, Guangdong 510640, China\\
mamliu@scut.edu.cn

\vspace{3 mm}

\noindent
A.M.:\newline
School of Mathematics and Statistics\newline
University of Sydney,
NSW 2006, Australia\newline
alexander.molev@sydney.edu.au

\end{document}